\documentclass[10pt,a4paper]{article}
\usepackage{graphicx}
\usepackage{amsmath, color, verbatim}
\usepackage[colorlinks=true]{hyperref}

\usepackage{amsfonts, mathtools}
\usepackage{titlesec}
\usepackage[hyperref, dvipsnames]{xcolor}
\hypersetup{
colorlinks=true,}

\usepackage{amssymb}
\usepackage{amsthm}
\usepackage{thmtools, thm-restate}

\renewcommand{\d}{\mathrm{d}}

\newcommand{\R}{\mathbb{R}}

\makeatletter
\newcommand*{\T}{%
  {\mathpalette\@T{}}%
}
\newcommand*{\@T}[2]{%
  \raisebox{\depth}{$\m@th#1\intercal$}%
}
\makeatother

\titleformat{\part}
  {\centering\Large\scshape}{}{1em}{Part \thepart: }

\newcommand{\N}{\mathbb N}

\allowdisplaybreaks

  \numberwithin{figure}{section}
\numberwithin{table}{section}

\usepackage{amsmath}
\usepackage{graphicx}
\usepackage{subcaption}
\usepackage[title,titletoc]{appendix}
\usepackage{lipsum}
\usepackage{epstopdf}
\usepackage{mathtools}
\usepackage{amssymb}
\usepackage{xcolor}
\usepackage{hyperref}
\usepackage{latexsym}
\usepackage[export] {adjustbox}
\usepackage{titlesec}
\usepackage[shortlabels]{enumitem}

\usepackage[T1]{fontenc}
\usepackage[utf8]{inputenc}
\usepackage{amsmath,stackengine}
\usepackage{comment}
\usepackage{amsmath,amssymb}
\selectfont
\numberwithin{equation}{section}
\usepackage[square,numbers,sectionbib]{natbib}
\usepackage{varioref}
\usepackage{listings}
\usepackage{color}
\definecolor{dkgreen}{rgb}{0,0.6,0}
\definecolor{gray}{rgb}{0.5,0.5,0.5}
\definecolor{mauve}{rgb}{0.58,0,0.82}

\usepackage{tikz}
\usetikzlibrary{shapes.geometric, arrows}
\tikzstyle{startstop} = [rectangle, rounded corners, minimum width=3cm, minimum height=1cm,text centered, draw=black, fill=red!30]
\tikzstyle{io} = [trapezium, trapezium left angle=70, trapezium right angle=110, minimum width=3cm, minimum height=1cm, text centered, draw=black, fill=blue!30]
\tikzstyle{process} = [rectangle, minimum width=3cm, minimum height=1cm, text centered, draw=black, fill=orange!30, text width=9.5cm]
\tikzstyle{process2} = [rectangle, minimum width=3cm, minimum height=1cm, text centered, draw=black, fill=orange!30]
\tikzstyle{decision} = [diamond, minimum width=3cm, minimum height=1cm, text centered, draw=black, fill=green!30]
\tikzstyle{arrow} = [thick,->,>=stealth]

\usepackage{empheq,etoolbox}
\patchcmd{\subequations}
{\theparentequation\alph{equation}}
{\theparentequation.\alph{equation}}
{}{}

\newcommand{\numberset}{\mathbb}
\renewcommand{\N}{\numberset{N}}
\renewcommand{\R}{\numberset{R}}
\newcommand{\vphi}{\boldsymbol{\varphi}}
\newcommand{\veta}{\boldsymbol{\eta}}
\newcommand{\vnu}{\boldsymbol{\nu}}
\newcommand{\vlambda}{\boldsymbol{\Lambda}}
\usepackage{amsthm}

\theoremstyle{plain}
\newtheorem{thm}{Theorem}[section]
\newtheorem{cor}[thm]{Corollary}
\newtheorem{lem}[thm]{Lemma}

\newtheorem{prob}[thm]{Problem}

\newcommand{\eps}{h}

\theoremstyle{definition}

\theoremstyle{remark}
\newtheorem{oss}[thm]{Remark}
\everymath{\displaystyle}

\newcommand{\vertiii}[1]{{\left\vert\kern-0.25ex\left\vert\kern-0.25ex\left\vert #1
		\right\vert\kern-0.25ex\right\vert\kern-0.25ex\right\vert}}
\newcommand{\norm}[1]{{\left\Vert #1
		\right\Vert}}	

\def\Xint#1{\mathchoice
{\XXint\displaystyle\textstyle{#1}}%
{\XXint\textstyle\scriptstyle{#1}}%
{\XXint\scriptstyle\scriptscriptstyle{#1}}%
{\XXint\scriptscriptstyle\scriptscriptstyle{#1}}%
\!\int}
\def\XXint#1#2#3{{\setbox0=\hbox{$#1{#2#3}{\int}$ }
\vcenter{\hbox{$#2#3$ }}\kern-.6\wd0}}

\def\dashint{\Xint-}

\usepackage{tabulary}
\usepackage{booktabs}
\usepackage{color}
\usepackage{bbm}
\usepackage{multicol}
\usepackage{framed}

\newcolumntype{2}{D{.}{}{2.0}}

\definecolor{light-gray}{gray}{0.9}

\renewenvironment{leftbar}[1][\hsize]
{%
    \MakeFramed{\hsize#1\advance\hsize-\width\FrameRestore}%
}
{\endMakeFramed}

\def\pphi{\boldsymbol\varphi}
\def\d{{\rm d}}

\def \HH{\mathcal{H}}
\def\VV{\mathcal{V}}

\def\dx{\ dx}

\def\tV{\widetilde{V}}
\def\peta{\boldsymbol{\eta}}
\def\pxi{\boldsymbol\xi}

\def\multibold #1{\def\arg{#1}%
	\ifx\arg\pto \let\next\relax
	\else
	\def\next{\expandafter
		\def\csname #1#1#1\endcsname{{\boldsymbol #1}}%
		\multibold}%
	\fi \next}

\def\pto{.}

\def\multical #1{\def\arg{#1}%
	\ifx\arg\pto \let\next\relax
	\else
	\def\next{\expandafter
		\def\csname #1#1\endcsname{{\cal #1}}%
		\multical}%
	\fi \next}

\def\multimathop #1 {\def\arg{#1}%
	\ifx\arg\pto \let\next\relax
	\else
	\def\next{\expandafter
		\def\csname #1\endcsname{\mathop{\rm #1}\nolimits}%
		\multimathop}%
	\fi \next}

\multibold
qwertyuiopasdfghjklzxcvbnmQWERTYUIOPASDFGHJKLZXCVBNM.

\multical
QWERTYUIOPASDFGHJKLZXCVBNM.

\multimathop
diag dist div dom mean meas sign supp .

\def\numberset{\mathbb}
\def\N{\numberset{N}}
\def\R{\numberset{R}}

\newcommand{\al}{\pmb\alpha}
\def\ww{\www}

\DeclareMathOperator*{\argmin}{arg\,min}
\everymath{\displaystyle}

\def \no#1#2#3 {{\bf #1} (#3), #2.}
\def \eds#1#2#3 {#1, #2, #3.}
\def\an #1{#1}
\renewcommand{\T}{\mathrm T}
\renewcommand{\d}{\mathrm d}

\newcommand{\tgrad}{\nabla_{\Gamma}}

\def \mmu{\ww}
\def\ol{\omega(\vphi_0,{u}_0)}



\def\Vb{\pmb{\VV}}
\def\Pb{\pmb{\PP}}
\def\mmm{\al}
\def\Andrea#1{{\color{black}#1}}

\begin{document}
\hypersetup{
  urlcolor     = blue,
  linkcolor    = Bittersweet,
  citecolor   = Cerulean
}
\title{Multi-component phase separation and small deformations \\ of a spherical biomembrane}
\author{Diogo Caetano\thanks{Mathematics Institute, University of Warwick, Coventry CV4 7AL, UK ({\tt Diogo.Caetano@warwick.ac.uk})}
\and Charles M. Elliott\thanks{Mathematics Institute, University of Warwick, Coventry CV4 7AL, UK ({\tt  C.M.Elliott@warwick.ac.uk})}
\and Maurizio Grasselli\thanks{Dipartimento di Matematica, Politecnico di Milano, Milano, 20133, Italy ({\tt maurizio.grasselli@polimi.it})}
\and Andrea Poiatti\thanks{Faculty of Mathematics, University of Vienna, Vienna, 1090, Austria
({\tt andrea.poiatti@univie.ac.at})}}
\date{}

\maketitle

\begin{abstract}\noindent
We focus on the derivation and analysis of a model for multi-component phase separation occurring on biological membranes, inspired by observations of lipid raft formation. The model integrates local membrane composition with local membrane curvature, describing the membrane's geometry through a perturbation method represented as a graph over an undeformed Helfrich minimising surface, such as a sphere. The resulting energy consists of  a small deformation functional coupled to a Cahn-Hilliard functional. By applying Onsager's variational principle, we obtain a multi-component Cahn-Hilliard equation for the vector $\vphi$ of protein concentrations coupled to an evolution equation
for the small deformation $u$ along the normal direction to the reference membrane. Then, in the case of a constant mobility matrix, we consider the Cauchy problem and
we prove that it is (globally) well posed in a weak setting. We also demonstrate that any weak solution regularises in finite time and satisfies the so-called ``strict separation property''. This property allows us
to show that any weak solution converges to a single stationary state by a suitable version of the {\L}ojasiewicz-Simon inequality.



\end{abstract}

\section{Introduction}
\subsection{The model}
The Canham-Helfrich energy is the  elastic bending energy for a biomembrane, modelled as a thin deformable surface. It is given by
(see, e.g., \cite{ZH1989} and references therein)
\begin{align*}
    \mathrm E(\Gamma) = \int_\Gamma \dfrac{\kappa}{2} (H-H_s)^2 + \sigma + \kappa_G K \, \d \Gamma.
\end{align*}
 In the above, $\Gamma$ is a 2-dimensional closed hypersurface in $\R^3$ enclosing an open domain; $\kappa>0$ and $\kappa_G$ are bending rigidities, and $\sigma\geq 0$ is the surface tension; $H$ and $K$ are, respectively, the mean curvature and Gauss curvature of $\Gamma$; $H_s$ is called the spontaneous curvature of $\Gamma$ and is a measure of stress within the membrane for the flat configuration. We are interested in modelling biomembranes that consist of multiple differing lipid types undergoing phase separation (see, for instance, \cite{HF2011,WBV2012} and their references). The importance of analysing multi-component biological systems is stressed, for instance, in \cite{ZL2022} (see also \cite{SB2021}).
Thus, motivated by the work in \cite{EllHat21} in the binary case, we introduce an order parameter $\vphi = (\varphi_1, \dots, \varphi_N)$ consisting of concentrations of each lipid with $0\le\vphi_i\le 1,~\sum_{i-1}^N\vphi_i=1$ and assume that the energy takes the form
\begin{align}
\label{eq:energy}
    \mathrm E(\Gamma, \vphi) = \int_\Gamma \dfrac{\kappa(\vphi)}{2} (H - H_s(\vphi))^2 + \sigma + \kappa_G(\vphi) K \, \d\Gamma+ b\int_\Gamma \dfrac{\varepsilon}{2}|\nabla_\Gamma \vphi|^2 + \dfrac{1}{\varepsilon} \Psi(\vphi) \, \d\Gamma,
\end{align}
where the free energy density $\Psi$ is defined as
\begin{align}\label{5_eq:logpot}
    \Psi(\vphi) = \sum_{i=1}^N \vphi_i \log \vphi_i - \dfrac{1}{2} \vphi \cdot  A \vphi=\mathbf\Psi^1(\vphi)+\mathbf\Psi^2(\vphi),
\end{align}
with $ A\in \R^{n\times n}$ symmetric and possesses a positive eigenvalue. Thus, it is assumed that the material parameters depend on the phase field variables $\varphi_i\in [0,1]$, each corresponding to the density of the $i$th component of the system. The second functional is a Ginsburg-Landau energy with coefficient $b>0$ which accounts for phase separation. We recall that $\Psi$ is the so-called Flory-Huggins potential density (see, e.g., \cite{RB2024} for its role in Cell Biology).

Our goal is to formulate and analyse a coupled system of equations for small deformations of a spherical surface coupled with a multicomponent phase field system. To derive the system we first  approximate the energy \eqref{eq:energy} in the case that the spontaneous curvature is scaled with a small parameter $\rho$ and the gradient energy  parameter $b$ scaled  with $\rho^2$ as  small perturbations of  a spherical surface that minimises the energy in the case of prescribed enclosed volume,  constant elastic energy parameters and  $\rho=0$.
Then, \Andrea{assuming for simplicity that $\kappa$ and $\kappa_G$ do not depend on $\vphi$ and using Onsager's variational principle as a tool to formulate dissipative dynamics,} we  obtain   the following  evolution system for state  variables
$(\vphi, u)$ holding on a fixed sphere $\Gamma$  
arising as a   critical point of a Rayleighian $\mathcal R$:
\begin{align}\label{MAIN}
\begin{split}
    \partial_t \vphi &= \nabla_\Gamma \cdot \left( \mathbf L(\vphi) \nabla_\Gamma \boldsymbol\mu\right) \\
    \boldsymbol\mu &= \dfrac b\varepsilon \left(\boldsymbol\psi'(\vphi) - A\vphi\right) + \dfrac{2\kappa u \vlambda}{R^2} - b\varepsilon \Delta_\Gamma \vphi + \kappa (\vlambda \cdot \vphi) \vlambda + \kappa \Delta_\Gamma u \, \vlambda \\
    \beta \partial_t u &=  -\kappa \Delta_\Gamma^2 u + \left(\sigma - \dfrac{2\kappa}{R^2} \right) \Delta_\Gamma u + \dfrac{2\sigma u}{R^2} - \kappa \vlambda \cdot \Delta_\Gamma \vphi - \dfrac{2\kappa \vlambda\cdot(\vphi-\pmb \alpha)}{R^2},
    \end{split}
\end{align}
\Andrea{where $\al=\dashint_{\Gamma} \vphi :=\frac{\int_{\Gamma} \vphi\d\Gamma}{\vert \Gamma\vert}$}, or equivalently
\begin{align}\label{5_eq:alternative}
\begin{split}
    \partial_t \vphi &= \nabla_\Gamma \cdot \left( \mathbf L \nabla_\Gamma \boldsymbol\mu\right) \\
    \boldsymbol\mu &= \dfrac b\varepsilon \left(\boldsymbol\psi'(\vphi) - {A}\vphi\right) + \dfrac{2\kappa u \vlambda}{R^2} - b\varepsilon \Delta_\Gamma \vphi + \kappa (\vlambda \cdot \vphi) \vlambda + \kappa \Delta_\Gamma u \, \vlambda \\
    \beta \partial_t u &=  \left(\kappa \Delta_\Gamma + \dfrac{2\kappa }{R^2} \right) \left( -\Delta_\Gamma u + \dfrac{\sigma}{\kappa} u - \boldsymbol\Lambda\cdot (\vphi-\pmb\alpha) \right),
    \end{split}
\end{align}
on $\Gamma \times (0,T)$. Here $\mathbf L(\cdot)$ is a mobility matrix \Andrea{(which will be assumed constant in the mathematical analysis of the problem)}, $\boldsymbol\mu$ is a chemical potential, $\vlambda$ is a spontaneous curvature parameter, and $R$ is the radius of the sphere. The notation is explained in Section \ref{Derive}.

\subsection{Contribution and outline of paper}
Our first contribution is to derive a model for multi-component  phase separation on a spherical biomembrane that allows for small deformations.
In the case of a constant mobility matrix, we first show that the associated Cauchy problem has a unique (globally defined) weak solution which continuously depends on the initial datum. Then, thirdly we prove that any
weak solution regularises in finite time and each $\vphi_i$ stays instantaneously and uniformly away from pure phases (i.e. $0$ and $1$). The latter result is known as ``strict separation property'' and plays an important role in the analysis of the properties of solutions to Cahn-Hilliard and Allen-Cahn type equations with singular potentials, e.g., of logarithmic type (see \cite{GalPoia} and references therein,
see also \cite{GGPS,AC2023} for multi-component systems). In particular, it is essential when one wants to show the convergence of a solution to a single stationary state. Finally, we prove this as a further result.

We recall that the original scalar model (see \cite{EllHat21}) was formulated by using a conserved Allen-Cahn equation. On account of \cite{AC2023}, we think that results similar to the ones obtained here can be proven for a multi-component version of that model as well.
It is also worth observing that there have been very few contributions to the theoretical analysis of multi-component Cahn-Hilliard equations with singular potential since the pioneering work \cite{EllLuc91}  (see \cite{EllGar97-b, GGPS, Gar05} and references therein, see also \cite{AGP2024} for multi-component fluids). Concerning the numerical analysis, the reader is referred, for instance, to \cite{LHXF2024}. {\color{black}In conclusion, we observe that recently there has been an increasing interest in the mathematical analysis of phase separation phenomena on evolving surfaces with prescribed surface evolution (see \cite{AGP1,AGP2,CE, CEGP, ES} and references therein). The model derived in the present contribution is an alternative approach to the description of the same biological phenomena, but incorporating the evolution of the surface in the equations, without prescribing it \textit{a priori} (although under the small deformation assumption).}

The plan of the paper is as follows. Section \ref{Derive} is devoted to the derivation of our approximating energy and the associated  evolution equations. The following section  introduces some useful notation as well as the main analytic assumptions. The statements of the main results are contained in Section \ref{mainres}.
Their proofs are contained in Sections \ref{proof1}-\ref{convec}.  Finally, an Appendix is devoted to prove Theorem \ref{5_thm:energy} and a pair of technical lemmas.

\section{Derivation of the model}\label{Derive}


\subsection{The approximate energy}
We begin by deriving an approximate energy, as in \cite{EllFriHob17-a,EllHat21},
making the following assumptions on the model:
\begin{itemize}
    \item[(A1)] the only material parameter that depends on the phase field $\vphi$ is the spontaneous curvature $H_s(\vphi)$, and it is of the form
    \begin{align*}
        H_s(\vphi) := \vlambda \cdot \vphi, \quad \vlambda = (\Lambda_1,\dots,\Lambda_N)\in \R^N;
    \end{align*}
    in particular, $\kappa(\vphi)\equiv \kappa$ and $\kappa_G(\vphi)\equiv \kappa_G$ are positive constants;
    \item[(A2)] we rescale the coefficients: $\vlambda$ is replaced by $\rho\vlambda$ and $b$ is replaced by $\rho^2 b$, where $0<\rho \ll 1$ is a small parameter;
    \item[(A3)] the volume enclosed by $\Gamma$ and the mass of each lipid type are preserved.
\end{itemize}
This allows us to rewrite \eqref{eq:energy} as
\begin{align*}
  \mathrm E_\rho(\Gamma, \vphi) = \int_\Gamma \dfrac{\kappa}{2} (H - \rho \vlambda\cdot\vphi)^2 + \sigma + \kappa_G K \, \d\Gamma+ \rho^2 b \int_\Gamma \dfrac{\varepsilon}{2}|\nabla_\Gamma \vphi|^2 + \dfrac{\Psi(\vphi)}{\varepsilon} \, \d\Gamma.
\end{align*}
   Since $\kappa_G$ is now independent of the phase field, and because $\Gamma$ is a closed surface, the Gauss-Bonnet theorem guarantees that
\begin{align*}
    \kappa_G \int_\Gamma K = 2\pi\kappa_G \chi(\Gamma)=c,
\end{align*}
where $\chi(\Gamma)$ denotes the Euler characteristic of $\Gamma$. Hence the term involving the Gauss curvature is simply a constant and, abusing notation, we drop it from the energy. Expanding also the square in the first term leads to
\begin{align}\label{eq:energy2}
    \mathrm E_\rho(\Gamma, \vphi) = \int_\Gamma (\dfrac{\kappa}{2} H^2 - \kappa \rho (\vlambda\cdot \vphi)H + \dfrac{\kappa \rho^2}{2} (\vlambda\cdot\vphi)^2 + \sigma )\, \d\Gamma +\rho^2 b \int_\Gamma \dfrac{\varepsilon}{2}|\nabla_\Gamma \vphi|^2 + \dfrac{\Psi(\vphi)}{\varepsilon}\, \d\Gamma.
\end{align}
We now obtain an energy $\mathcal E$ which approximates \eqref{eq:energy2} using the perturbation method as in \cite{EllHat21}.
By consideration of forcing terms  $\mathcal F_1=\mathcal F_1(\Gamma, \vphi)$ and $\mathcal F_2=\mathcal F_2(\Gamma, \vphi)$ defined by
\begin{align*}
    \mathcal F_1(\Gamma, \vphi) = -\int_\Gamma H \vlambda \cdot \vphi \d\Gamma, \quad \mathcal F_2(\Gamma, \vphi) = \int_\Gamma \dfrac {b\varepsilon} 2 |\nabla_\Gamma \vphi|^2 + \dfrac{b}{\varepsilon}\Psi(\vphi) + \dfrac{\kappa (\vlambda\cdot \vphi)^2}{2}\d\Gamma
\end{align*}
we rewrite  the free energy given by \eqref{eq:energy2} as
$$ \mathrm E_\rho(\Gamma, \vphi) = \kappa \mathcal W(\Gamma) + \sigma \mathcal A(\Gamma)+\rho \kappa \mathcal F_1(\Gamma, \vphi) + \rho^2 \mathcal F_2(\Gamma, \vphi)$$
where we define the functionals
 \begin{align*}
    \mathcal W(\Gamma) = \dfrac 12 \int_\Gamma H^2\d\Gamma, \quad \mathcal A(\Gamma) = \int_\Gamma 1 \d\Gamma, \quad \mathcal V(\Gamma) = \dfrac{1}{3}\int_\Gamma \text{Id}_\Gamma \cdot \nu \d\Gamma.
\end{align*}
Here   $\mathcal V(\Gamma)$ is the volume enclosed by $\Gamma$  and $\vnu$ denotes the unit outward normal to $\Gamma$.

Now, let us  consider the Lagrangian
\begin{align}\label{5_eq:L}
    \mathcal L(\Gamma, \lambda) = \kappa \mathcal W(\Gamma) + \sigma \mathcal A(\Gamma) + \lambda (\mathcal V(\Gamma) - V_0)
\end{align}

 \noindent in which the volume  constraint   $ \mathcal V(\Gamma):=V_0=\frac{4}{3}R^3>0$ is incorporated by introduction of  the Lagrange multiplier $\lambda\in \R$.
 It was shown in \cite{EllFriHob17-a} that $\mathcal L$ in \eqref{5_eq:L} has a critical point $(\Gamma_0, \lambda_0)$ defined by
\begin{align}\label{5_eq:criticals}
    \Gamma_0 = \text{ sphere of radius } R \text{ centered at the origin} \quad \text{ and } \quad \lambda_0 = -\dfrac {2\sigma}R.
\end{align}
We also set $\nu_0$, $H_0$ as to be the outer unit normal and the curvature of $\Gamma_0$, respectively.
Minimisers of \eqref{eq:energy2} are the critical points of the perturbed functional
\begin{align}\label{5_eq:Lrho}
    \mathcal L_\rho(\Gamma, \lambda; \vphi) = \mathcal L(\Gamma, \lambda) + \rho \kappa \mathcal F_1(\Gamma, \vphi) + \rho^2 \mathcal F_2(\Gamma, \vphi).
\end{align}
The evolution problem studied in this paper is obtained from  a second order approximation to $\mathcal L_\rho$.  More precisely, setting  $\vphi$ in the neighbourhood of $\Gamma_0$ by extending  constantly in the normal direction
\begin{align*}
    \vphi \left(x + \rho u(x) \vnu(x)\right) = \vphi(x), \quad  x\in\Gamma_0.
\end{align*}
we have (see Appendix for the proof)

\begin{thm}\label{5_thm:energy}
Let $\vphi\in C^2(\Gamma_0)$, $\lambda_1\in \R$ and $u\in C^2(\Gamma_0)$ satisfy $\dashint_{\Gamma_0}u=0$. Then
\begin{align}\label{5_eq:taylor_expansion}
    \mathcal L_\rho(\Gamma_0 + \rho u \nu_0, \lambda_0 + \rho\lambda_1; \vphi) = C_1 + \rho C_2 + \rho^2 \mathcal E(\vphi, u) + \mathcal O(\rho^3),
\end{align}
for some constants $C_i=C_i(\kappa, \Gamma_0, \sigma, \mathbf \Lambda, \alpha)$, $i=1,2$, and
$$\mathcal E(u,\vphi) = \mathrm E_{\mathrm H} (u,\vphi) + \mathrm E_{\mathrm {CH}}(\vphi),$$ where
\begin{align}\label{5_eq:energy1}
    \mathrm E_{\mathrm H}(u,\vphi) &= \int_{\Gamma_0} \dfrac{\kappa  (\Delta_\Gamma u)^2}{2} + \left(\sigma - \dfrac{2\kappa}{R^2}\right)\dfrac{|\tgrad  u|^2}{2}- \dfrac{\sigma u^2}{R^2}+ \kappa (\vlambda \cdot \vphi) \Delta_\Gamma u + \dfrac{2\kappa u \vlambda \cdot \vphi}{R^2} +\dfrac{\kappa (\vlambda\cdot\vphi)^2}{2} \d\Gamma
\end{align} and
\begin{align}\label{5_eq:energy2}
    \mathrm E_{\mathrm{CH}}(\vphi) &= b\int_{\Gamma_0} \dfrac{\varepsilon}{2}|\tgrad  \vphi|^2 + \dfrac{ \Psi(\vphi)}{\varepsilon}\d\Gamma.
\end{align}
  \end{thm}
\Andrea{
\begin{oss}
    Notice that the condition $\dashint_{\Gamma_0} u=0$ is necessary to guarantee that the volume of $\Gamma_0+\rho u\nu_0$ is the same as the volume of $\Gamma_0$. Indeed (see Lemma \ref{5_lem:variations_standard}), it holds $\mathcal V'(\Gamma_0) [\rho u\nu_0] = \rho\int_{\Gamma_0} u\d\Gamma$, and we need this variation to be zero.
\end{oss}}

 We will use mass conservation and variational considerations   in order to obtain the  system of evolution equations, i.e., (\ref{MAIN}).

\subsection{Conservation of lipid concentrations}
For the sake of simplicity, from now on we use  the generic notation $\Gamma$ in place of  the sphere of radius $R$, $\Gamma_0$. Our model is that we have $N$  lipid components with concentration  $\vphi_i$ such that
 each component  is conserved  and evolves via a flux. The  continuity equation for the  order parameter $\vphi$ is  written
\begin{align}\label{5_eq:continuity}
    \partial_t \vphi + \nabla_\Gamma\cdot \mathbf J(\vphi) = 0, \quad \text{ meaning } \,\, \partial_t \vphi_i + \nabla_\Gamma\cdot \mathbf{J}_i(\vphi) = 0, \,\, \text{ for all } i=1,\dots,N,
\end{align}
where $\mathbf J$ is the flux 
such that $$\mathbf J := [\mathbf{J}_1, \dots, \mathbf{J}_N]^T.$$  We assume $\mathbf z := [\mathbf{z}_1,\dots,\mathbf{z}_N]^T$, where each $\mathbf z_i$ is a thermodynamic field tangential to $\Gamma$ representing the dissipation by drag on component $i$, and 
 the constitutive relation
\begin{align*}
     \,  \mathbf{J}_i :=
     \sum_{j=1}^N  L_{ij}\mathbf z_j,
\end{align*}
where the mobility matrix $\mathbf L(\pphi)\equiv\{L_{ij}(\pphi)\} \in \mathbb R^{N\times N}$ is symmetric, positive semidefinite, and satisfies
\begin{align*}
    \sum_{j=1}^N  L_{ij} = 0 \,\,\, \forall i, \quad \text{ and } \quad \mathbf L  {\pmb \eta} \cdot {\pmb \eta} \geq C |
    \veta|^2, \quad \forall {\pmb \eta} \in \{ {\pmb \eta} = (\eta_1, \dots, \eta_N) \colon \sum_{j=1}^N \eta_j = 0\}.
\end{align*}
Allowing  for concentration  dependent mobilities, possibly degenerate, we  write  $\mathbf L = \mathbf L(\vphi)$, so that $\mathbf J = \mathbf J(\vphi)$. By construction, we observe that  $\mathbf J\nu=\mathbf 0$.
We can thus rewrite the continuity equation \eqref{5_eq:continuity} together with the conservation of mass as
\begin{align}\label{ConsMass}\begin{cases}
    \partial_t \vphi_k + \nabla_\Gamma\cdot (\mathbf L (\vphi) \mathbf z)_k = 0,~k=1,2,..N\\
   \dashint_{\Gamma} \vphi =\dashint_{\Gamma} \vphi(0) =: \pmb \alpha .
    \end{cases}
\end{align}
\subsection{First variation of the energy}
For fixed state variables $\left(u,\vphi\right)$, we set $\delta{\mathcal E}(\zeta,\veta)$ to be the first variation of $\mathcal E$ in the direction $(\zeta,\veta)$. A straightforward calculation and an integration by parts yields

\begin{align*}\begin{cases}
     \delta{\mathcal E}(\zeta,\veta)=\delta{\mathcal E}_1(\zeta,\veta)+\delta{\mathcal E}_2(\zeta,\veta)+\delta{\mathcal E}_3(\zeta,\veta)\\
      \delta{\mathcal E}_1(\zeta,\veta):=\int_\Gamma \left( \kappa \Delta_\Gamma^2 u - \left(\sigma - \dfrac{2\kappa}{R^2} \right) \Delta_\Gamma u - \dfrac{2\sigma u}{R^2}\right)\zeta \d\Gamma \\
        \delta{\mathcal E}_2(\zeta,\veta):=\int_\Gamma \left(\kappa \vlambda \cdot \Delta_\Gamma \vphi + \dfrac{2\kappa \vlambda\cdot\vphi}{R^2}  \right) \zeta+ \int_\Gamma\left( \kappa (\vlambda \cdot \vphi) \vlambda +\kappa \Delta_\Gamma u  \, \vlambda + \dfrac{2\kappa u \vlambda}{R^2}\right)\cdot\veta \d\Gamma
  \\ \delta{\mathcal E}_3(\zeta,\veta):=\int_\Gamma \left( -b\varepsilon \Delta_\Gamma \vphi + \dfrac b\varepsilon \left(\boldsymbol\psi'(\vphi) - A\vphi\right)  \right)
    \cdot \veta \d\Gamma.
    \end{cases}
 \end{align*}
Note that, for the reader's convenience, we have separated the contributions arising from  the pure elastic membrane energy, the spontaneous curvature coupling and the Cahn-Hlliard energy. Now it is convenient to introduce a chemical potential
\begin{align}
\boldsymbol\mu:= \dfrac b\varepsilon \left(\boldsymbol\psi'(\vphi) - A\vphi\right) + \dfrac{2\kappa u \vlambda}{R^2} - b\varepsilon \Delta_\Gamma \vphi + \kappa (\vlambda \cdot \vphi) \vlambda + \kappa \Delta_\Gamma u \, \vlambda.
\end{align}
and an elastic restoring force
\begin{align}
\mathbf F:= \left(\kappa \Delta_\Gamma^2 u - \left(\sigma - \dfrac{2\kappa}{R^2} \right) \Delta_\Gamma u - \dfrac{2\sigma u}{R^2} + \kappa \vlambda \cdot \Delta_\Gamma \vphi + \dfrac{2\kappa \vlambda\cdot\vphi}{R^2}\right)\vnu,
\end{align}
denoting by $F$ its modulus, and write
\begin{align}
\label{variation_of_energy_onsager}
 \delta{\mathcal E}(\zeta,\veta):=\int_\Gamma F\zeta \d\Gamma+\int_\Gamma     \boldsymbol\mu \cdot \veta \d\Gamma.
\end{align}
It follows that
\begin{align} \label{D_tEnergy}\dfrac{d}{dt}\mathcal E(u,\vphi)=\int_\Gamma Fu_t \d\Gamma+\int_\Gamma  \boldsymbol\mu \cdot\vphi_t \d\Gamma.\end{align}

\subsection{Dissipative potential, Rayleighian and Onsager's principle}
Here we follow arguments in  \cite{ArrWalTorr18, doi11, ons31a}.  Given the state variables we wish to find evolution equations  for the  dynamical variables, $(V,\mathbf z)$, where $V$ is the normal velocity of the surface and $\mathbf z= [\mathbf{z}_1,\dots,\mathbf{z}_N]^T$ is the collection of drag velocities within  the flux $\mathbf L(\vphi)\mathbf z$. This is achieved using  an  application of Onsager's variational principle by minimising a potential involving dissipative forces. We begin by using (\ref{D_tEnergy}) and respecting the conservation of mass constraint (\ref{ConsMass}), and set for a scalar field $v$ and a tangential vector field $\mathbf v$
\begin{align}
\dot{\mathcal E}(v,\mathbf v)=\int_\Gamma Fv \d\Gamma-\int_\Gamma  \boldsymbol\mu \cdot\nabla_\Gamma( \mathbf L(\vphi) \mathbf v) \d\Gamma.\end{align}
Integrating by parts the second term, we get 
\begin{align*}
     \int_\Gamma \nabla_\Gamma \boldsymbol\mu \colon \mathbf L(\vphi) \mathbf v \d\Gamma
\end{align*}
and this leads to the expression
\begin{align}
\dot{\mathcal E}(v,\mathbf v)=\int_\Gamma Fv \d\Gamma+ \int_\Gamma  \mathbf L(\vphi)\nabla_\Gamma  \boldsymbol\mu \cdot \mathbf v \d\Gamma.\end{align}

The   dissipation potential is chosen to be
\begin{align*}
    \mathcal D(v, \mathbf v) = \mathbf D_1(v) + \mathbf D_2(\mathbf v)    &= \dfrac{\beta}{2} \int_\Gamma \vert v \vert^2 \d\Gamma+\dfrac{1}{2} \int_\Gamma \mathbf v \cdot\mathbf J \d\Gamma,
\end{align*}
consisting of a viscous kinetic energy term $\mathbf D_1$ associated with the normal velocity of the membrane, $v$,  where $\beta\geq 0$, and a   quadratic energy  $\mathbf D_2$,
with protein velocities, $\mathbf v $,   on the membrane. For the state $(u,\vphi)$ the Rayleighian potential $\mathcal R$ is then defined as
\begin{align}
\mathcal R(v, \mathbf v) = \dot{\mathcal E}(v,\mathbf v) + \mathcal D(v, \mathbf v).
\end{align}
Onsager's variational principle yields
\begin{align}
(V,\mathbf z)=\mathop{\argmin}\limits_{v,\mathbf v} \mathcal R(v,\mathbf v)
\end{align}
 subject to the conservation of mass (\ref{ConsMass}) and
 the volume constraint
$\int_\Gamma v \d\Gamma=\dfrac{d}{dt}\int_\Gamma u \d\Gamma = 0$.
Introducing a Lagrange multiplier we are led to finding a critical point $(V,\mathbf w, \lambda ^*)$ of the Lagrangian
\begin{align*}
    \mathcal L(v, \mathbf v, \lambda) := \mathcal R(v,\mathbf v) + \lambda\int_\Gamma v \d\Gamma.\end{align*}
Using the above calculations, we get
\begin{align}
 &\int_\Gamma  \mathbf L(\vphi)\left( \nabla_\Gamma \boldsymbol\mu + \mathbf z \right)\cdot \mathbf v \d\Gamma= 0,
    \\
    &\int_{\Gamma} \Bigg(F + \beta V+\lambda^* \Bigg) v \d\Gamma= 0,\\
    &\int_\Gamma V\d\Gamma=0.
\end{align}
From the above we obtain
$$  \mathbf L(\vphi)\mathbf z=-\mathbf L(\vphi)\nabla_\Gamma  \boldsymbol\mu,~~F + \beta V+\lambda^*=0, ~~ \lambda^* = - \dfrac{2\kappa \vlambda\cdot \pmb \alpha}{R^2}$$
and writing
$$u_t=V, ~~\partial_t\vphi_k+(\nabla \cdot (\mathbf L(\vphi)\mathbf z))_k=0,\quad k=1,\ldots,N.$$
we deduce (\ref{MAIN}).

\section{Notation and main assumptions}

	\label{setting}
	\subsection{Function spaces}
Let $\Gamma$ be a two dimensional
closed, connected, orientable, and sufficiently smooth hypersurface (embedded in $\R^3$).
In particular, in our model we assume that $\Gamma=\Gamma_0$ (see \eqref{5_eq:criticals}).
 We denote the usual Sobolev spaces by $W^{k,p}(\Gamma )$%
	, where $k\in \mathbb{N}$ and $1\leq p\leq +\infty $, with norm $\Vert \cdot
	\Vert _{W^{k,p}(\Gamma )}$. The Hilbert space $W^{k,2}(\Gamma )$ is identified with $H^{k}(\Gamma )$ with norm $\Vert \cdot \Vert _{H^{k}(\Gamma )}$.
	Furthermore, given a (real) vector space $X$, we set $\mathbf{X}=X^3$. We then denote by $(\cdot
	,\cdot )$ the inner product in $L^{2}(\Gamma )$ and by $\Vert \cdot \Vert $
	the corresponding induced norm. By $(\cdot ,\cdot )_{X}$ and $\Vert \cdot
	\Vert _{X}$ we indicate the canonical inner product and its induced norm in a generic real Hilbert
	space $X$, respectively. Moreover, we set $\VVV=\HHH^1(\Gamma)$ and $\HHH=\LLL^2(\Gamma)$, as well as $H=L^2(\Gamma)$ and $V=H^1(\Gamma)$. 
	We then introduce the
	affine hyperplane
	\begin{equation}
	\Sigma :=\left\{ \mathbf{c}^{\prime }\in \mathbb{R}^{N}:%
	\sum_{i=1}^{N}c_{i}^{\prime }=1\right\} ,  \label{sigma}
	\end{equation}%
	the Gibbs simplex
	\begin{equation}
	\mathbf{G}:=\left\{ \mathbf{c}^{\prime }\in \mathbb{R}^{N}:%
	\sum_{i=1}^{N}c_{i}^{\prime }=1,\quad c_{i}^{\prime }\geq 0,\quad i=1,\ldots
	,N\right\} ,  \label{Gibbs}
	\end{equation}%
	and the corresponding tangent space to $\Sigma $
	\begin{equation}
	T\Sigma :=\left\{ \mathbf{d}^{\prime }\in \mathbb{R}^{N}:%
	\sum_{i=1}^{N}d_{i}^{\prime }=0\right\} .  \label{sigma2}
	\end{equation}%
	Then, we set:
	{\color{black}\begin{equation*}
\an{\HHH_0}:=\{\fff\in \HHH:\ \int_{\Gamma
}\fff\\dx =\mathbf{0}\text{ and } \fff(x)\in
T\Sigma \;\text{ for a.a. }x\in \Gamma \},
\end{equation*}%
\begin{equation*}
\an{\widetilde{\HHH}_{0}}:=\{\fff\in \HHH:\fff(x)\in T\Sigma \,\text{ for a.a. }x\in
\Gamma \},
\end{equation*}%
\begin{equation*}
\VVV_{0}:=\{\fff\in \VVV:\int_{\Gamma }%
\fff\\dx =\mathbf{0}\text{ and }\fff(x)\in T\Sigma\, \text{ for a.a. }x\in
\Gamma \},
\end{equation*}%
\begin{equation*}
\widetilde{\VVV}_{0}:=\{\fff\in \VVV:\ \fff%
(x)\in T\Sigma \, \text{ for a.a. }x\in \Gamma \}.
\end{equation*}
}	%
	

 Given a
closed subspace $W$ of a real Hilbert space $Z$, we denote by $W^{\perp_Z}$ its
orthogonal complement with respect to the $Z$-topology. We then indicate by
the symbol $W^\perp$ the annihilator of $W$, i.e.
\begin{equation*}
W^\perp:=\{x\in Z^{\prime}:\ \langle x,y\rangle_{Z^{\prime},Z}=0\quad
\forall\, y\in W\}.
\end{equation*}
We denote by $\mathcal{B}(X,Y)$ ($\mathcal{B}(X)$ when $X=Y$) the set of linear
bounded operators from the real Banach space $X$ to the real Banach space $Y$. Furthermore, given an operator $T$ from $X$ to $Y$, we define
by $T^\prime: Y^{\prime}\to X^{\prime}$ the adjoint of $T$.
In case of real Hilbert spaces, we denote the Hilbert adjoint of $T$ by $T^*: Y\to X$.

In order to use the arguments presented in \cite{GalPoia}, here we also recall the following well known Sobolev estimate, which can be proved by the definition of Sobolev spaces on smooth compact
manifold (and thus on the sphere $\Gamma$ we are considering in our setting) and standard partition of unity arguments:
\begin{lem}
Let $\Gamma$ be a smooth compact 2-dimensional manifold, then there exists $C>0$ such that
\begin{align}
\norm{f}_{L^p(\Gamma)}\leq C\sqrt{p}\norm{f}_{H^{1}(\Gamma)},\quad \forall f\in H^1(\Gamma).
	\label{M}
\end{align}
\end{lem}
\subsection{Constitutive assumptions}
\label{constassum}
 {\bf The mobility matrix}
 We now assume that the mobility $\mathbf{L}$
	is constant and positive definite over $T\Sigma $. Moreover, we suppose that:
	
	\begin{itemize}
		\item[(\textbf{M0})] there exists $l_{0}>0$
		such that
		\begin{equation}
		\mathbf{L}\boldsymbol{\eta }\cdot \boldsymbol{\eta }\geq l_{0}%
		\boldsymbol{\eta }\cdot \boldsymbol{\eta },\quad\forall\,\boldsymbol{%
			\eta }\in T\Sigma. \label{pos}
		\end{equation}
  \end{itemize}
	For the sake of simplicity we will adopt the compact notation $\mathbf{v}%
	\geq k$, with $\mathbf{v}\in \mathbb{R}^{N}$ and $k\in \mathbb{R}$ to
	indicate the relations $v_{i}\geq k$, $i=1,\ldots ,N$.
	\newcommand{\op}{\Delta_{\Gamma, L}^{-1}}
	
	Let us introduce the weighted inverse Laplacian as 
	the operator
	$\textstyle{-\op \colon \VVV_0' \to \VVV_0}$, which for each $\mathbf g\in \VVV_0'$ returns the unique element $\mathbf f := -\op\mathbf g\in \VVV_0$ defined by
	\begin{align}
		(\mathbf L \tgrad \mathbf f, \tgrad \boldsymbol\eta) = \langle \mathbf g, \boldsymbol\eta\rangle_{\VVV', \, \VVV}, \quad \text{ for all }  \boldsymbol\eta\in \VVV_0.
		\label{multilap}
	\end{align}
	The operator is well defined thanks to property \eqref{pos} of the matrix $\mathbf{L}$.
	Using the operator $\op$ we can also introduce an equivalent norm in $\VVV_0'$ by setting
	\begin{align*}
		\| \ggg \|_{-1,L}^2 := \left(\mathbf L \tgrad (-\op \ggg ), \tgrad (-\op \ggg)\right) = \langle \ggg, -\op \ggg\rangle_{\VVV', \, \VVV}.
	\end{align*}
	{\bf The entropy density}
	In order to include a fairly large admissible class of entropy densities, following the recent ideas presented in \cite{GalPoia}, we assume that%
	$$
	\boldsymbol\psi(\mathbf s):=\sum_{i=1}^N\psi(s_i),
	$$
	and
		\begin{equation}
		\left( \boldsymbol{\psi }'(\mathbf{s})\right) _{i}:=\psi ^{\prime }({s}_{i}),\quad i=1,\ldots ,N,
		\label{boldphi}\end{equation}%
	where
	
	\begin{equation*}
	\psi \in C\left[ 0,1\right] \cap C^{2}(0,1]
	\end{equation*}%
	has the following properties:
	
	\begin{itemize}
		\item[(\textbf{E0})] $\psi^{\prime \prime }(s)\geq \zeta >0,$ for all $s\in
		(0,1];$
		
		\item[(\textbf{E1})] $\lim_{s\rightarrow 0^{+}}\psi^{\prime }\left( s\right)
		=-\infty ;$
		
		\item[(\textbf{E2})]
		
As $\delta \rightarrow 0^{+}$, we assume, for some $\iota
>1/2$,
\begin{equation}
	-\frac{1}{\psi^{\prime }(2\delta )}=O\left( \frac{1}{|\ln (\delta )|^{\iota }}%
	\right) .  \label{est}
\end{equation}
	\end{itemize}

\begin{oss}
We point out that assumption (\textbf{E2}) above is much lighter than the usual assumption adopted, for instance, in \cite[(\textbf{E2})]{GGPS}.
However, thanks to \cite{GalPoia}, the results obtained in \cite{GGPS} still hold under this more general assumption on $\Psi$.
\end{oss}


 \an{As in \cite{GGPS}, we set $\psi(s)=+\infty ,$ for all $s\in (-\infty ,0)$, and we extend $\psi$ for all $s\in \lbrack 1,+\infty )$ so that $\psi$ is a $C^2$ function on $(0,+\infty)$ and  $(\textbf{E0})$ holds for any $s>0$. More precisely, we define
\begin{align}
\psi(s):=As^3+Bs^2+Ds,\quad\text{ for all }s\geq 1,
    \label{psiext}
\end{align}
with
\begin{align*}
    \begin{cases}
      A=\psi(1)-\psi'(1)+\frac 1 2 \psi''(1),\\
      B=-3\psi(1)+3\psi'(1)-\psi''(1),\\
      D=3\psi(1)-2\psi'(1)+\frac 1 2\psi''(1).
    \end{cases}
\end{align*}
}
 As it is immediate to verify, the logarithmic potential satisfies all the properties above.
 We also refer the reader to \cite{GalPoia} for some other important classes
	of mixing potentials that are singular at $0$.
	Furthermore, following the  scheme developed in \cite[Section 2]%
	{GGPS}, by (\textbf{E0})-(\textbf{E1}) we can define an approximation of the
	potential $\psi $ by means of a sequence $\left\{ \psi _{h
	}\right\} _{h >0}$ of everywhere defined non-negative functions.
	More precisely, let
	\begin{equation}
	\psi _{h }(s)=\frac{h }{2}|\mathbb{T}_{h
	}s|^{2}+\psi (J_{h }(s)),\qquad s\in \mathbb{R},\ h >0,
	\label{approx}
	\end{equation}%
	where $J_{h }=(I+h {\color{black} \mathbb{T}})^{-1}: \an{\R\to (0,+\infty)}$ is the resolvent
	operator and $\mathbb{T}_{h }=\frac{1}{h }%
	(I-J_{h })$ is the Yosida approximation of $\mathbb{T}\left(
	s\right) :=\psi^{^{\prime }}\left( s\right) ,$ for all $s\in \mathfrak{D}%
	\left( \mathbb{T}\right) =(0,1]$. We recall that the following properties hold:
	
	\begin{itemize}
		\item[(i)] $\psi _{h }$ is convex and $\psi _{h
		}(s)\nearrow \psi(s)$, for all $s\in \mathbb{R}$, as $h $ goes to $%
		0$;
		
		\item[(ii)] $\psi _{h }^{\prime }(s)=\mathbb{A}_{h }(s) $
		and \an{$\phi _{h }:=\psi _{h }^{\prime }$} is globally Lipschitz with constant $\frac{1}{h }$;
		
		\item[(iii)] $|\psi _{h }^{\prime }(s)|\nearrow |\psi^{\prime
		}(s)| $ for all $s\in (0,1]$ and $|\psi _{h }^{\prime
		}(s)|\nearrow +\infty ,$ for all $s\in (-\infty ,0]$, as $h $ goes
		to $0$;
		
		
		\item[(iv)] for any $h \in (0,1]$, there holds
		\begin{equation*}
		\psi _{h }^{\prime \prime }(s)\geq \frac{\zeta}{1+\zeta},\quad
		\text{for all}\ s\in \mathbb{R};
		\end{equation*}
		
		\item[(v)] for any compact subset $M\subset (0,1]$, $\psi _{h
		}^{\prime }$ converges uniformly to $\psi^{\prime }$ on $M$;
		
		\item[(vi)] for any $h _{0}>0$ there exists $\tilde{K}=\tilde{K}%
		(h _{0})>0$ such that
		\begin{equation*}
		\sum_{i=1}^{N}\psi _{h }(r_{i})\geq \frac{1}{4h _{0}}|%
		\mathbf{r}|^{2}-\tilde{K},\quad \forall \mathbf{r}\in \mathbb{R}^{N},\quad
		\forall 0<h <h _{0}.
		\end{equation*}
\end{itemize}
Let us introduce
		\begin{equation}
		\boldsymbol\Psi _{h }(\mathbf{r}):=\sum_{i=1}^{N}\psi _{h }(r_{i})-%
		\frac{1}{2}\mathbf{r}^{T}A\mathbf{r}=\boldsymbol\Psi_{h }^1(\mathbf{r})-%
		\frac{1}{2}\mathbf{r}^{T}A\mathbf{r},
  \label{Psidef}
		\end{equation}%
  \an{where, as presented in the Introduction, ${A}$ is a symmetric $N\times N$ matrix with $\lambda _{A}>0$ as the largest eigenvalue. }
		We then have that, for any $h _{0}>0$ sufficiently small, there
		exist $K=K(h _{0})>0$ and $C=C(h _{0})>0$, $%
		C(h _{0})\nearrow +\infty $ as $h _{0}\rightarrow 0$,
		such that
		\begin{equation}
\label{coercive}
		\boldsymbol\Psi _{h }(\mathbf{r})\geq C(h _{0})|\mathbf{r}%
		|^{2}-K,\quad \forall \mathbf{r}\in \mathbb{R}^{N},\quad \forall
		h \in(0,h _{0}).
		\end{equation}%
		In particular, this comes from the fact that $-\frac{1}{2}\mathbf{r}\cdot
		\mathbf{Ar}\geq -\frac{\lambda _{A}}{2}|\mathbf{r}|^{2}$ and $%
		h _{0}$ has to be small enough so that, e.g., $C(h _{0})=%
		\frac{1}{4h _{0}}-\frac{\lambda _{A}}{2}>0$.
\subsection{A recursive  inequality}
We recall the key tool for the application of De Giorgi's iteration argument to obtain the validity of the strict separation property. This lemma can be found, e.g., in \cite[Ch.2, Lemma 5.6]{Lady} (see also \cite{P} for a proof)
\begin{lem}
	\label{conv}
	Let $\{y_n\}_{n\in\N\cup \{0\}}\subset \R^+$ satisfy the recursive inequalities
	\begin{align}
		y_{n+1}\leq Cb^ny_n^{1+\gamma},
		\label{ineq}\qquad \forall n\geq 0,
	\end{align}
	for some $C>0$, $b>1$ and $\gamma>0$. If
	\begin{align}
		\label{condition}
		y_0\leq \theta:= C^{-\frac{1}{\gamma}}b^{-\frac{1}{\gamma^2}},
	\end{align}
	then
	\begin{align}
		y_n\leq \theta b^{-\frac{n}{\gamma}},\qquad \forall n\geq 0,
		\label{yn}
	\end{align}
	and consequently $y_n\to 0$ for $n\to +\infty$.
\end{lem}

\section{Main results} \label{mainres}

Here we first state a weak formulation of the Cauchy problem associated with system \eqref{MAIN} in the case $\mathbf L(\vphi)\equiv \mathbf L$ constant. Then we will state our main results.

Let us set
\begin{align*}
    \mathcal K :=
    \left\{
    (\vphi, u)\in \VVV \times H^2(\Gamma) \colon \;\dashint_{\Gamma} \vphi = \pmb \alpha \,\, \text{ and } \,\, u\in \operatorname{span }\{1,\nu_1, \nu_2, \nu_3\}^\perp
    \right\},
\end{align*}
where $\al\in(0,1)^N$, $\al\in \Sigma$, \Andrea{and $\nu_i$ are the components of the normal vector $\nu$ to $\Gamma_0$}.
and define the slices
\begin{align*}
    \mathcal K_1 := \left\{
    \vphi\in \VVV \colon \; \dashint_{\Gamma} \vphi = \pmb \alpha
    \right\}, \quad \mathcal K_2 := \left\{
    u\in H^2(\Gamma) \colon\;  u\in \operatorname{span }\{1,\nu_1, \nu_2, \nu_3\}^\perp
    \right\}.
\end{align*}
The projection of $\R^N$ onto $T\Sigma$ (see \eqref{sigma2}) is given by
\begin{align}\label{5_eq:proj}
    \mathbf P\colon \R^N \to T\Sigma, \quad \mathbf P\mathbf v = \mathbf v - \left(\dfrac{1}{N}\sum_{i=1}^N v_j\right) \mathbf e,
\end{align}
where $\mathbf e = (1,\dots,1)^T$. Componentwise this reads, for $l=1,\dots,N$,
\begin{align*}
    (\mathbf P \mathbf v)_l = \dfrac{1}{N} \sum_{m=1}^N (v_l - v_m).
\end{align*}
We then notice that, due to the property $$\sum_{j=1}^N L_{kj} = 0,\quad \forall k=1,\ldots, N,$$ the diffusion equation for $\vphi_k$ becomes
\begin{align*}
    \partial_t \vphi_k = (\nabla_\Gamma \cdot \left( \mathbf L \nabla_\Gamma \pmb\mu\right))_k = \nabla_{\Gamma} \cdot \left(\sum_{j=1}^N L_{kj} \nabla_\Gamma \mu_j \right) =  \sum_{j=1}^N L_{kj} \Delta_\Gamma \dfrac{1}{N}\sum_{m=1}^N (\mu_j-\mu_m);
\end{align*}
in other words, the diffusion equations are determined by the chemical potential differences $(\mu_j-\mu_k)$, and it is thus convenient to introduce the \textit{vector of generalised chemical potential differences}
\begin{align*}
    \mathbf w = \dfrac{1}{N} \left( \sum_{m=1}^N (\mu_j - \mu_m)\right)_{j=1,\dots,N} = \mathbf P \boldsymbol \mu,
\end{align*}
where $\mathbf P$ is given by \eqref{5_eq:proj}. We therefore rewrite system \eqref{MAIN} as follows
\begin{align}
\begin{split}
    \partial_t \vphi &=\nabla_\Gamma\cdot ( \mathbf L \nabla_\Gamma \mathbf w) \\
    \mathbf w &= - b\varepsilon \Delta_\Gamma \vphi+ \mathbf P\left( \dfrac b\varepsilon \left(\boldsymbol\psi'(\vphi) - A\vphi\right) + \dfrac{2\kappa u \vlambda}{R^2} + \kappa (\vlambda \cdot \vphi) \vlambda + \kappa \Delta_\Gamma u \, \vlambda\right) \\
    \beta \partial_t u &=  -\kappa \Delta_\Gamma^2 u + \left(\sigma - \dfrac{2\kappa}{R^2} \right) \Delta_\Gamma u + \dfrac{2\sigma u}{R^2} - \kappa \vlambda \cdot \Delta_\Gamma \vphi - \dfrac{2\kappa \vlambda\cdot(\vphi-\pmb \alpha)}{R^2}.
    \end{split}
\end{align}
Notice that, since $\sum_{i=1}^N\vphi_i=1$, it holds $\textbf{P}\Delta_\Gamma\vphi=\Delta_\Gamma\vphi$.

\noindent We can now introduce the following weak formulation of the Cauchy problem for our system:

\begin{prob}\label{5_pr}
Given $(\vphi_0,u_0)\in \mathcal K$ such that $0\leq \vphi_0\leq 1$, and $\sum_{i=1}^N\vphi_{0,i}=1$, find $(\vphi, \pmb\mu, u)$ satisfying:
\begin{itemize}
    \item[(i)] the properties
    \begin{align}
    \begin{cases}
    \label{wregprop}
        &\vphi \in L^\infty(0,T; \VVV) \cap L^2(0,T; \HHH^2(\Gamma)), \,\, \text{ with } \,\, \partial_t \vphi\in L^2(0,T; \VVV'); \\
        &0 < \vphi < 1 \,\,\quad \sum_{i=1}^N\vphi_i=1,\ \text{ a.e. in } \Gamma \times (0,T); \\
        &\mathbf w \in L^2(0,T; \VVV); \\
        &u \in L^\infty(0,T; H^2(\Gamma))\cap L^2(0,T; H^4(\Gamma)),\ \text{with} \,\, \partial_t u\in L^2(0,T; H);
    \end{cases}
    \end{align}
    \item[(ii)] the equations
    \begin{align}
    \label{weakeq}
    \begin{cases}
        \langle \partial_t \vphi, \peta\rangle + (\mathbf L \nabla_\Gamma \mathbf w, \nabla_\Gamma \peta) = 0, \quad \forall \peta \in \VVV \\
        (\mathbf w, \pxi) = \dfrac{b}{\varepsilon} (\mathbf P\boldsymbol\psi'(\vphi), \pxi) - \dfrac{b}{\varepsilon}(\mathbf PA\vphi, \pxi) + \dfrac{2\kappa}{R^2} (\mathbf Pu\boldsymbol\Lambda, \pxi) + b\varepsilon (\nabla_\Gamma \vphi, \nabla_\Gamma \pxi) \\
        + \kappa( \mathbf P(\boldsymbol\Lambda\cdot\vphi)\boldsymbol\Lambda, \pxi) + \kappa ( \mathbf P\boldsymbol\Lambda\Delta_\Gamma u, \pxi), \quad \forall \pxi \in \VVV \\
        \beta \langle \partial_t u, \zeta\rangle = -(\Delta_\Gamma u, \Delta_\Gamma \zeta) - \left( \sigma - \dfrac{2\kappa}{R^2}\right) (\nabla_\Gamma u, \nabla_\Gamma \zeta) + \dfrac{2\sigma}{R^2} (u, \zeta) \\+ \kappa( \nabla_\Gamma(\boldsymbol\Lambda\cdot\vphi), \nabla_\Gamma \zeta) - \dfrac{2\kappa}{R^2} (\boldsymbol\Lambda \cdot(\vphi-\pmb\alpha), \zeta), \quad \forall \zeta\in H^2(\Gamma),
    \end{cases}
    \end{align}
    almost everywhere in $(0,T)$;
    \item[(iii)] the initial conditions $\vphi(0)=\vphi_0$, $u(0)=u_0$.
\end{itemize}
\end{prob}

Problem \ref{5_pr} is well posed. Indeed, we have

\begin{thm}
\label{maina}
Given $(\vphi_0,u_0)\in \mathcal K$ such that $0\leq \vphi_0\leq 1$, and $\sum_{i=1}^N\vphi_{0,i}=1$, there exists a unique solution $(\vphi,\pmb\mu,u)$ to Problem \ref{5_pr} with initial data $(\vphi_0, u_0)$. Furthermore, the following time-weighted estimates hold, namely, there exists $C>0$, possibly depending on $T>0$ and $r\geq2$, such that
 \begin{itemize}
	\item[(i)] $\|\sqrt{t}\mathbf w\|_{L^\infty(0,T;\VVV))}+\|\sqrt{t}\boldsymbol\psi'(\vphi)\|_{L^\infty(0,T; \HHH)} \leq C$;
	\item[(ii)] $\| \sqrt{t}{{\vphi}}\|_{L^\infty(0,T; \HHH^2(\Gamma))} + \|\sqrt{t}\partial_t \vphi\|_{L^2(0,T;\VVV)}
+\|\sqrt{t} u\|_{L^\infty(0,T; H^4(\Gamma))} \leq C$;
	\item[(iii)] $\|\sqrt{t}\mathbf w\|_{L^\infty(0,T;\LLL^r(\Gamma)))}+\|\sqrt{t}\boldsymbol\psi'(\vphi)\|_{L^\infty(0,T; \LLL^r(\Gamma))}\leq C\sqrt{r},\qquad \forall r\geq 2;$
	\item[(iv)] $\| \sqrt{t}{{\vphi}}\|_{L^\infty(0,T;\WWW^{2,r}(\Gamma))}\leq C(r),\qquad \forall r\geq 2$.
\end{itemize}
Moreover, $\vphi \in C([0,T];\VVV)$ and the solution satisfies the energy identity
\begin{align}
    \mathcal E (u(t),\vphi(t))  +\int_s^t\beta\|\partial_t u(r)\|^2 dr+\int_s^t (\mathbf L \tgrad\mathbf w(r), \tgrad \mathbf w(r))dr=\mathcal E (\vphi(s), u(s)),
    \label{enerident}
\end{align}
for any $t\in (0,T]$ and almost any $s\geq 0$, $s=0$ included.

Let now $(\vphi_i,\mathbf w_i, u_i)$, $i=1,2$ be solutions to Problem \ref{5_pr} with initial data
\begin{align*}
	\vphi_i(0) = \vphi_{i,0}\in\mathcal K_1,\quad 0\leq \vphi_{i,0}\leq1, \quad \sum_{j=1}^N\vphi^j_{i,0}\equiv 1, \quad u_i(0) = u_{i,0}\in\mathcal K_2, \quad i=1,2.
\end{align*}
Then, for any $T>0$, there exists $C(T)>0$ such that, for almost all $t\in (0,T)$, the following continuous dependence estimate holds
\begin{align}\label{5_eq:contdep}
	\|\vphi_1(t) - \vphi_2(t)\|_{-1,L}^2 + &\| u_1(t) - u_2(t)\|^2 \leq C\left( \|\vphi_{1,0} - \vphi_{2,0}\|_{-1,L}^2 + \| u_{1,0} - u_{2,0}\|^2 \right).
\end{align}
\end{thm}

\begin{oss}
Note that $u\in C([0,T];H^2(\Gamma))$, thanks to a classical embedding.
Moreover, the above weighted estimates combined with elliptic regularity results entail that the weak solution
gets instantaneously strong. In particular, we also have $\sqrt{t}\mathbf w \in L^2(0,T;\mathbf H^2(\Gamma))$ and
$\sqrt{t}u \in L^2(0,T;H^4(\Gamma))$.
\end{oss}

\begin{oss}
\label{5_constraint}
Notice that, arguing as in \cite[Prop.2.1]{EllLuc91}, we easily obtain that $\sum_{i=1}^N\vphi_i=1$ and $\sum_{i=1}^N\mathbf w_i=0$. Moreover, by choosing $\eta\equiv 1$
we get that the total mass is preserved, i.e., $\vphi(t)\in \mathcal{K}_1$ for any $t\in [0,T]$. In conclusion, choosing $\zeta=1,\nu_1,\nu_2,\nu_3$, respectively,
we infer that $u(t)\in \mathcal{K}_2$ for any $t\in[0,T]$.
\end{oss}

%

Thanks to the energy identity, the (unique) weak solution given by Theorem \ref{maina}
can be easily defined over the time interval $[0,+\infty)$. This fact also allows us
to prove that the instantaneous separation property from pure phases holds globally in time. More precisely, we have
\begin{thm}\label{mainb}
Let $(\vphi_0,u_0)\in \mathcal K$ such that $0\leq \vphi_0\leq 1$, and $\sum_{i=1}^N\vphi_{0,i}=0$. Then the unique solution $(\vphi,\mathbf w)$ to Problem \ref{5_pr}
is defined in $[0,+\infty)$  and is such that, for any $\tau>0$,
	\begin{align}
\label{globreg}
\begin{cases}
	&\vphi\in L^\infty(\tau,+\infty;\textbf{W}^{\,2,p}(\Gamma)),\quad \text{for } p\in[2,+\infty),\\&
 \partial_t\vphi\in L^2(t,t+1;\VVV),\quad \forall t\geq \tau,\\&
\mathbf w\in L^\infty(\tau,+\infty;\VVV),\\&
\boldsymbol\psi^\prime(\vphi)\in L^\infty(\tau,+\infty;\LLL^r(\Gamma)),\quad \text{for } r\in[2,+\infty),
	\\&
       u \in L^\infty(\tau,+\infty; H^4(\Gamma))\cap L^2(t,t+1; H^4(\Gamma)),\quad\forall t\geq \tau,\\&
      \partial_t u\in L^\infty(\tau,+\infty; H^2(\Gamma)).
      \end{cases}
	\end{align}
 Moreover, \eqref{enerident} is valid for any $t>0$ and almost any $s\geq 0$ ($s=0$ included) and
 the strict separation property holds, i.e., for any $0<\tau<T$, there exists $0<{\delta}=\delta(\tau)<1/N$ such that
 \begin{equation}
 \label{globalssp}
 \delta<\vphi< 1-(N-1)\delta, \quad \text{ in } \Gamma\times[\tau,+\infty).
 \end{equation}
\end{thm}

Finally, we state the result about the convergence of any weak solution to a single stationary state. In order
to do that, we shall use a suitable version of the well-known {\L}ojasiewicz-Simon inequality (see Lemma \ref{Loja} below).

Let us first introduce the space
\begin{align*}
	\mathcal{V}_{\al}:=&\{(\vphi,u)\in \VVV\times H^2(\Gamma):\ 0\leq
	\vphi(x)\leq 1, \; \text{ for a.a. }x\in \Gamma ,\quad \dashint_\Gamma \vphi =\al,\quad \sum_{i=1}^{N}\vphi_i = 1,\\&\quad u\in \text{span}\{1,\nu_1,\nu_2,\nu_3\}^\perp\},
\end{align*}%
endowed with the $\VVV\times H^2(\Gamma)$-topology. This is a complete metric space. Given $(\vphi_0,{u}_{0})\in \mathcal{V}_{\al}$, we consider the (unique) weak solution $(\vphi,u)$ given by Theorem \ref{maina} and introduce
the $\omega $-limit set $\omega (\vphi_0,{u}_{0})$ as
\begin{align*}
&\omega (\vphi_0,{u}_{0})\\
&=\{(\mathbf{z},\tilde{z})\in (\mathbf{H}^{2r}(\Gamma )\times H^{4r}(\Gamma))\cap
	\mathcal{V}_{\al}:\exists\,t_{n} \nearrow +\infty \text{ s.t. } (\vphi(t_{n}),u(t_n))\rightarrow (\mathbf{z},\tilde{z}) \text{ in }\mathbf{H}^{2r}(\Gamma)\times H^{4r}(\Gamma)\},
\end{align*}%
where $r\in \lbrack \tfrac{1}{2},1)$. In particular, we
fix $r\in (\tfrac{1}{2},1)$. In order to show that $\ol$ is actually a singleton, since the set of stationary points is
uncountable, the standard tool is the {\L}ojasiewicz-Simon inequality (see Lemma \ref{Loja} below)
which requires the real analyticity of the nonlinearities. In our case, we assume
\begin{itemize}
\item[(\textbf{E3})] $\psi$ is (real) analytic in $(0,1)$.
\end{itemize}
Note that this assumption is satisfied if, for instance, $\psi(s)=s\ln s$, $s\in (0,1]$, i.e., for the standard logarithmic potential.

We then have the following

\begin{thm}
\label{maincon} {\color{black} Let the assumptions listed in Subsection \ref{constassum} hold, but (\textbf{E2}). In addition, assume (\textbf{E3})}.
Then, for any $(\vphi_0,u_0)\in \mathcal{V}_{\al}$, it holds $\ol=\{(\vphi_\infty,{u}_\infty)\}$, where $(\vphi_\infty,u_\infty)$ is a stationary solution, that is,
$(\vphi_\infty,u_\infty)\in H^2(\Gamma) \times H^4(\Gamma)$ solves the system
\begin{align}
\label{steady2b1}
\begin{cases}
  &- b\varepsilon \Delta_\Gamma \vphi_\infty+\frac{b}{\varepsilon}\mathbf{P}\boldsymbol\psi'(\vphi_\infty) = \pmb{f}+\mathbf{P}\left(\dfrac b\varepsilon  A\vphi_\infty-\dfrac{2\kappa u_\infty \vlambda}{R^2}  - \kappa (\vlambda \cdot \vphi_\infty) \vlambda - \kappa \Delta_\Gamma u_\infty \, \vlambda\right) \\
  &-\kappa \Delta_\Gamma^2 u_\infty + \left(\sigma - \dfrac{2\kappa}{R^2} \right) \Delta_\Gamma u_\infty + \dfrac{2\sigma u_\infty}{R^2} - \kappa \vlambda \cdot \Delta_\Gamma \vphi_\infty
  - \dfrac{2\kappa \vlambda\cdot(\vphi_\infty-\pmb \alpha)}{R^2} =0
    \end{cases}
\end{align}
almost everywhere on $\Gamma$, with
\begin{equation*}
\pmb{f}=\dashint_\Gamma\frac{b}{\varepsilon}\mathbf{P}\boldsymbol\psi'(\vphi_\infty) -\dashint_\Gamma\mathbf{P}\left(\dfrac b\varepsilon  A\vphi_\infty-\dfrac{2\kappa u_\infty \vlambda}{R^2}  - \kappa (\vlambda \cdot \vphi_\infty) \vlambda - \kappa \Delta_\Gamma u_\infty \, \vlambda\right).
\end{equation*}
Moreover, there exists $\delta>0$ so
that
\begin{equation*}
\delta<\vphi_\infty(x) <1-(N-1)\delta,\quad \forall\,x\in \Gamma,
\end{equation*}
and
\begin{align}
(\vphi(t),u(t))\underset{ t\to +\infty}{\longrightarrow}
(\vphi_\infty,{u}_\infty)\quad \text{in }\mathbf{H}^{2r}(\Gamma)\times H^{4r}(\Gamma),  \label{sing}
\end{align}
for any $r\in(0,1)$.
\end{thm}

\begin{oss}
We have all the ingredients to study our problem as a dissipative dynamical system in a suitable phase space (see \cite{Temam}), proving the
existence of the global attractor as well as of an exponential attractor (cf. \cite[Rem.3.10]{GGPS}).
\end{oss}

\section{Proof of Theorem \ref{maina}}\label{proof1}

We first prove the continuous dependence estimate \eqref{5_eq:contdep}, which entails uniqueness.
\paragraph{Proof of estimate \eqref{5_eq:contdep}.}
Set $\pxi^{\vphi} = \vphi^1 - \vphi^2$, $\pxi^{\mathbf w} = \mathbf w_1 - \mathbf w_2$, $\xi^u = u_1 - u_2$. Recalling the definition of the weighted inverse Laplacian given in \eqref{multilap},
we test the equation satisfied by $\pxi^{\vphi}$ with $(\Delta^L_\Gamma)^{-1} \pxi^{\vphi}\in \VVV_{0}$ and noting that $\pxi^{\vphi} = \mathbf P \pxi^{\vphi}$, we get
\begin{align*}
	\dfrac 12 \dfrac{d}{dt} \|\pxi^{\vphi}\|_{-1,L}^2 + \dfrac{b}{\varepsilon} (\boldsymbol\psi'(\vphi^1)-\boldsymbol\psi'(\vphi^2), \pxi^{\vphi}) &+ \dfrac{2\kappa}{R^2} (\xi^u \pmb \Lambda, \pxi^{\vphi}) - ( A \mathbf \pxi^{\vphi}, \pxi^{\vphi}) \\
	&+ b\varepsilon \|\tgrad \pxi^{\vphi}\|^2 + \kappa \|\pmb \Lambda\cdot \pxi^{\vphi}\|^2 + \kappa (\Delta_\Gamma \xi^u\pmb \Lambda, \pxi^{\vphi}) = 0.
\end{align*}
Due to the convexity of $\psi$, for each component $\vphi^1_i, \vphi^2_i$, we have
\begin{align*}
	\int_\Gamma \left(\psi'(\vphi^1_i) - \psi'(\vphi^2_i) \right) \left( \vphi_i^1 - \vphi_i^2\right) \geq 0.
\end{align*}
Observe now that
\begin{align*}
	( A \mathbf \pxi^{\vphi}, \pxi^{\vphi}) \leq C_1\|\pxi^{\vphi}\|^2 = C_1 ( \mathbf L\tgrad (\Delta_L^{-1} \pxi^{\vphi}, \tgrad \pxi^{\vphi}) &\leq C_2 \| \pxi^{\vphi}\|_{-1,L} \, \|\tgrad \pxi^{\vphi}\| \\
	&\leq C_3 \| \pxi^{\vphi}\|_{-1,L}^2 + \dfrac{b\varepsilon}{4} \|\tgrad \pxi^{\vphi}\|^2.
\end{align*}
and similarly
\begin{align*}
	(\xi^u \pmb \Lambda, \pxi^{\vphi}) \leq C_1\|\xi^u\|^2 + C_2 \|\pxi^{\vphi}\|^2 \leq  C_1\|\xi^u\|^2 +C_3 \| \pxi^{\vphi}\|_{-1,L}^2 + \dfrac{b\varepsilon}{4} \|\tgrad \pxi^{\vphi}\|^2
\end{align*}
as well as
\begin{align*}
	(\Delta_\Gamma \xi^u\pmb \Lambda, \pxi^{\vphi}) \leq \dfrac 14 \|\Delta_\Gamma \xi^u\|^2 + C_1\|\pxi^{\vphi}\|^2 \leq \dfrac 14 \|\Delta_\Gamma \xi^u\|^2 + C_2 \| \pxi^{\vphi}\|_{-1,L}^2 + \dfrac{b\varepsilon}{4} \|\tgrad \pxi^{\vphi}\|^2.
\end{align*}
Putting everything together leads to
\begin{align}
\label{diffpsi}
	\dfrac 12 \dfrac{d}{dt} \|\pxi^{\vphi}\|_{-1,L}^2 + \dfrac{b\varepsilon}{4} \|\tgrad \pxi^{\vphi}\|^2 \leq C_1 \left(\| \pxi^{\vphi}\|_{-1,L}^2 + \|\xi^u\|^2\right) + \dfrac 14\|\Delta_\Gamma \xi^u\|^2.
\end{align}
Testing the equation satisfied by $\xi^u$ with $\xi^u$ itself and integrating by parts, we get
\begin{align*}
	\dfrac{\beta}{2} \dfrac{d}{dt} \|\xi^u\|^2 + \|\Delta_\Gamma \xi^u\|^2 &\leq C_1 \| \xi^u\|^2 + \dfrac{1}{4} \|\Delta_\Gamma \xi^u\|^2 + C_2 \|\xi^u\|^2 + C_3 \|\pxi^{\vphi}\|^2 + \dfrac 14 \|\Delta_\Gamma \xi^u\|^2 \\
	&\leq C_3 \|\xi^u\|^2 + \dfrac 12 \|\Delta_\Gamma \xi^u\|^2 + C_4 \|\pxi^{\vphi}\|_{-1,L}^2 + \dfrac{b\varepsilon}{8} \|\tgrad \pxi^{\vphi}\|^2,
\end{align*}
that is,
\begin{align}
\label{diffu}
	\dfrac{\beta}{2}\dfrac{d}{dt} \|\xi^u\|^2 + \dfrac 12 \|\Delta_\Gamma \xi^u\|^2  \leq C_3 \|\xi^u\|^2 + C_4 \|\pxi^{\vphi}\|_{-1,L}^2 + \dfrac{b\varepsilon}{8} \|\tgrad \pxi^{\vphi}\|^2
\end{align}
Adding \eqref{diffpsi} and \eqref{diffu} together, we obtain
\begin{align*}
	\dfrac{d}{dt} \left[ \dfrac 12 \|\pxi^{\vphi}\|_{-1,L}^2 +  \dfrac{\beta}{2} \|\xi^u\|^2 \right] + \dfrac{b\varepsilon}{8}\|\tgrad \pxi^{\vphi}\|^2 + \dfrac 14 \|\Delta_\Gamma \xi^u\|^2 \leq C \left( \dfrac 12 \|\pxi^{\vphi}\|_{-1,L}^2 +  \dfrac{\beta}{2} \|\xi^u\|^2 \right).
\end{align*}
Gronwall's inequality implies \eqref{5_eq:contdep}. The proof is finished.

\paragraph{Existence of solutions: an approximate problem.}
In the sequel we will adopt the notation $\boldsymbol\psi'_\eps$ and $\boldsymbol\Psi^2_{,\vphi}$ for the vectors whose components are defined as follows
\begin{align*}
 (\boldsymbol\psi'_\eps)_i=\dfrac{\partial\mathbf\Psi^1_h}{\partial\vphi_i} \quad \text{ and } \quad (\mathbf\Psi_{,\vphi}^2)_i=\dfrac{\partial\mathbf\Psi^2}{\partial\vphi_i}
\end{align*}
for $i=1,\ldots,N$, respectively. Note that we then have
$$
(\boldsymbol\psi'_\eps)_i=\psi^\prime_\eps(\vphi_i),\quad \mathbf\Psi^2_{,
\vphi}(\vphi)=-A\vphi.
$$
We also define the approximate energy
\begin{align*}
    &\mathcal E_\eps(u,\vphi) = \mathrm E_{\mathrm H}(u,\vphi) + \mathrm E_{\mathrm{CH},\eps}(\vphi) := \mathrm E_{\mathrm H}(u,\vphi) + \int_{\Gamma_0}\dfrac{b\varepsilon}{2}|\tgrad  \vphi|^2 + \dfrac{b}{\varepsilon} \boldsymbol\Psi_\eps(\vphi)\d\Gamma .
\end{align*}
The approximate regularised problem is then formulated as follows
\begin{prob}\label{5_prob:regularised}
Let $\eps>0$ be sufficiently small. Given $(\vphi_0,u_0)\in \mathcal K$, with $0\leq \vphi_0\leq 1,\  \sum_{i=1}^N\vphi_{0,i}=1$, find $(\vphi_\eps, \mathbf w_\eps, u_\eps)$ satisfying
\begin{itemize}
    \item[(i)] the regularity properties
    \begin{align*}
        &\vphi_\eps \in L^\infty(0,T; \VVV) \cap L^2(0,T; \HHH^2(\Gamma)) \,\, \text{ with } \,\, \partial_t \vphi_\eps\in L^2(0,T; \VVV'); \\
        &\mathbf w_\eps \in L^2(0,T; \VVV); \\
        &u_\eps \in L^\infty(0,T; H^2(\Gamma))\cap L^2(0,T; H^4(\Gamma)) \,\, \text{ with } \,\, \partial_t u_\eps\in L^2(0,T; H)
    \end{align*}
    \item[(ii)] the equations
    \begin{align*}
        \langle \partial_t \vphi_\eps, \peta\rangle + (L \nabla_\Gamma \mathbf w_\eps, \nabla_\Gamma \peta) = 0, &\quad \forall \peta \in \VVV \\
        (\mathbf w_\eps, \pxi) = \dfrac{b}{\varepsilon} (\mathbf P\boldsymbol\psi_\eps'(\vphi_\eps), \pxi) - \dfrac{b}{\varepsilon}(\mathbf PA\vphi_\eps, \pxi) + \dfrac{2\kappa}{R^2} (\mathbf Pu_\eps\boldsymbol\lambda, \pxi) + b\varepsilon (\nabla_\Gamma \vphi_\eps, \nabla_\Gamma \pxi) \\
        + \kappa( \mathbf P(\boldsymbol\lambda\cdot\vphi_\eps)\boldsymbol\lambda, \pxi) + \kappa (\mathbf P\boldsymbol\lambda\Delta_\Gamma u_\eps, \pxi), &\quad \forall \pxi \in \VVV \\
        \beta \langle \partial_t u_\eps, \zeta\rangle = -(\Delta_\Gamma u_\eps, \Delta_\Gamma \zeta) - \left( \sigma - \dfrac{2\kappa}{R^2}\right) (\nabla_\Gamma u_\eps, \nabla_\Gamma \zeta) + \dfrac{2\sigma}{R^2} (u_\eps, \zeta) \\+ \kappa( \nabla_\Gamma(\boldsymbol\lambda\cdot\vphi_\eps), \nabla_\Gamma \zeta) - \dfrac{2\kappa}{R^2} (\Lambda \cdot(\vphi_\eps-\pmb\alpha), \zeta), &\quad \forall \zeta\in H^2(\Gamma);
    \end{align*}
    \item[(iii)] the initial conditions $\vphi_\eps(0)=\vphi_0$, $u_\eps(0)=u_0$.
\end{itemize}
\end{prob}

\begin{oss}
Remark \ref{5_constraint} also holds for problem \ref{5_prob:regularised}.
\label{5_constraint2}
\end{oss}

\paragraph{Galerkin scheme.}

We establish well-posedness for problem \ref{5_prob:regularised} using a Galerkin scheme. This is standard (see also the proof of \cite[Theorem 3.1]{AC2023}), so we only describe it briefly. Denote by $\{\zzz_j\colon j\in \N\}$ the family of smooth eigenfunctions of the component-wise Laplace-Beltrami operator $-\Delta_\Gamma$ on $\Gamma$ subject to the constraint of zero integral average and sum of the $N$ components equal zero (corresponding to the eigenvalues $\{\widetilde{\lambda}_j\}_{j\in\N}$). Namely, each $\zzz_j$ solves
\begin{align*}
	\begin{cases}
		-\Delta_\Gamma z_{j,i}=\widetilde{\lambda}_jz_{j,i},\quad \forall i=1,\ldots,N,\\
		\dashint_\Gamma\zzz_j=\textbf{0},\\
		\sum_{i=1}^Nz_{j,i}=0.
	\end{cases}
\end{align*}
Following what developed in \cite[Appendix 8.1]{GGPS}, which is valid also if we substitute the domain $\Omega$ with the closed surface $\Gamma$, $\{\zzz_j\}_{j\in\N}$ forms an orthogonal basis for $\VVV_0$ and an orthonormal basis for $\HHH_0$. Defining the finite-dimensional spaces
\begin{align*}
	\Vb_m := \{\zzz_1,\dots, \zzz_m\}, \quad m\in \N,
\end{align*}
and the corresponding $\LLL^2$-projection operators
\begin{align*}
	\Pb_m \colon \HHH_0\to \Vb_m, \,\,
\end{align*}
we have the convergence
\begin{align*}
	\Pb_m \pmb\chi \to \pmb\chi \, \text{ in } \HHH_0, \,\,\, &\forall \pmb\chi\in \HHH_0, \\
	\Pb_m \pmb\chi \to \pmb\chi \, \text{ in } \VVV_0, \,\,\, &\forall \pmb\chi\in \VVV_0,
\end{align*}
as well as the uniform bounds
\begin{align}
	\|\Pb_m \pmb\chi\| \leq \|\pmb\chi\|, \quad \quad \|\Pb_m \pmb\chi\|_{\VVV} \leq C\|\pmb\chi\|_{\VVV}.
	\label{proj}
\end{align}
We also denote by $\{f_j\colon j\in \N\}$ the family of smooth eigenfunctions of the Laplace-Beltrami operator $-\Delta_\Gamma$ on $\Gamma$ (corresponding to the eigenvalues $\{\lambda_j\}_{j\in\N}$), forming an orthogonal basis for $V$ and an orthonormal basis for $H$. We the finite-dimensional spaces
\begin{align*}
	\WW_m := \{f_1,\dots, f_m\}, \quad m\in \N
\end{align*}
and observe that $f_1=const$ and corresponds to $\lambda_1=0$. Moreover, since $-\Delta_\Gamma\nu_i=\frac 2 {R^2}\nu_i$, for $i=1,2,3$, there exists $\overline{m}$ such that $\lambda_{\overline{m}}=\frac 2 {R^2}$ and $\nu_1,\nu_2,\nu_3$ are orthogonal eigenfunctions belonging to the corresponding eigenspace. This means that, for $m\geq\overline{m}$, $\{1,\nu_1,\nu_2,\nu_3\}\in \WW_m$. From now on we thus consider $m\geq\overline{m}$.
We also introduce the corresponding $L^2$-projection operator as
\begin{align*}
	\PP_m \colon H\to \mathcal W_m, \,\,
\end{align*}
and we have the convergence
\begin{align*}
	\mathcal P_m \chi \to \chi \, \text{ in } H, \,\,\, &\forall \chi\in H, \\
	\mathcal P_m \chi \to \chi \, \text{ in } V, \,\,\, &\forall \chi\in V, \\
	\mathcal P_m \chi \to \chi \, \text{ in } H^2(\Gamma), \,\,\, &\forall \chi\in H^2(\Gamma),
\end{align*}
as well as the uniform bounds
\begin{align}
	\|\mathcal P_m \chi\| \leq \|\chi\|, \quad \quad \|\mathcal P_m \chi\|_{V} \leq \|\chi\|_{V} \quad \text{ and } \quad \|\mathcal P_m \chi\|_{H^2(\Gamma)} \leq \|\chi\|_{H^2(\Gamma)}.
	\label{proj2}
\end{align}

We now set up the Galerkin approximation to \eqref{5_prob:regularised}. First, we recall that $\pmb\alpha=\dashint_\Gamma\vphi_0$. Then we consider

\begin{prob}
\label{5_gal}
For each $m\in\N$, find $(\vphi^m_\eps, \mathbf w^m_\eps, u^m_\eps)$ of the form
\begin{align}
    \vphi^m_{\eps,i}(t) = \sum_{k=1}^m \alpha^i_k(t) \zzz_k, \quad \mathbf w^m_{\eps,i}(t)-\dashint_\Gamma{\mathbf w}_{\eps,i}^m(t) = \sum_{k=1}^m \beta^i_k(t)  \zzz_k, \quad u^m_\eps(t) = \sum_{k=1}^m \gamma_k(t) f_k,
    \label{5_local}
\end{align}
with $\alpha_k^i,\beta_k^i,\gamma_k\in C^1([0,T])$ for $i=1,\ldots,N,\ k=1,\ldots,m$, such that
\begin{align}
       \label{e1} &( \partial_t \vphi_\eps^m, \peta^m) + (\mathbf L \nabla_\Gamma \mathbf w_\eps^m, \nabla_\Gamma \peta^m) = 0, \\
        &\label{e2}(\mathbf w_\eps^m-\dashint_\Gamma{\mathbf w}_\eps^m, \pxi^m) = \dfrac{b}{\varepsilon} (\mathbf P\boldsymbol\psi_\eps'(\vphi_\eps^m+\al), \pxi^m) - \dfrac{b}{\varepsilon}(\mathbf PA(\vphi_\eps^m+\al), \pxi^m) \\
        &\nonumber+ \dfrac{2\kappa}{R^2} (\mathbf Pu_\eps^m\boldsymbol\Lambda, \pxi^m) + b\varepsilon (\nabla_\Gamma \vphi_\eps^m, \nabla_\Gamma \pxi^m)
       + \kappa( \mathbf P(\boldsymbol\Lambda\cdot(\vphi_\eps^m+\al))\boldsymbol\Lambda, \pxi^m) + \kappa (\mathbf P\boldsymbol\Lambda \Delta_\Gamma u_\eps^m, \pxi^m),\\\label{e3}
       &\dashint_\Gamma{\mathbf w}_{\eps}^m=\dashint_\Gamma\left({\frac b\varepsilon\mathbf P\boldsymbol\psi'_\eps(\vphi_\eps^m+\al)-\frac b \varepsilon\mathbf PA(\vphi_\eps^m+\al)+\dfrac{2\kappa}{R^2} \mathbf Pu_\eps^m\boldsymbol\Lambda+ \kappa \mathbf P(\boldsymbol\Lambda\cdot(\vphi_\eps^m+\al))\boldsymbol\Lambda}\right),
         \\\label{e4}\nonumber
        &\beta \langle \partial_t u_\eps^m, \zeta^m\rangle = -(\Delta_\Gamma u_\eps^m, \Delta_\Gamma \zeta^m) - \left( \sigma - \dfrac{2\kappa}{R^2}\right) (\nabla_\Gamma u_\eps^m, \nabla_\Gamma \zeta^m) + \dfrac{2\sigma}{R^2} (u_\eps^m, \zeta^m) \\
        &+ \kappa( \nabla_\Gamma(\boldsymbol\Lambda\cdot\vphi_\eps^m), \nabla_\Gamma \zeta^m) - \dfrac{2\kappa}{R^2} (\Lambda \cdot\vphi_\eps^m, \zeta^m),
    \end{align}
for all $\peta^m, \pxi^m \in \boldsymbol\VV^m$ and $\zeta^m\in V^m$, in $(0,T)$, and subject to the initial conditions
    \begin{align*}
        \vphi_\eps^m(0) = \boldsymbol{\mathcal P}_m (\vphi_0-\al) \quad \text{ and } \quad u_\eps^m(0) = \mathcal P_m u_0.
    \end{align*}
\end{prob}

Problem \ref{5_gal} can be viewed as a Cauchy problem for a system of ODEs for the functions $\alpha_k^i$, $\gamma_k$, with $i=1,\dots,N$ and $k=1,\dots,m$. This problem has a unique (local) solution by Cauchy-Lipschitz theorem. Moreover, thanks to the regularity of $\boldsymbol\Psi_\eps$, $\alpha_k^i,\gamma_k\in C^2([0,t_n))$ for $i=1,\ldots,N,\ k=1,\ldots,m$ for some $t_n>0$, this also entails that $\beta_k^i\in C^2([0,t_n))$ for $i=1,\ldots,N,\ k=1,\ldots,m$. Note that by construction $\vphi_\eps^m,\mathbf w_\eps-\dashint_\Gamma{\mathbf w}_\eps^m\in \boldsymbol{\mathcal V}_m$.

We now establish uniform bounds with respect to $m$ on this approximating solution in order to pass to the limit and recover the wanted weak solution.
These bounds follow from a number of technical lemmas on the approximated energy.

First of all, we prove that the approximated energy is bounded from below.

\begin{lem}
\label{5_eps}
There exists $\eps_0=\eps_0(\kappa,\vert\Lambda\vert,b,\varepsilon)>0$ such that, for $0<\eps\leq\eps_0(\kappa,\vert\Lambda\vert,b,\varepsilon)$, there holds
\begin{align*}
    \mathcal E_\eps(u,\vphi) \geq K_1 \|u\|_{H^2(\Gamma)}^2 + K_2 \|\tgrad\vphi\|^2 - K_3, \quad \forall (\vphi,u)\in \mathcal{K},
\end{align*}
where $K_1, K_2, K_3>0$ are independent of $\eps$.
\end{lem}

\begin{proof}
Given $u\in\mathcal{K}_2$, the following Poincar\'e-type inequality holds (see, e.g., \cite{EllFriHob17-a}):
\begin{align}
    \int_\Gamma u^2 \d\Gamma \leq \dfrac{R^2}{6}\int_\Gamma |\tgrad u|^2  \d\Gamma \leq \dfrac{R^4}{36}\int_\Gamma (\Delta_\Gamma u)^2 \d\Gamma.
    \label{5_poin}
\end{align}
This allows us to obtain
\begin{align*}
    \mathcal E_\eps(u,\vphi)
    &\geq \int_{\Gamma} \dfrac{\kappa}{3} (\Delta_{\Gamma} u)^2 + \dfrac{2\sigma u^2}{R^2} + \kappa (\vlambda \cdot \vphi) \Delta_{\Gamma}u + \dfrac{2\kappa u \vlambda \cdot \vphi}{R^2} +\dfrac{\kappa (\vlambda\cdot\vphi)^2}{2}\d\Gamma + \mathrm E_{\mathrm{CH},\eps}(\vphi).
\end{align*}
Applying Young's inequality
\begin{align*}
    \kappa \vlambda\cdot\vphi \Delta_\Gamma u \geq -\dfrac{\kappa}{6}(\Delta_\Gamma u)^2 - \dfrac{3\kappa}{2} (\vlambda\cdot\vphi)^2,
\end{align*}
\Andrea{and using the definition of $\mathrm E_{\mathrm{CH}, \eps}$}, we get
\begin{align*}
    \mathcal E_\eps(u,\vphi) &\geq \int_{\Gamma} \dfrac{\kappa}{6} (\Delta_{\Gamma} u)^2 + \dfrac{2\sigma u^2}{R^2} - \kappa (\vlambda\cdot\vphi)^2 + \dfrac{2\kappa u \vlambda \cdot \vphi}{R^2} + \dfrac{b\varepsilon}{2}|\nabla_{\Gamma} \vphi|^2 + \dfrac{b}{\varepsilon} \boldsymbol\Psi_\eps(\vphi) \d\Gamma.
\end{align*}
Similarly, from
\begin{align*}
    \dfrac{2\kappa (\vlambda\cdot\vphi) u}{R^2} \geq -\dfrac{\kappa}{3} (\vlambda\cdot\vphi)^2 - \dfrac{3\kappa u^2}{R^4}
\end{align*}
combined with
\begin{align*}
    \int_\Gamma -\dfrac{3\kappa u^2}{R^4}\d\Gamma \geq \int_\Gamma -\dfrac{\kappa}{12} (\Delta_\Gamma u)^2 \d\Gamma,
\end{align*}
we deduce
\begin{align*}
    \mathcal E_\eps(u,\vphi) &\geq \int_{\Gamma} \dfrac{\kappa}{12} (\Delta_{\Gamma} u)^2 + \dfrac{2\sigma u^2}{R^2} - \dfrac{4\kappa (\vlambda\cdot\vphi)^2}{3} + \dfrac{b\varepsilon}{2}|\nabla_{\Gamma} \vphi|^2 + \dfrac{b}{\varepsilon} \boldsymbol\Psi_\eps(\vphi)\d\Gamma \\
    &\geq C_1 \| u\|_{H^2(\Gamma)}^2 + \int_\Gamma  \dfrac{b\varepsilon}{2}|\nabla_{\Gamma} \vphi|^2 + \dfrac{b}{\varepsilon} \boldsymbol\Psi_\eps(\vphi)- \dfrac{4\kappa (\vlambda\cdot\vphi)^2}{3}\d\Gamma \\
    &\geq C_1 \| u\|_{H^2(\Gamma)}^2 + \int_\Gamma  \dfrac{b\varepsilon}{2}|\nabla_{\Gamma} \vphi|^2 + \dfrac{b}{\varepsilon} \boldsymbol\Psi_\eps(\vphi)- \dfrac{4\kappa|\vlambda|^2 |\vphi|^2}{3}\d\Gamma.
\end{align*}
 We now apply \eqref{coercive} with $K_0=4\varepsilon\kappa|\vlambda|^2/(3b)$ to get the existence of $\eps_0$ depending on $K_0$ and of a constant $C_2>0$ (depending on $\eps_0$) such that
$$
\dfrac{b}{\varepsilon} \boldsymbol\Psi_\eps(\vphi)- \dfrac{4\kappa|\vlambda|^2 |\vphi|^2}{3}\geq -C_2,\quad \forall 0<\eps\leq \eps_0,
$$
which implies
\begin{align*}
    \mathcal{E}_\eps(u, \vphi) \geq K_1 \| u\|_{H^2(\Gamma)}^2 + K_2\Vert\nabla_\Gamma\vphi\Vert^2-K_3,\quad \forall 0<\eps\leq \eps_0,
\end{align*}
for $K_1,K_2,K_3>0$ constants independent of the specific $\eps$. This concludes the proof.
\end{proof}

The following lemma gives a bound from above for the approximated energy.

\begin{lem}\label{5_lem:initial_energy_estimate}
There exists $C>0$ such that, for any $\eps>0$, there exists $M(\eps)>0$ such that
 \begin{align*}
     \mathcal{E}_\eps (\boldsymbol{\mathcal P}_m(\vphi_0-\al)+\al, \mathcal P_m u_0) \leq C,\qquad \forall m\geq M(\eps).
 \end{align*}
\end{lem}

\begin{proof}
The only part which is not immediate from properties \eqref{proj}-\eqref{proj2} is the control of the terms involving $\boldsymbol\Psi_\eps(\vphi_0)$, which can be obtained by using the convexity of $\boldsymbol\Psi_\eps^1$.
Indeed, we have
\begin{align*}
	\int_\Gamma \boldsymbol\Psi_\eps^1(\boldsymbol{\mathcal P}_m(\vphi_0-\al)+\al)\d\Gamma &\leq \int_\Gamma \boldsymbol\Psi_\eps^1(\vphi_0)\d\Gamma + \int_\Gamma \boldsymbol\psi_\eps'(\boldsymbol{\mathcal P}_m(\vphi_0-\al)+\al)\cdot(\Pb_m (\vphi_0-\al)+\al - \vphi_0)\d\Gamma 
\\&\leq \int_\Gamma \boldsymbol\Psi^1(\vphi_0)\d\Gamma  + (C+C(\eps))\Vert\Pb_m (\vphi_0-\al) +\al- \vphi_0\Vert^2, 
\end{align*}
where we used $0\leq\vphi_0 \leq 1$ and properties (i)-(iv) of $\psi_\eps$. For any fixed $\eps>0$, the last term converges to $0$ because $\Pb_m (\vphi_0-\al)\to\vphi_0-\al$ as $m\to\infty$ strongly in both $\HHH$ and $\VVV$.
This leads to the existence of a constant $\textstyle{C=\int_\Gamma \boldsymbol\Psi^1(\vphi_0)+1}$ such that, for any fixed $\eps>0$, there exists $M(\eps)>0$ so that
\begin{align*}
	\int_\Gamma \boldsymbol\Psi_\eps(\Pb_m (\vphi_0-\al)+\al)\d\Gamma \leq C, \quad \forall m\geq M(\eps),
\end{align*}
from which the conclusion follows.
\end{proof}

Let us then fix $0<\eps\leq \eps_0$ ($\eps_0$ given in Lemma \ref{5_eps}) and $m>M(\eps)$ given in Lemma \ref{5_lem:initial_energy_estimate}. We can prove the validity of the following approximated
energy inequality deriving from the energy identity. More precisely, it holds

\begin{lem}\label{5_lem:energy_estimate}
We have the energy identity
\begin{align}\label{5_eq:derivative_energy}
    \dfrac{d}{dt} \mathcal E_\eps( u_\eps^m, \vphi_\eps^m+\al) = -\|\partial_t u_\eps^m\|^2 - \|\tgrad \mathbf{w}_\eps^m\|^2.
\end{align}
Consequently, there exists $C>0$, independent of $\eps$ and $m$, such that
\begin{align}\label{5_eq:energy_estimate}
    \sup_{t\in [0,T]} \mathcal E_\eps( u_\eps^m(t),\vphi_\eps^m(t)+\al) + \int_0^T \|\partial_t u_\eps^m\|^2 + \int_0^T \|\tgrad \mathbf{w}_\eps^m\|^2 \leq C.
\end{align}
\end{lem}

\begin{proof}
Energy identity \eqref{5_eq:derivative_energy} holds due to the gradient flow structure of the problem. As a consequence, integrating in time over an interval $(0,t)$, we find
\begin{align*}
     \mathcal E_\eps( u_\eps^m,\vphi_\eps^m+\al) + \int_0^t \|\partial_t u_\eps^m\|^2 + \int_0^t \|\tgrad \mathbf{w}_\delta^m\|^2 = \mathcal E_\eps(\mathcal P_m u_0,\boldsymbol{\mathcal P}_m (\vphi_0-\al)+\al) \leq C,
\end{align*}
where the inequality follows from Lemma \ref{5_lem:initial_energy_estimate}. The above inequality, together with Lemma \ref{5_eps} entails \eqref{5_eq:energy_estimate}.
\end{proof}

From the above lemmas we infer the first uniform bounds which allow to extend the approximating solution to the time interval $[0,+\infty)$. We now state the following
\begin{cor}
\label{5_lem:estimates1}
There exists $C>0$, independent of $\eps$ and $m$, such that:
 \begin{itemize}
     \item[(i)] we have
     \begin{align*}
     \|u_\eps^m\|_{L^\infty(0,T; H^2(\Gamma))} + \|\partial_t u_\eps^m\|_{L^2(0,T; H)} + \|\nabla_\Gamma \vphi_\eps^m\|_{L^\infty(0,T; \HHH)} + \|\nabla _\Gamma \mathbf w_\eps\|_{L^2(0,T; \HHH)}\leq C
     \end{align*}
     and
     \begin{align*}
         \sup_{t\in [0,T]} \int_\Gamma \boldsymbol\Psi_\eps(\vphi_\eps^m+\al)\d\Gamma \leq C;
    \end{align*}
     \item[(ii)] $ \|\vphi_\eps^m\|_{L^\infty(0,T; \VVV)} \leq C$;
    \item[(iii)] $\|\partial_t \vphi_\eps^m\|_{L^2(0,T; \VVV')} \leq C$.
    \end{itemize}
\end{cor}

\begin{proof}
(i). The second statement is an immediate consequence of the energy estimate \eqref{5_eq:energy_estimate}. To obtain the first combine \eqref{5_eq:energy_estimate} with  Lemma \ref{5_eps}. \medskip

(ii). This follows immediately from Poincar\'{e}'s inequality and the fact that the total mass is conserved over time. \medskip

(iii). Note that, for any $\peta\in L^2(0,T; \VVV)$, the first equation of Problem \ref{5_gal} yields
\begin{align*}
     \left| \int_0^T \int_\Gamma (\partial_t \vphi_\eps^m, \peta) \right|
    &= \left| \int_0^T \int_\Gamma (\partial_t \vphi_\eps^m, \boldsymbol{\mathcal P}_m \peta) \d\Gamma\right|   \\
    &= \left|\int_0^T \int_\Gamma \mathbf L\nabla_\Gamma \mathbf w_\eps^m \colon \nabla_\Gamma \boldsymbol{\mathcal P}_m\peta\d\Gamma \right| \\
    &\leq \|\mathbf L \|_\infty \, \|\nabla_\Gamma \mathbf w_\eps^m \|_{L^2(0,T; \HHH)} \, \| \boldsymbol{\mathcal P}_m \peta \|_{L^2(0,T; \VVV)}\\
    &\leq \| \mathbf L \|_\infty \, \|\nabla_\Gamma \mathbf w_\eps^m \|_{L^2(0,T; \HHH)} \, \| \mathcal \peta \|_{L^2(0,T; \VVV)}.
\end{align*}
Indeed, $\partial_t\vphi^m_\eps=\boldsymbol{\mathcal P}_m\partial_t\vphi^m_\eps$ and $\boldsymbol{\mathcal P}_m$ is selfadjoint.
Then the first estimate on $\nabla_\Gamma \mathbf w_\eps^m$ implies
\begin{align}\label{5_eq:third_estimate}
    \|\partial_t \vphi_\eps^m\|_{L^2(0,T; \VVV')} \leq C.
\end{align}
The proof is finished.
\end{proof}
Exploiting the previous lemmas, we can also prove the following time-weighted estimates:

\begin{lem}
\label{5_lem:estimates_weighteds}
There exists $C>0$, independent of $\eps$ and $m$, such that:
\begin{itemize}
     \item[(i)] we have
     \begin{align}
     &\|\sqrt{t}\tgrad\mathbf w_\eps^m\|_{L^\infty(0,T;\HHH)} + \|\sqrt{t}\tgrad\partial_t \vphi_\eps^m\|_{L^2(0,T;\HHH)} \notag \\
     \label{wb}
     &+ \|\sqrt{t}\partial_t u_\eps^m\|_{L^\infty(0,T;H)} +\|\sqrt{t}\Delta _\Gamma \partial_t u_\eps\|_{L^2(0,T;H)}\leq C;
     \end{align}
  \item[(ii)] moreover it holds $$\|\sqrt{t} \partial_t u_\eps\|_{L^2(0,T; H^2(\Gamma))}+ \|\sqrt{t}\partial_t \vphi_\eps^m\|_{L^2(0,T; \VVV)}\leq C.$$
    \end{itemize}
\end{lem}

\begin{proof}
(i). Let us choose $\peta^m=\partial_t(\mathbf w_\eps^m-\dashint_\Gamma{\mathbf w}_\eps^m)\in \boldsymbol\VV_m$ in \eqref{e1}, and then differentiate with respect to time equation \eqref{e3}, choosing then $\xi_m=\partial_t u^m_\eps$ (this is possible due to the regularity of the Galerkin approximation). Adding together the obtained identities and recalling that $\partial_t\vphi_\eps^m\in T\Sigma$, we find
\begin{align*}
    &\dfrac{1}{2}\dfrac{d}{dt}\left( (\textbf{L}\tgrad\mathbf w_\eps^m,\tgrad \mathbf w_\eps^m)+\beta\Vert\partial_t u_\eps^m\Vert^2\right)+\frac{b}{\varepsilon}\sum_{i=1}^N(\psi^{\prime\prime}_\eps(\vphi_{\eps,i}^m+\al_i),\vert\partial_t\vphi_{\eps,i}^m\vert^2)\\
    &\hskip1cm -\frac{b}{\varepsilon}(A\partial_t\vphi^m_\eps,\partial_t\vphi_\eps^m)+\frac{2\kappa}{R^2}(\partial_t u_\eps^m\boldsymbol\Lambda, \partial_t\vphi_\eps^m)+b\varepsilon\Vert\tgrad\partial_t\vphi_\eps^m\Vert^2+\kappa\Vert\boldsymbol\Lambda\cdot\partial_t\vphi_\eps^m\Vert^2\\
    &\hskip1cm+\kappa(\boldsymbol\Lambda\Delta_\Gamma\partial_t u_\eps^m,\partial_t\vphi_\eps^m)+\Vert\Delta_\Gamma\partial_t u_\eps^m\Vert^2+\left(\sigma-\frac{2\kappa}{R^2}\right)\Vert\tgrad\partial_t u_\eps^m\Vert^2\\
    &\hskip1cm-\frac{2\sigma}{R^2}\Vert\partial_t u_\eps^m\Vert^2-\kappa(\boldsymbol\Lambda\tgrad\partial_t\vphi_\eps^m,\nabla_\Gamma\partial_t u_\eps^m)+\frac{2\kappa}{R^2}(\boldsymbol\Lambda\cdot\partial_t\vphi_\eps^m,\partial_t u_\eps^m)=0.
\end{align*}
Observe that $\sum_{i=1}^N(\psi_\eps^{\prime\prime}(\vphi_{\eps,i}^m+\al_i),\vert\partial_t\vphi_{\eps,i}^m\vert^2)\geq 0$ since $\psi_\eps$ is convex. Moreover, we have
\begin{align}
\nonumber(\partial_t\vphi^m_\eps,\partial_t\vphi_\eps^m)=(\textbf{L}\tgrad\mathbf w_\eps^m,\tgrad\partial_t\vphi_\eps^m)&\nonumber\leq C\Vert\tgrad\mathbf w_\eps^m\Vert\,\Vert\tgrad\partial_t\vphi_\eps^m\Vert\\
&\leq C\Vert\tgrad\mathbf w_\eps^m\Vert^2+\frac{b\varepsilon}{4}\Vert\tgrad\partial_t\vphi_\eps^m\Vert^2.\label{interp}
\end{align}
On the other hand, via interpolation, we get
\begin{align}
\Vert \tgrad\partial_t u_\eps^m\Vert^2\leq C\Vert \partial_t u_\eps^m\Vert\Vert \partial_t u_\eps^m\Vert_{H^2(\Gamma)}\leq C\Vert \partial_t u_\eps^m\Vert\Vert \Delta_\Gamma\partial_t u_\eps^m\Vert,
    \label{5_interp}
\end{align}
where we used the fact that $\partial_t u_\varepsilon^m\in \mathcal{K}_2$ and thus \eqref{5_poin} applies.
Using now standard inequalities,  \eqref{5_interp}, and $\Vert \partial_t \vphi^m_\eps\Vert\leq C\Vert \tgrad\partial_t \vphi^m_\eps\Vert$, we obtain
\begin{align*}
    &-\frac{2\kappa}{R^2}(\partial_t u_\eps^m\boldsymbol\vlambda, \partial_t\vphi_\eps^m)-\kappa(\boldsymbol\vlambda\Delta_\Gamma\partial_t u_\eps^m,\partial_t\vphi_\eps^m)-\left(\sigma-\frac{2\kappa}{R^2}\right)\Vert\tgrad\partial_t u_\eps^m\Vert^2\\
    &+\frac{2\sigma}{R^2}\Vert\partial_t u_\eps^m\Vert^2+\kappa(\Lambda\tgrad\partial_t\vphi_\eps^m,\nabla_\Gamma\partial_t u_\eps^m)-\frac{2\kappa}{R^2}(\boldsymbol\Lambda\cdot\partial_t\vphi_\eps^m,\partial_tu_\eps^m)\\&\leq \dfrac{b\varepsilon}{4}\Vert\tgrad\partial_t\vphi_\eps^m\Vert^2+\frac{1}{2}\Vert\Delta_\Gamma\partial_t u_\eps^m\Vert^2+C\Vert\partial_t u_\eps^m\Vert^2.
    \end{align*}
    Therefore we find in the end, using also \eqref{interp},
    \begin{align*}
    \dfrac{1}{2}\dfrac{d}{dt}\left( (\textbf{L}\tgrad\mathbf w_\eps^m,\tgrad \mathbf w_\eps^m)+\beta\Vert\partial_t u_\eps^m\Vert^2\right)+\frac{b\varepsilon}{4}\Vert\tgrad\partial_t\vphi_\eps^m\Vert^2+&\frac{1}{2}\Vert\Delta_\Gamma\partial_t u_\eps^m\Vert^2\\
    &\leq C\left(\Vert\tgrad\mathbf w_\eps^m\Vert^2+\beta\Vert\partial_t u_\eps^m\Vert^2\right).
\end{align*}
We now multiply the above inequality by $s\in (0,T)$. This gives
   \begin{align*}
    \dfrac{1}{2}\dfrac{d}{dt}\left( s(\textbf{L}\tgrad\mathbf w_\eps^m,\tgrad \mathbf w_\eps^m)+\beta s\Vert\partial_t u_\eps^m\Vert^2\right)+&\frac{b\varepsilon s}{4}\Vert\tgrad\partial_t\vphi_\eps^m\Vert^2+\frac{s}{2}\Vert\Delta_\Gamma\partial_t u_\eps^m\Vert^2\\
    &\leq C(1+s)\left(\Vert\tgrad\mathbf w_\eps^m\Vert^2+\beta\Vert\partial_t u_\eps^m\Vert^2\right).
\end{align*}
Let us now integrate over $(0,t)$ and recall the fact that $\textbf{L}$ is strictly positive definite over $T\Sigma$, to obtain
  \begin{align*}
    t\left( C\Vert\tgrad \mathbf w_\eps^m\Vert^2+\beta \Vert\partial_t u_\eps^m\Vert^2\right)+\frac{b\varepsilon }{4}\int_0^t &s\Vert\tgrad\partial_t\vphi_\eps^m\Vert^2ds+\frac{1}{2}\int_0^t s\Vert\Delta_\Gamma\partial_t u_\eps^m\Vert^2ds\\
    &\leq C\int_0^t(1+s)\left(\Vert\tgrad\mathbf w_\eps^m\Vert^2+\beta\Vert\partial_t u_\eps^m\Vert^2\right)ds.
\end{align*}
The integral Gronwall's Lemma, together with Corollary \ref{5_lem:estimates1}, gives \eqref{wb}. \medskip

(ii). The result easily follows by recalling that $\partial_t u_\eps^m\in \mathcal{K}_2$ and $\int_\Gamma\partial_t\vphi_\eps^m\equiv 0$, via standard Poincar\'e's inequality (see \eqref{5_poin}) .
\end{proof}

Additional bounds, for $\eps>0$ fixed, are given by

\begin{lem}
\label{epps}
Let $\eps\in(0,\eps_0)$ be fixed. There exists $C_\eps>0$, independent of $m$ but possibly depending on $\eps$, such that
 $$\| \boldsymbol\psi_{\eps}'(\vphi_\eps^m+\al)\|_{L^2(0,T; \HHH)} \leq C_\eps\quad \text{ and } \quad \|\mathbf w_\eps^m\|_{L^2(0,T; \HHH)}\leq C_\eps .$$
\end{lem}

\begin{proof}
Noticing that there exist $C_{1,\eps}, C_{2,\eps}>0$ such that
\begin{align*}
    |\boldsymbol\psi_{\eps}'(r)| \leq C_{1,\eps} |r| + C_{2,\eps}
\end{align*}
immediately implies the estimate for $\boldsymbol\psi_{\eps}'(\vphi_\eps^m)$. As for $\mathbf w_\eps^m$, we can control its mean value by
\begin{align*}
    \left\vert \dashint_\Gamma{\mathbf w}_\eps^m\right\vert&\leq C \left( \| \mathbf P\boldsymbol\psi_{\eps}'(\vphi_\eps^m+\al)\| + \|\vphi_{\eps}^m\|_{\VVV} + \|u_\eps^m\|_{H^2(\Gamma)} \right) \\
    &\leq C \left(1+ \|\boldsymbol\psi_{\eps}'(\vphi_\eps^m+\al)\| + \|\vphi_{\eps}^m\|_{\VVV} + \|u_\eps^m\|_{H^2(\Gamma)} \right),
\end{align*}
and the conclusion follows from the previous a priori estimates.
\end{proof}

We now have all the ingredients to pass to the limit as $m$ goes to $\infty$ along a suitable subsequence.

\paragraph{Passage to the limit $m\to\infty$.} By standard compactness arguments, the previous bounds imply the existence of $(\vphi_\eps, \mathbf w_\eps, u_\eps)$, defined on $[0,+\infty)$, which solves Problem \ref{5_prob:regularised} and satisfies, for some $C>0$ independent of $\eps$, the energy estimate
\begin{align*}
    \sup_{t\in [0,T]} \mathcal E_\eps(\vphi_\eps(t), u_\eps(t)) + \int_0^T \|\partial_t u_\eps\|^2 + \int_0^T \|\tgrad \mathbf{w}_\eps\|^2 \leq C
\end{align*}
as well as the estimates
\begin{align}
    \|u_\eps\|_{L^\infty(0,T;H^2(\Gamma))} + \|\partial_t u_\eps\|_{L^2(0,T;H)} &\leq C, \label{5_eq:first_estimatec} \\
    \|\vphi_\eps\|_{L^\infty(0,T;\VVV)} + \|\partial_t\vphi_\eps\|_{L^2(0,T;\VVV')} +\Vert\boldsymbol\Psi_\eps(\vphi_\eps)\Vert_{L^\infty(0,T;\LLL^1(\Gamma))}&\leq C, \label{5_eq:first_estimatea} \\
    \|\nabla _\Gamma \mathbf w_\eps\|_{L^2(0,T;\HHH)} &\leq C, \\
    \|\sqrt{t}\tgrad\mathbf w_\eps\|_{L^\infty(0,T; \HHH)} + \|\sqrt{t}\partial_t \vphi_\eps\|_{L^2(0,T; \VVV)} + \|\sqrt{t}\partial_t u_\eps\|_{L^\infty(0,T; H^2(\Gamma))}&\leq C.\label{5_eq:first_estimateb}
\end{align}
To be precise, we find that $\vphi_{\eps}^m$ converges in suitable norms to a function $\widetilde{\vphi}_{\eps}(t)\in \VVV_0$ (for almost any $t\geq0$) as $m\to \infty$ up to a subsequence. We then define $\vphi_\eps:=\widetilde{\vphi}_\eps+\al$ to obtain the results above. We recall that it is now easy to show:
\begin{align*}
	\sum_{i=1}^N\vphi_{\eps,i}\equiv 1,\quad \dashint_\Gamma\vphi_\eps=\al=\dashint_\Gamma\vphi_0,\quad \sum_{i=1}^N\mathbf w_{\eps,i}\equiv 0.
\end{align*}
Note that, by elliptic regularity in the definition of $\mathbf w_\eps$, recalling Lemma \ref{epps}, which is valid also as $m\to\infty$ by lower semicontinuity of the norms involved, we have
\begin{align}
	\norm{\vphi_\eps}_{L^2(0,T;\mathbf H^2(\Omega))}\leq C_\eps.
	\label{H2}
	\end{align}

In order to let $\eps\to 0$ we need to improve the estimates in Lemma \ref{epps}, which are not independent of $\eps$. We do so in the next lemma.

\begin{lem}\label{5_lem:estimates2}

 There exists $C>0$, independent of $\eps$, such that:
 \begin{itemize}
     \item[(i)] $\|\mathbf w_\eps\|_{L^2(0,T; \VVV)} \leq C$;
     \item[(ii)] $\|\boldsymbol\psi_{\eps}'(\vphi_\eps)\|_{L^2(0,T; \HHH)} \leq C$.
 \end{itemize}
\end{lem}

\begin{proof}
(i). We use here an idea developed in \cite{Gar05} and adapted in \cite{GGPS,AC2023}: observe that
\begin{align*}
    &\boldsymbol\lambda_\eps := \dashint_\Gamma \mathbf w_\eps= \dashint_\Gamma \left(\dfrac{b}{\varepsilon} \mathbf P \boldsymbol\psi_{\eps}'(\vphi_\eps)-\dfrac{b}{\varepsilon} \mathbf P \boldsymbol\Psi_{,\vphi}^2(\vphi_\eps) + \dfrac{2\kappa}{R^2}\mathbf P u_\eps \mathbf \Lambda + \kappa \mathbf P(\mathbf \Lambda \cdot\vphi_\eps) \mathbf \Lambda\right)\\&=\dashint_\Gamma \left(\dfrac{b}{\varepsilon} \mathbf P \boldsymbol\psi_\eps'(\vphi_\eps)-\dfrac{b}{\varepsilon}\mathbf PA\vphi_\eps +  \kappa \mathbf P(\mathbf \Lambda \cdot\vphi_\eps) \mathbf \Lambda\right),
\end{align*}
where we used (see Remark \ref{5_constraint2})
$$
\left(\dashint_\Gamma \dfrac{2\kappa}{R^2}\mathbf P u_\eps \mathbf \Lambda\right)_i=\dfrac{2\kappa}{R^2}\left(\mathbf\Lambda_i-\dfrac{\sum_{j=1}^N\mathbf\Lambda_j}{N}\right)\dashint_\Gamma u_\eps=0.
$$
Set now $\widehat{\mathbf w}_\eps = \mathbf w_\eps - \boldsymbol\lambda_\eps$. The second equation of the system is then equivalent to
\begin{align}
\begin{split}
\label{5_eq:second_equation_hat}
    (\widehat{\mathbf w}_\eps + \boldsymbol\lambda_\eps, \pxi) = \dfrac{b}{\varepsilon} (\mathbf P\boldsymbol\psi_{\eps}'(\vphi_\eps), \pxi) - \dfrac{b}{\varepsilon}(\mathbf PA\vphi_\eps, \pxi) + \dfrac{2\kappa}{R^2} (\mathbf Pu_\eps\boldsymbol\Lambda, \pxi) + b\varepsilon (\nabla_\Gamma \vphi_\eps, \nabla_\Gamma \mathbf \pxi) \\
        + \kappa( \mathbf P(\boldsymbol\Lambda\cdot\vphi_\eps)\boldsymbol\Lambda, \pxi) - \kappa (\nabla_\Gamma u_\eps, \nabla_\Gamma (\mathbf P\boldsymbol\Lambda\cdot \pxi)).
\end{split}
\end{align}
Let $\mathbf k\in \mathbf G$ ($\mathbf G$ is the Gibbs simplex) be a constant in space, but possibly depending on time. Since $\vphi_\eps\in  \Sigma$, we have $\mathbf k - \vphi_\eps\in T\Sigma$, and thus by convexity of $\boldsymbol\Psi_\eps^1$ it holds
\begin{align*}
    \int_\Gamma \boldsymbol\Psi_\eps^1(\mathbf k)\d\Gamma &\geq \int_\Gamma \boldsymbol\Psi_\eps^1(\vphi_\eps) \d\Gamma+ \int_\Gamma \boldsymbol\psi_{\eps}'(\vphi_\eps)\cdot(\mathbf k - \vphi_\eps)\d\Gamma \\
    &= \int_\Gamma \boldsymbol\Psi_\eps^1(\vphi_\eps)\d\Gamma + \int_\Gamma \mathbf P\boldsymbol\psi_{\eps}'(\vphi_\eps)\cdot(\mathbf k - \vphi_\eps)\d\Gamma,
\end{align*}
on account of the identity $\mathbf k - \vphi_\eps =\textbf{P}(\mathbf k - \vphi_\eps)$ and recalling that $\textbf{P}$ is selfadjoint.
Using \eqref{5_eq:second_equation_hat} with $\pxi= \mathbf k - \vphi_\eps$ and property (i) of $\psi _{\varepsilon }$, we deduce
\begin{align}
\begin{split}
\label{5_eq:F}
   C\geq\int_\Gamma \boldsymbol\Psi_\eps^1(\mathbf k)\d\Gamma \geq \int_\Gamma \boldsymbol\Psi_\eps^1(\vphi_\eps)\d\Gamma + \dfrac{b}{\varepsilon}\int_\Gamma \left(\widehat{\mathbf w}_\eps + \boldsymbol\lambda_\eps\right)\cdot (\mathbf k - \vphi_\eps) \d\Gamma+ \int_\Gamma &\mathbf PA\vphi_\eps\cdot(\mathbf k - \vphi_\eps)\d\Gamma \\
    &+
    F(u_\eps, \vphi_\eps, \mathbf k - \vphi_\eps)
\end{split}
\end{align}
where
\begin{align*}
    F:=-\dfrac{b}{\varepsilon}(\dfrac{2\kappa}{R^2} (\mathbf Pu_\eps\boldsymbol\Lambda, \mathbf k - \vphi_\eps) + b\varepsilon (\nabla_\Gamma \vphi_\eps, \nabla_\Gamma \mathbf (\mathbf k - \vphi_\eps)) &+ \kappa( \mathbf P(\boldsymbol\Lambda\cdot\vphi_\eps)\boldsymbol\Lambda, \mathbf k - \vphi_\eps)  \\
         &- \kappa (\nabla_\Gamma u_\eps, \nabla_\Gamma(\mathbf P\boldsymbol\Lambda\cdot (\mathbf k - \vphi_\eps)))).
\end{align*}
Observe that, recalling the previous lemma and using the Cauchy-Schwarz inequality, we get
\begin{align*}
    |  F  | &\leq C\left( \|u_\eps\|_{V} + \|\vphi_\eps\|_{\VVV}\right) \|\mathbf k - \vphi_\eps\|_{\VVV}\leq C.
\end{align*}
We aim to obtain an estimate on
\begin{align*}
    \int_\Gamma \boldsymbol\lambda_\eps \cdot (\mathbf k - \vphi_\eps)\d\Gamma.
\end{align*}
Rearranging the terms in \eqref{5_eq:F}, using the estimate for $ F$, Poincar\'{e}'s inequality for the term involving $\widehat{\mathbf w}_\eps$ (which has null mean)
and the uniform bounds stated in the previous lemma we obtain
\begin{align*}
    \int_\Gamma \boldsymbol\lambda_\eps \cdot (\mathbf k - \vphi_\eps)\d\Gamma &\leq  \|\boldsymbol\Psi_\eps^1(\vphi_\eps)\|_{\LLL^1(\Gamma)} + \left(\| \widehat{\mathbf w}_\eps\| + \|\mathbf PA\vphi_\eps\| \right) \left( \|\mathbf k\| + \|\vphi_\eps\|\right) + | \mathbf F(u_\eps, \vphi_\eps, \mathbf k - \vphi_\eps) |\\
    &\leq C\left(1+\| \nabla_\Gamma \mathbf w_\eps\|\right).
\end{align*}
Observe now that there exists $\rho\in (0,1)$ such that, for all $i=1,\dots,N$ and almost all $t\in (0,T)$,
\begin{align*}
    \rho <\dashint_\Gamma \vphi_{i,\eps} < 1-\rho,
\end{align*}
which is a consequence of the conservation of mass and the fact that $\dashint_\Gamma \vphi_0=\al\in(0,1)^N$, and $\al\in \Sigma$.
Then, choose, for each $k,j=1,\dots,N$,
\begin{align*}
    \mathbf k = \mathbf k_{k,j} := \dashint_\Gamma \vphi_\eps(t) + \rho \,\text{sign }(\lambda_{k,\eps} - \lambda_{j,\eps})(\mathbf e_k - \mathbf e_j)\in \mathbf G,
\end{align*}
where $\mathbf e_j$ is the $j$-th element of the canonical basis of $\R^N$.
As a consequence, we have
\begin{align*}
    \int_\Gamma \boldsymbol\lambda_\eps \cdot (\mathbf k_{k,j} - \vphi_\eps)\d\Gamma \leq C\left( 1 + \| \nabla_\Gamma \mathbf w_\eps\| \right),
\end{align*}
with $C$ independent of $\eps$, $k$ and $j$. On the other hand, we note that
\begin{align*}
    \int_\Gamma \boldsymbol\lambda_\eps \cdot (\mathbf k_{k,j} - \vphi_\eps) \d\Gamma= \rho |\Gamma| \, |\lambda_{k,\eps} - \lambda_{j,\eps}| + \int_\Gamma \boldsymbol\lambda_\eps \cdot \left( \vphi_\eps - \dashint_\Gamma \vphi_\eps\right)\d\Gamma = \rho |\Gamma| \, |\lambda_{k,\eps} - \lambda_{j,{\eps}}|,
\end{align*}
and hence expanding the left hand side we are led to
\begin{align*}
    |\lambda_{k,\eps} - \lambda_{j,\eps}| \leq \dfrac{ C}{\rho |\Gamma|}\left( 1 + \| \nabla_\Gamma \mathbf w_\eps\| \right).
\end{align*}
Using now the identity
\begin{align*}
    \lambda_{k,\eps} = \dfrac{1}{N} \sum_{j=1}^N (\lambda_{k,\eps} - \lambda_{j,\eps}),
\end{align*}
which is valid since $\boldsymbol\lambda_\eps\in T\Sigma$, we find, for almost any $t\in(0,T)$, that
\begin{align}
\label{5_lam}
    |\boldsymbol\lambda_\eps(t)|^2 = \sum_{k=1}^N|\boldsymbol\lambda_{k,\eps}(t)|^2
    \leq \dfrac{1}{N}\sum_{k=1}^N\sum_{j=1}^N \left| \lambda_{k,\eps}(t) - \lambda_{j,\eps}(t)\right|^2 \leq C_1 + C_2 \|\nabla_\Gamma \mathbf w_\eps(t)\|^2.
\end{align}
Therefore
\begin{align*}
    \int_0^T |\boldsymbol\lambda_\eps|^2 \d\Gamma\leq C.
\end{align*}
Combining the above bound with Poincar\'{e}'s inequality and the previous bounds, we are led to
\begin{align*}
    \int_0^T \int_\Gamma |\mathbf w_\eps|^2\d\Gamma \leq 2 \int_0^T \int_\Gamma |\mathbf w_\eps - \boldsymbol\lambda_\eps|^2\d\Gamma + 2 \int_0^T \int_\Gamma  |\boldsymbol\lambda_\eps|^2\d\Gamma \leq 2C_P \int_0^T \int_\Gamma |\nabla_\Gamma \mathbf w_\eps|^2\d\Gamma + C \leq \tilde C,
\end{align*}
as desired.
\vskip 3mm

(ii). Set $ \boldsymbol\psi_\eps'(\vphi_\eps)=(\psi'_\eps(\vphi_{\eps,i}))_{i=1,\ldots,N}$ and note that
\begin{align*}
    \nabla_\Gamma \psi'_\eps(\vphi_{\eps,i}) = \psi''_\eps(\vphi_{\eps,i}) \nabla_\Gamma \vphi_{\eps,i} \in \HHH,\quad i=1,\ldots,N,
\end{align*}
because $\psi''_\eps(\vphi_{\eps,i})$ is bounded for almost any $t\in(0,T)$ (thanks to $\eps>0$ and \eqref{H2} recalling that $\mathbf H^2(\Omega)\hookrightarrow \mathbf L^\infty(\Omega)$). Therefore we can test the equation for $\mathbf w_\eps$ with $\boldsymbol\psi'_\eps(\vphi_\eps)$ to obtain
\begin{align*}
    (\mathbf w_\eps, \boldsymbol\psi'_\eps(\vphi_\eps)) = \dfrac{b}{\varepsilon} (\mathbf P&\boldsymbol\psi'_\eps(\vphi_\eps), \boldsymbol\psi'_\eps(\vphi_\eps)) - \dfrac{b}{\varepsilon}(\mathbf PA\vphi_\eps, \boldsymbol\psi'_\eps(\vphi_\eps)) + \dfrac{2\kappa}{R^2} (\mathbf Pu_\eps\boldsymbol\Lambda, \boldsymbol\psi'_\eps(\vphi_\eps))  \\
        &+ b\varepsilon (\nabla_\Gamma \vphi_\eps, \nabla_\Gamma \boldsymbol\psi'_\eps(\vphi_\eps)) + \kappa( \mathbf P(\boldsymbol\Lambda\cdot\vphi_\eps)\boldsymbol\Lambda, \boldsymbol\psi'_\eps(\vphi_\eps)) + \kappa (\Delta_\Gamma u_\eps, \mathbf P\boldsymbol\Lambda \boldsymbol\psi'_\eps(\vphi_\eps)).
\end{align*}
On account of
$$
(\nabla_\Gamma \vphi_\eps, \nabla_\Gamma \boldsymbol\psi'_\eps(\vphi_\eps)) = \sum_{i=i}^N\int_\Gamma \psi''_\eps(\vphi_{\eps,i})|\nabla_\Gamma \vphi_{\eps,i}|^2\d\Gamma\geq 0
$$
we then deduce
\begin{align}
\begin{split}
\label{5_eq:nonlinear1}
    \dfrac{b}{\varepsilon} (\mathbf P\boldsymbol\psi'_\eps(\vphi_\eps), \boldsymbol\psi'_\eps(\vphi_\eps)) &\leq  (\mathbf w_\eps, \boldsymbol\psi'_\eps(\vphi_\eps)) + \dfrac{b}{\varepsilon}(\mathbf PA\vphi_\eps, \boldsymbol\psi'_\eps(\vphi_\eps)) - \dfrac{2\kappa}{R^2} (\mathbf Pu_\eps\boldsymbol\Lambda, \boldsymbol\psi'_\eps(\vphi_\eps)) \\
    &\hskip 1cm - \kappa( \mathbf P(\boldsymbol\Lambda\cdot\pmb\vphi_\eps)\boldsymbol\Lambda, \boldsymbol\psi'_\eps(\vphi_\eps)) - \kappa (\Delta_\Gamma u_\eps, \mathbf P\boldsymbol\Lambda\cdot \boldsymbol\psi'_\eps(\vphi_\eps))\\
    &\leq  C\left( \|\mathbf w_\eps\|^2 + \|\vphi_\eps\|^2 + \|u_\eps\|^2 + \|\Delta_\Gamma u_\eps\|^2\right) + \dfrac{1}{4N} \int_\Gamma \max_{i=1,\ldots,N} \psi'_\eps(\vphi_{\eps,i})^2.
\end{split}
\end{align}
Observe now that
\begin{align*}
(\mathbf P\boldsymbol\psi'_\eps(\vphi_\eps), \boldsymbol\psi'_\eps(\vphi_\eps)) &= \sum_{k=1}^N  \left( \psi'_\eps(\vphi_{\eps,k}) - \dfrac 1N \sum_{j=1}^N \psi'_\eps(\vphi_{\eps,j})
    \right) \psi'_\eps(\vphi_{\eps,k}) \\
    &= \dfrac{1}{N} \sum_{k,j=1}^N \left(\psi'_\eps(\vphi_{\eps,k}) - \psi'_\eps(\vphi_{\eps,j})\right) \psi'_\eps(\vphi_{\eps,k}) \\
    &= \dfrac{1}{N} \sum_{k<j} \left(\psi'_\eps(\vphi_{\eps,k}) - \psi'_\eps(\vphi_{\eps,j})\right) \psi'_\eps(\vphi_{\eps,k}) + \dfrac{1}{N} \sum_{k>j} \left(\psi'_\eps(\vphi_{\eps,k}) - \psi'_\eps(\vphi_{\eps,j})\right) \psi'_\eps(\vphi_{\eps,k}) \\
    &= \dfrac{1}{N} \sum_{k<j} \left(\psi'_\eps(\vphi_{\eps,k}) - \psi'_\eps(\vphi_{\eps,j})\right)^2\geq 0.
\end{align*}
On the other hand, recalling that $\sum_{i=1}^N\vphi_{\eps,i} = 1$  (see Remark \ref{5_constraint2}), we find
\begin{align*}
    \vphi_{\eps, m} := \min_{i=1,\ldots,N} \vphi_{\eps, i} \leq \dfrac{1}{N} \leq \max_{i=1,\ldots,N} \vphi_{\eps, i} =: \vphi_{\eps, M}
\end{align*}
and using the monotonicity of $\psi'_\eps$ we infer
\begin{align}
\begin{split}
\label{5_eq:nonlinear2}
     \dfrac{1}{N} \sum_{k<j} \left(\psi'_\eps(\vphi_{\eps,k}) - \psi'_\eps(\vphi_{\eps,j})\right)^2 &\geq \dfrac{1}{N} \left(\psi'_\eps(\vphi_{\eps,m}) - \psi'_\eps(\vphi_{\eps,M})\right)^2 \\
     &\geq \dfrac{1}{N}\max_{i} \left( \psi'_\eps(\vphi_{\eps,i}) - \psi'_\eps(1/N)\right)^2 \\
     &\geq \dfrac{1}{N} \max_i \left( \dfrac{1}{2} \psi'_\eps(\vphi_{\eps,i})^2 - \psi'_\eps(1/N)^2\right) \\
     &\geq \dfrac{1}{2N} \max_i \psi'_\eps(\vphi_{\eps, i})^2 - C,
     \end{split}
\end{align}
where we used the basic inequality $(a-b)^2 \geq a^2/2 - b^2$. Notice
that $C$ is independent of $\eps $ provided that we choose $%
\eps $ sufficiently small. Indeed, since we have the pointwise
convergence $\psi' _{\eps }(\tfrac{1}{N})\rightarrow \psi' (\tfrac{1}{N})$
as $\eps \rightarrow 0$, then there exists $C>0$, independent of $\eps$,
such that $|\psi' _{\eps }(\tfrac{1}{N})|\leq C$ for any $%
h \in (0,\eps _{0})$, with $h_0>0$ sufficiently small. Combining \eqref{5_eq:nonlinear1} and \eqref{5_eq:nonlinear2} with the previous uniform estimates, we get
\begin{align}\label{5_psiprime1}
    \dfrac{1}{4N}\int_\Gamma \max_i \psi'_\eps(\vphi_{\eps,i})^2\d\Gamma \leq C \left(1 + \|\mathbf w_\eps\|^2\right).
\end{align}
Integrating in time and using part (i), we infer
\begin{align*}
    \|\boldsymbol\psi'_\eps(\vphi_\eps)\|_{L^2(0,T; \HHH)} \leq C
\end{align*}
and the proof is finished.
\end{proof}

We can now prove some other useful uniform time-weighted estimates, namely,

\begin{lem}
 \Andrea{For any $T>0$ there exists $C=C(T)>0$}, independent of $\eps$, such that:
 \begin{itemize}
     \item[(i)] $\|\sqrt{t}\mathbf w_\eps\|_{L^\infty(0,T; \VVV)}+\|\sqrt{t}\boldsymbol\psi_{\eps}'(\vphi_\eps)\|_{L^\infty(0,T; \HHH)} \leq C$;
     \item[(ii)] $\| \sqrt{t}{{\vphi}_\eps}\|_{L^\infty(0,T; \HHH^2(\Gamma))}+\|\sqrt{t} u_\eps\|_{L^\infty(0,T; H^4(\Gamma))} \leq C$;
     \item[(iii)] $\|\sqrt{t}\mathbf w_\eps\|_{L^\infty(0,T; \LLL^r(\Gamma))}+\|\sqrt{t}\boldsymbol\psi_{\eps}'(\vphi_\eps)\|_{L^\infty(0,T; \LLL^r(\Gamma))}\leq C\sqrt r,\qquad \forall r\geq 2;$
        \item[(iv)] $\| \sqrt{t}{{\vphi}_\eps}\|_{L^\infty(0,T; \WWW^{2,r}(\Gamma))}\leq C(r),\qquad \forall r\geq 2$.

 \end{itemize}
 \label{weigh}
\end{lem}
\Andrea{
\begin{oss}
\label{globalbounds}
    This lemma only gives local regularization estimates, since the constants depend on the time horizon $T>0$. Nevertheless, a similar proof can be performed, using this time the Uniform Gronwall Lemma, to obtain global in time uniform estimates (see also Section \ref{proof3} below).
\end{oss}}
\begin{proof}
(i). Let us recall \eqref{5_lam}: thanks to this estimate, together with \eqref{5_eq:first_estimateb} and the previous lemma, we immediately infer
\begin{align}
\|\sqrt{t}\mathbf w_\eps\|_{L^\infty(0,T; \VVV)}\leq C.
\label{5_muw}
\end{align}
Then, thanks to \eqref{M}, we also obtain
\begin{align}
	\|\sqrt{t}\mathbf w_\eps\|_{L^\infty(0,T; \LLL^r(\Gamma))}\leq C\sqrt r.
	\label{5_muw_r}
\end{align}
Concerning $\boldsymbol\psi'_\eps$, we infer from \eqref{5_psiprime1}, together with \eqref{5_muw}, that it holds
\begin{align}
\|\sqrt{t}\boldsymbol\psi_{\eps}'(\vphi_\eps)\|_{L^\infty(0,T; \HHH)}\leq C.
\label{5_muw1}
\end{align}

(ii). Due to the property $\Vert\textbf{P}\textbf{v}\Vert\leq \Vert\textbf{v}\Vert$, thanks to \eqref{5_eq:first_estimatea}, \eqref{5_muw} and \eqref{5_muw1}, by standard elliptic regularity results we obtain
\begin{align}
\| \sqrt{t}{{\vphi}_\eps}\|_{L^\infty(0,T; \HHH^2(\Gamma))} \leq C.
\label{5_phih2}
\end{align}
Therefore, since we have, in the weak sense,
\begin{align*}
&\sqrt{t}\kappa \Delta_\Gamma^2 u_\eps=\sqrt{t}\left(-\beta \partial_t u_\eps  + \left(\sigma - \dfrac{2\kappa}{R^2} \right) \Delta_\Gamma u_\eps + \dfrac{2\sigma u_\eps}{R^2} - \kappa \vlambda \cdot \Delta_\Gamma \vphi_\eps - \dfrac{2\kappa \vlambda\cdot(\vphi_\eps-\pmb \alpha)}{R^2}\right),
\end{align*}
and the right-hand side belongs to $L^\infty(0,T;H)$, we infer by elliptic regularity, thanks to \eqref{5_eq:first_estimatec}, \eqref{5_eq:first_estimateb} and \eqref{5_phih2}, that
\begin{align}
\| \sqrt{t}{{u}_\eps}\|_{L^\infty(0,T; H^4(\Gamma))} \leq C.
\label{5_phih2bis}
\end{align}

(iii). We perform an argument that is similar to the one used in the proof of Lemma \ref{5_lem:estimates2} part (ii).  Here we set
 $$
\zeta_\eps^r(s)=\psi'_\eps(s)\vert\psi'_\eps(s)\vert^{r-2},
$$
for $r\geq 2$ and  $\boldsymbol\zeta_\eps^r(\vphi_\eps) =(\zeta_\eps^r(\vphi_{\eps,i}))_{i=1,\ldots,N}$. Notice that, combining again the facts that $\psi_\eps^{\prime\prime}(\vphi_{\eps,i})$ and $\psi_\eps'(\vphi_{\eps,i})$ are bounded due to $\vphi_{\eps,i}\in {L}^\infty(\Gamma)$, $i=1,\ldots,N$, which holds thanks to \eqref{5_phih2} and the embedding $H^2(\Gamma)\hookrightarrow L^\infty(\Gamma)$, we thus get, for almost any $t\in(0,T)$,
$$
\tgrad\zeta_\eps^r(\vphi_{\eps,i})=(r-1)\psi_\eps^{\prime\prime}(\vphi_{\eps,i})\vert\psi'_\eps(\vphi_{\eps,i})\vert^{r-2}\nabla_\Gamma\vphi_{\eps,i}\in \HHH
$$
for almost any $t\in(0,T)$. Therefore, we can test the equation for $\mathbf w_\eps$ with $\boldsymbol\zeta^r_\eps(\vphi_\eps)$ to obtain
\begin{align*}
    (\mathbf w_\eps, \boldsymbol\zeta_\eps^r(\vphi_\eps)) = \dfrac{b}{\varepsilon} (\mathbf P&\boldsymbol\zeta_\eps(\vphi_\eps), \boldsymbol\zeta_\eps^r(\vphi_\eps)) - \dfrac{b}{\varepsilon}(\mathbf PA\vphi_\eps, \boldsymbol\zeta^r_\eps(\vphi_\eps)) + \dfrac{2\kappa}{R^2} (\mathbf Pu_\eps\boldsymbol\Lambda, \boldsymbol\zeta^r_\eps(\vphi_\eps)) \\
    &+ b\varepsilon (\nabla_\Gamma \vphi_\eps, \nabla_\Gamma \boldsymbol\zeta^r_\eps(\vphi_\eps))
        + \kappa( \mathbf P(\boldsymbol\Lambda\cdot\vphi_\eps)\boldsymbol\Lambda, \boldsymbol\zeta^r_\eps(\vphi_\eps)) + \kappa (\Delta_\Gamma u_\eps, \mathbf P\boldsymbol\Lambda \cdot\boldsymbol\zeta^r_\eps(\vphi_\eps)).
\end{align*}
Observing that
$$
{(\nabla_\Gamma \vphi_\eps, \nabla_\Gamma \boldsymbol\zeta^r_\eps(\vphi_\eps)) = \sum_{i=i}^N\int_\Gamma (r-1)\vert\psi'_\eps(\vphi_{\eps,i})\vert^{r-2}\psi_\eps^{\prime\prime}(\vphi_{\eps,i})|\nabla_\Gamma \vphi_{\eps,i}|^2\d\Gamma\geq 0},
$$
from the above inequality multiplied by by $t^{\frac{r}{2}}$, we infer
\begin{align}
\begin{split}
\label{5_eq:nonlinear1bis}
    t^{\frac{r}{2}}\dfrac{b}{\varepsilon} (\mathbf P\boldsymbol\zeta_\eps(\vphi_\eps), \boldsymbol\zeta^r_\eps(\vphi_\eps)) &\leq t^{\frac{r}{2}}\left(  (\mathbf w_\eps, \boldsymbol\zeta^r_\eps(\vphi_\eps)) + \dfrac{b}{\varepsilon}(\mathbf PA\vphi_\eps, \boldsymbol\zeta^r_\eps(\vphi_\eps)) - \dfrac{2\kappa}{R^2} (\mathbf Pu_\eps\boldsymbol\Lambda, \boldsymbol\zeta^r_\eps(\vphi_\eps))\right. \\
    &\left.
    \hskip 1cm - \kappa( \mathbf P(\boldsymbol\Lambda\cdot\pmb\vphi_\eps)\boldsymbol\Lambda, \boldsymbol\zeta^r_\eps(\vphi_\eps)) - \kappa (\Delta_\Gamma u_\eps, \mathbf P\boldsymbol\Lambda\cdot \boldsymbol\zeta^r_\eps(\vphi_\eps))\right)\\
    &\leq  C\left( \|\sqrt{t}\mathbf w_\eps\|_{\LLL^r(\Gamma)}^r + \|\sqrt{t}\vphi_\eps\|_{\LLL^r(\Gamma)}^r + \|\sqrt{t}u_\eps\|^r_{L^r(\Gamma)} + \|\sqrt{t}\Delta_\Gamma u_\eps\|_{L^r(\Gamma)}^r\right) \\
    &\hskip 1cm+ \dfrac{t^{\frac r2}}{4N} \int_\Gamma \max_{i=1,\ldots,N} \vert\psi'_\eps(\vphi_{\eps,i})\vert^r\d\Gamma\\&\leq C(\sqrt r)^r+\dfrac{t^{\frac r2}}{4N} \int_\Gamma \max_{i=1,\ldots,N} \vert\psi'_\eps(\vphi_{\eps,i})\vert^r \d\Gamma,
\end{split}
\end{align}
where we applied H\"{o}lder's and Young's inequalities, together with \eqref{5_muw_r}, \eqref{5_phih2} and \eqref{5_phih2bis} (by means of standard Sobolev embeddings).
Note now again that
\begin{align*}
(\mathbf P\boldsymbol\psi'_\eps(\vphi_\eps), \boldsymbol\zeta^r_\eps(\vphi_\eps)) &= \sum_{k=1}^N  \left( \psi'_\eps(\vphi_{\eps,k}) - \dfrac 1N \sum_{j=1}^N \psi'_\eps(\vphi_{\eps,j})
    \right)\zeta_\eps^r(\vphi_{\eps,k}) \\
    &\hskip -5mm = \dfrac{1}{N} \sum_{k,j=1}^N \left(\psi'_\eps(\vphi_{\eps,k}) - \psi'_\eps(\vphi_{\eps,j})\right) \zeta_\eps^r(\vphi_{\eps,k}) \\
    &\hskip -5mm = \dfrac{1}{N} \sum_{k<j} \left(\psi'_\eps(\vphi_{\eps,k}) - \psi'_\eps(\vphi_{\eps,j})\right) \zeta_\eps^r(\vphi_{\eps,k}) + \dfrac{1}{N} \sum_{k>j} \left(\psi'_\eps(\vphi_{\eps,k}) - \psi'_\eps(\vphi_{\eps,j})\right) \zeta_\eps^r(\vphi_{\eps,k}) \\
    &\hskip -5mm = \dfrac{1}{N} \sum_{k<j} \left(\psi'_\eps(\vphi_{\eps,k}) - \psi'_\eps(\vphi_{\eps,j})\right)\left(\zeta_\eps^r(\vphi_{\eps,k}) - \zeta_\eps^r(\vphi_{\eps,j})\right)\geq 0,
\end{align*}
where the last inequality is due the fact that both $\psi'_\eps$ and $\zeta_\eps^r$ are monotone increasing so that
$$
\left(\psi'_\eps(\vphi_{\eps,k}) - \psi'_\eps(\vphi_{\eps,j})\right)\left(\zeta_\eps^r(\vphi_{\eps,k}) - \zeta_\eps^r(\vphi_{\eps,j})\right)\geq 0, \quad \forall k\leq j,\quad k=1,\ldots,N.
$$
Recalling that $\sum_{i=1}^N\vphi_{\eps,i} = 1$  (see again Remark \ref{5_constraint2}), we observe again that we have
\begin{align}
    \vphi_{\eps, m} := \min_{i=1,\ldots,N} \vphi_{\eps, i} \leq \dfrac{1}{N} \leq \max_{i=1,\ldots,N} \vphi_{\eps, i} =: \vphi_{\eps, M}
    \label{5_minmax}
\end{align}
and, using the monotonicity of $\psi'_\eps$ and $\zeta_\eps^r$, we infer
\begin{align*}
&\frac{1}{N}\sum_{k<l}^N\left(\psi'_\eps(\vphi_{\eps,k})-\psi'_\eps(\vphi_{\eps,l})\right)(\zeta_\eps^r(\vphi_{\eps,k})-\zeta_\eps^r(\vphi_{\eps,l}))\\
&\geq \frac{1}{N}\left(\psi'_\eps(\vphi_{\eps,M})-\psi'_\eps(\vphi_{\eps,M})\right)(\zeta_\eps^r(\vphi_{\eps,M})-\zeta_\eps^r(\vphi_{\eps,M}))\\
&= \frac{1}{N}\left\vert\psi'_\eps(\vphi_{\eps,M})-\psi'_\eps(\vphi_{\eps,M})\right\vert\left\vert\zeta_\eps^r(\vphi_{\eps,M})-\zeta_\eps^r(\vphi_{\eps,M})\right\vert\\
&\geq
\frac{1}{N}\max_{k=1,\ldots,N}\left\vert\psi'_\eps(\vphi_{\eps,k})-\psi'_\eps\left(\frac{1}{N}\right)\right\vert\left\vert\zeta_\eps^r(\vphi_{\eps,k})-\zeta_\eps^r\left(\frac{1}{N}\right)\right\vert\\
 &\geq
\frac{1}{N}\max_{k=1,\ldots,N}\left\vert\vert\psi'_\eps(\vphi_{\eps,k})\vert^{r}-\zeta_\eps^r\left(\frac{1}{N}\right)\psi'_\eps(\vphi_{\eps,k})-\psi'_\eps\left(\frac{1}{N}\right)\zeta_\eps^r(\vphi_{\eps,k})+\left\vert\psi'_\eps\left(\frac{1}{N}\right)\right\vert^{r}\right\vert\\
&\geq \frac{1}{N}\max_{k=1,\ldots,N}\left (\vert\psi'_\eps(\vphi_{\eps,k})\vert^{r}-\zeta_\eps^r\left(\frac{1}{N}\right)\psi'_\eps(\vphi_{\eps,k})-\psi'_\eps\left(\frac{1}{N}\right)\zeta_\eps^r(\vphi_{\eps,k})+\left\vert\psi'_\eps\left(\frac{1}{N}\right)\right\vert^{r}\right)\\
&\geq \frac{1}{N}\max_{k=1,\ldots,N}\left(\vert\psi'_\eps(\vphi_{\eps,k})\vert^{r}-\left\vert\psi'_\eps\left(\frac{1}{N}\right)\right\vert^{r-1}\left\vert\psi'_\eps(\vphi_{\eps,k})\right\vert-\left\vert\psi'_\eps\left(\frac{1}{N}\right)\right\vert\left\vert\psi'_\eps(\vphi_{\eps,k})\right\vert^{r-1}+\left\vert\psi'_\eps\left(\frac{1}{N}\right)\right\vert^{r}\right)\\
&\geq \frac{1}{2N}\max_{k=1,\ldots,N}\vert\psi'_\eps(\vphi_{\eps,k})\vert^{r}-C,
\end{align*}
where in the last step we applied Young's inequality to the middle terms.
The constant $C>0$ can be taken independently of $\eps>0$. Indeed, thanks to the pointwise
convergence $\psi' _{\eps }(\tfrac{1}{N})\rightarrow \psi' (\tfrac{1}{N})$
as $\eps \rightarrow 0$, there exists $C>0$, independent of $\eps$,
such that $|\psi' _{\eps }(\tfrac{1}{N})|\leq C$ for any $%
h \in (0,\eps _{0})$, with $h_0>0$ sufficiently small. Combining \eqref{5_eq:nonlinear1bis} and the estimate above leads to
\begin{align}
    \dfrac{t^{\frac r2}}{4N}\int_\Gamma \max_i \vert \psi'_\eps(\vphi_{\eps,i})\vert ^r \d\Gamma\leq C(\sqrt r)^r,
    \label{5_psiprime}
\end{align}
i.e., as desired,
\begin{align}
    \|\sqrt{t}\boldsymbol\psi^\prime_\eps(\vphi_\eps)\|_{L^\infty(0,T; \LLL^r(\Gamma))} \leq C\sqrt r, \quad \forall r\geq 2,
    \label{5_psibis}
\end{align}
for some positive constant $C$ independent of $r$. \medskip

(iv). This bound comes directly from elliptic regularity, by noticing that the right-hand side multiplied by $\sqrt{t}$ belongs to $L^\infty(0,T;\LLL^r(\Gamma))$ (indeed, observe that, by the definition of $\textbf{P}$, $\Vert \textbf{P}\textbf{v}\Vert_{\LLL^r(\Gamma)}\leq C\Vert \textbf{v}\Vert_{\LLL^r(\Gamma)}$ for some positive constant $C$ depending on $N$). Hence, the wanted bound can be deduced from \eqref{5_muw}, \eqref{5_phih2}, \eqref{5_phih2bis} and \eqref{5_psibis}, by applying standard Sobolev embeddings. This concludes the proof.
\end{proof}

\paragraph{Passage to the limit $\eps\to 0$.} Using the previous uniform a priori estimates, we obtain the following convergences (up to a subsequence), as $\eps\to 0$:
\begin{itemize}
    \item[(i)] $u_\eps \overset{\ast}{\rightharpoonup} u$ in $L^\infty(0,T; H^2(\Gamma))$, $\partial_t u_\eps \rightharpoonup \partial_t u$ in $L^2(0,T; H)$, $\sqrt{t}\partial_t u_\eps \overset{\ast}{\rightharpoonup} \sqrt{t}\partial_t u$ in $L^\infty(0,T; H^2(\Gamma))$ and $\sqrt{t}u_\eps \overset{\ast}{\rightharpoonup} \sqrt{t}u$ in $L^\infty(0,T; H^4(\Gamma))$ for any $r\geq2$;
    \item[(ii)] $\vphi_\eps \overset{\ast}{\rightharpoonup} \vphi$ in $L^\infty(0,T; \VVV)$, $\sqrt{t}\vphi_\eps \overset{\ast}{\rightharpoonup} \sqrt{t}\vphi$ in $L^\infty(0,T; \WWW^{2,r}(\Gamma))$ for any $r\geq2$, $\partial_t \vphi_\eps \rightharpoonup \partial_t\vphi$ in $L^2(0,T; \VVV')$ and $\sqrt{t}\partial_t \vphi_\eps \rightharpoonup \sqrt{t}\partial_t\vphi$ in $L^2(0,T; \VVV)$;
    \item[(iii)] $\vphi_\eps \to \vphi$ in $L^2(0,T; \HHH)$ and $\vphi_\eps \to \vphi$ almost everywhere in $\Gamma\times (0,T)$;
    \item[(iv)] $\mathbf w_\eps \rightharpoonup \mathbf w$ in $L^2(0,T; \VVV)$ and $\sqrt t\mathbf w_\eps \overset{\ast}{\rightharpoonup}\sqrt t\mathbf w$ in $L^\infty(0,T; \VVV)$.
\end{itemize}
Most of the terms in the equation are linear so passing to the limit is straightforward. Let us analyse the behavior of the terms involving $\mathbf \Psi_\eps$ and $\boldsymbol\psi'_\eps$.
This is given by the following lemmas.

\begin{lem}\label{5_phiprime}
The following holds true:
\begin{itemize}
    \item[(i)] The limit function $\vphi$ satisfies $\vphi_k \in (0,1)$ for all $k=1,\dots,N$.
    \item[(ii)] As $\eps\to 0$, up to a subsequence, we have
\begin{align}
    \boldsymbol\psi'_\eps(\vphi_\eps) \rightharpoonup  \boldsymbol\psi'(\vphi) \,\, \text{ in } \,\, L^2(0,T; \HHH).
    \label{5_conv}
\end{align}
\end{itemize}
\end{lem}

\begin{proof}
(i). The same argument as in \cite{Gar05} shows that
\begin{align*}
    \lim_{\eps\to 0} \boldsymbol\psi'_\eps(\vphi_\eps) =
    \begin{cases}
     \boldsymbol\psi'(\vphi) \, &\text{ if } \lim_{\eps\to 0}
    \vphi_\eps=\vphi>0,\\
    \infty &\text{ otherwise. }
    \end{cases}
\end{align*}
Because of the estimate (ii) in Lemma \ref{5_lem:estimates2}, by a standard argument, it must be that, in the limit as $\eps\to0$, $\vphi_k>0$ for all $k=1,\dots,N$ almost everywhere in $\Gamma\times(0,T)$. Combining this fact with
\begin{align*}
    \sum_{k=1}^N \vphi_k = 1,
\end{align*}
it follows that $\vphi\in (0,1)^N$ almost everywhere in $\Gamma\times(0,T)$, as desired.
\vskip 3mm

(ii). The previous part shows that
\begin{align}
\boldsymbol\psi'_\eps(\vphi_\eps)\to \boldsymbol\psi'(\vphi)
\label{5_almostev}
\end{align}
almost everywhere in $\Gamma\times (0,T)$. Using Fatou's Lemma and Lemma \ref{5_lem:estimates2}, we have
\begin{align*}
    \int_0^T\int_\Gamma |\boldsymbol\psi'(\vphi)|^2\d\Gamma =
    \int_0^T\int_\Gamma \liminf_{\eps\to 0}|\boldsymbol\psi'_\eps(\vphi_\eps)|^2\d\Gamma \leq \liminf_{\eps\to 0}  \int_0^T\int_\Gamma |\boldsymbol\psi'_\eps(\vphi_\eps)|^2 \d\Gamma\leq C.
\end{align*}
The generalized Lebesgue theorem then implies
\begin{align*}
    \boldsymbol\psi'_\eps(\vphi_\eps) \rightharpoonup \boldsymbol\psi'(\vphi) \,\, \text{ in } \,\, L^2(0,T; \HHH),
\end{align*}
concluding the proof.
\end{proof}

\begin{lem}\label{important}
There exists a constant $C=C(r)>0$ such that
\begin{equation}
    \Vert\sqrt t\boldsymbol\psi^\prime(\vphi)\Vert_{L^\infty(0,T;\LLL^r(\Gamma))}\leq C,\quad \forall r\geq 2.
\end{equation}
\end{lem}

\begin{proof}

Since we know by Lemma \ref{5_lem:estimates2} part (iii) that, up to subsequences, for any $r\in[2,\infty)$ and $k=1,\ldots,N$, there holds
$$
\sqrt t\psi'_\eps(\vphi_{\eps,k})\rightharpoonup\xi\qquad\text{weak* in }L^\infty(0,T;L^r(\Gamma)),
$$
by \eqref{5_conv} we deduce $\xi=\sqrt t\psi'(\vphi_{k})$ and, by weak lower sequential semicontinuity, we infer that
\begin{align}
\Vert\sqrt{t}\psi'(\vphi_{k})\Vert_{L^\infty(0,T;L^r(\Gamma))}\leq C\sqrt t,
\end{align}
for any $k=1,\ldots,N$.
\end{proof}
The convergence results we have obtained so far are sufficient to pass to the limit in Problem \ref{5_pr} and getting the existence of a weak solution satisfying \eqref{wregprop} and \eqref{weakeq} along with the initial conditions.

\paragraph{Energy identity.} Let us now prove that the weak solution satisfies \eqref{enerident}. We argue as in the proof of \cite[Theorem 31]{GGPS}. Let us set $\mmm:=\dashint_\Gamma{\vphi}_{0}$ and introduce the functional
\begin{equation*}
	\mathcal{J}(\vphi):=\frac{b\varepsilon}{2}\int_{\Gamma}\vert \nabla_{\Gamma}\vphi\vert^2\d\Gamma+\frac b\varepsilon\int_\Gamma \boldsymbol \Psi^1(\vphi+\al)\d\Gamma,
\end{equation*}%
whose effective domain in $\HHH_{0}$ is
\begin{equation*}
	\mathfrak{D}(\mathcal{J})=\left\{\mathbf{v}\in \VVV_{0}:\boldsymbol\Psi^1(\mathbf{v}%
	+\mmm)\in \mathbf{L}^{1}(\Gamma )\right\}.
\end{equation*}%
We set $\mathcal J=+\infty$ outside its effective domain.
Clearly $\mathcal{J}$ is proper and convex, but also lower semicontinuous
(w.r.t. the $\HHH$ topology). Moreover, we can immediately
see  (following, for instance, \cite{GGPS}) that, for almost any fixed $t\geq 0$, the solution $(\vphi(t),u(t))$ is
such that $\mathbf{z}(t):=\vphi(t)-\mmm\in
\mathfrak{D}(\mathcal{J})$ and, defining
\begin{align*}
	\mathbf{q}(t)&:=\mathbf w(t)+\mathbf{P}\left( A\vphi(t)-\frac{2\kappa u\vlambda}{R^2}-\kappa(\vlambda\cdot \vphi)\vlambda-\kappa\Delta_\Gamma u\vlambda\right)\\&-\dashint_\Gamma{\mathbf w}%
	(t)-\dashint_\Gamma{\mathbf P\left( A\vphi(t)-\frac{2\kappa u\vlambda}{R^2}-\kappa(\vlambda\cdot \vphi)\vlambda-\kappa\Delta_\Gamma u\vlambda\right)},
\end{align*}%
we easily see that
$$
\mathbf q(t)\in \partial \mathcal{J}(\mathbf{z}%
(t)).
$$

Then, thanks to the regularity of the weak solution and using the Hilbert triplet
$\VVV_{0}\hookrightarrow \hookrightarrow \HHH_{0}\equiv \HHH_{0}^{\prime }\hookrightarrow \VVV
_{0}^{\prime }$, we have:
\begin{itemize}
	\item $\mathcal{J}:\HHH_0\to (-\infty,+\infty]$ is a proper, convex,
	lower semicontinous functional;
	
	\item $\mathbf{z}\in  H^1(0,T;%
	\VVV_0^\prime)\cap L^2(0,T;\VVV_0)$;
	
	\item $\mathbf{q}(t)\in \partial\mathcal{J}(\mathbf{z}(t))$ for almost any $%
	t\in(0,T)$ and $\mathbf{q}\in L^2(0,T;\VVV_0)$;
	
	\item by Poincaré's inequality (whenever $\rrr\in \VVV_0$), and since $\boldsymbol\Psi ^{1}$ is bounded
	below, there exist $k_{1},k_{2}>0$ such that
	\begin{equation*}
		\mathcal{J}(\mathbf{r})\geq k_{1}\Vert \mathbf{r}\Vert ^{2}_{\HHH_0}-k_{2},\quad
		\forall\, \mathbf{r}\in \HHH_{0}.
	\end{equation*}
\end{itemize}

Therefore, we can apply \cite[Lemma 4.1]{RS}, with $H=\HHH_{0}$ and $V=%
\VVV_{0}$ and conclude that the function $\mathcal{J}:t\mapsto
\mathcal{J}(\mathbf{z}(t))\in AC([0,T])$ and
\begin{equation*}
	\int_{s}^{t}\langle \partial _{t}\mathbf{z}(r),\mathbf{q}(r)\rangle _{%
		\VVV_{0}^{\prime },\VVV_{0}}dr=\mathcal{J}(\mathbf{z}(t))-%
	\mathcal{J}(\mathbf{z}(s)).
\end{equation*}%
As a consequence, note that $\vphi \in C([0,T];\VVV)$. Let us now introduce the translated functional
\begin{equation*}
	\tilde{\mathcal{J}}(\mathbf{v}):=\mathcal{J}(\mathbf{v}%
	-\mmm),
\end{equation*}%
and observe that $\mathcal{J}(\mathbf{z}(t))=\tilde{\mathcal{J}}(\vphi(t))$. Thus, we obtain that $\tilde{\mathcal{J}}(\vphi(\cdot))\in AC([0,T])$
and
\begin{align}
\nonumber&	\int_{s}^{t}\left\langle \partial _{t}\vphi(r),\mmu(r)+\mathbf P \left( A\vphi(t)-\frac{2\kappa u\vlambda}{R^2}-\kappa(\vlambda\cdot \vphi)\vlambda-\kappa\Delta_\Gamma u\vlambda\right)\right\rangle _{\VVV^{\prime },\VVV}dr\\&=%
	\tilde{\mathcal{J}}(\vphi(t))-\tilde{\mathcal{J}}(\vphi(s)),
	\label{int}
\end{align}%
recalling the definition of $\mathbf{q}$. This entails
\begin{align*}
	\frac{d}{dt}\tilde{\mathcal{J}}(\vphi(t))&=\left\langle \partial _{t}\vphi(t),\mmu(t)+\mathbf P\left( A\vphi(t)-\frac{2\kappa u\vlambda}{R^2}-\kappa(\vlambda\cdot \vphi)\vlambda-\kappa\Delta_\Gamma u\vlambda\right)\right\rangle _{\VVV^{\prime },\VVV}\\&
	=\left\langle \partial _{t}\vphi(t),\mmu(t)\right\rangle _{\VVV^{\prime },\VVV}\\&+\frac{d}{dt}\int_\Gamma\left(\frac12\vphi(t)\cdot A\vphi(t) -\frac{2\kappa u(t)\vlambda\cdot \vphi(t)}{R^2}-\frac12(\vlambda\cdot\vphi(t))^2-\kappa\Delta_\Gamma u(t) \vlambda\cdot \vphi(t)\right)\\&
	+\int_\Gamma\frac{2\kappa \vlambda\cdot\vphi(t) \partial_t u(t)}{R^2}\d\Gamma+\kappa\int_\Gamma \partial_t u (t)\vlambda\cdot \Delta_\Gamma\vphi(t)\d\Gamma,
\end{align*}
where, since $u\in H^1(0,T;L^2(\Gamma))\cap L^2(0,T;H^4(\Gamma))\hookrightarrow C([0,T];H^2(\Gamma))$, it holds
\begin{align*}
	\frac d{ dt}\int_\Gamma \kappa\Delta_\Gamma u \vlambda\cdot \vphi\d\Gamma=\int_\Gamma \kappa \partial_t u  \vlambda\cdot \Delta_\Gamma\vphi\d\Gamma+\left\langle \partial _{t}\vphi(t),\mathbf P\kappa\Delta_\Gamma u\vlambda\right\rangle _{\VVV^{\prime },\VVV}.
\end{align*}
Now observe that, as an immediate consequence of the regularity of $u$, the energy ${E}_H$ is also absolutely continuous on $[0,T]$, since $\norm{\Delta_\Gamma u}\in AC([0,T])$. It is then easy to show that the energy identity holds,
by testing the equation for $\vphi(t)$ with $\mathbf w(t)$ and the equation for $u(t)$ by $\partial_t u(t)$. Note that in the  equation for $u(t)$ one needs to integrate by parts first, to obtain
\begin{align*}
	\frac12\norm{\Delta_\Gamma u}^2=\int_\Gamma \Delta_\Gamma^2 u \partial_t u \d\Gamma\in L^1(0,T).
\end{align*}

\section{Proof of Theorem \ref{mainb}}\label{proof3}
The weak solution given by Theorem \ref{maina} can be easily extended to $[0,+\infty)$ thanks to the energy identity \eqref{enerident}.
We can argue, for instance, as in \cite[5.4]{GGGP2023}. On the other hand, going back to the proof of Theorem \ref{maina},
it is easy to realize that \eqref{globreg} holds because of the Uniform Gronwall Lemma (see, e.g., \cite[Lemma 1.1]{Temam}). We refer the reader to \cite[Theorem 3.7]{GGPS} and \cite[Theorem 3.10]{AC2023}
for more details on similar results.  It is clear then that \eqref{enerident} holds globally.

We are left to prove the instantaneous strict separation property \eqref{globalssp}.
This is done through a De Giorgi iteration scheme introduced in \cite[Theorem 3.3]{GalPoia}.
Let us fix $\tau>0$. We begin to observe that, thanks to the uniform control $\norm{\mathbf w}_{L^\infty(\tau,\infty;\VVV)}\leq C$, we can deduce, by the Gagliardo-Nirenberg inequality \eqref{M}, that
\begin{align}
	\norm{\mathbf w}_{L^\infty(\tau,\infty;\LLL^r(\Gamma))}\leq C\sqrt r, \quad \forall r\geq 2.
	\label{pp1}
\end{align}
Then, by analogous computations leading to Lemma \ref{important}, we deduce
that we have also
\begin{align}
	\norm{\boldsymbol\psi'(\vphi)}_{L^\infty(\tau,\infty;\LLL^r(\Gamma))}\leq C\sqrt r, \quad \forall r\geq 2,
	\label{pp2}
\end{align}
so that we can directly control $\boldsymbol \mu$, and not only $\mathbf w=\mathbf P \boldsymbol \mu$ in the same space as \eqref{pp1}. This means, {\color{black} combining the uniform bounds associated with \eqref{globreg} (cf. also Remark \ref{globalbounds} )} with \eqref{pp1} and \eqref{pp2}, that
\begin{align}
	&\label{e1b}\norm{u}_{L^\infty(\tau,+\infty;H^4(\Gamma))}\leq C(\tau),\\&
		\norm{\boldsymbol{\mu}}_{L^\infty(\tau,+\infty;\LLL^r(\Gamma))}
+\norm{{\boldsymbol\psi'(\vphi)}}_{L^\infty(\tau,+\infty;\LLL^r(\Gamma))}\leq C(\tau)\sqrt r,\quad \forall r\geq 2.
	\label{e2b}
\end{align}
On the other hand, for any $i=1,\ldots,N$, it holds
\begin{align}
	\boldsymbol\mu_i = \dfrac b\varepsilon \left(\psi'(\vphi_i) - (A\vphi)_i\right) + \dfrac{2\kappa u \vlambda_i}{R^2} - b\varepsilon \Delta_\Gamma \vphi_i + \kappa (\vlambda\cdot \vphi) \vlambda_i + \kappa \Delta_\Gamma u \, \vlambda_i.
	\label{mu_i}
\end{align}
Arguing in the proof of \cite[Theorem 3.3]{GalPoia}, let us fix $\delta \in \left(0,\frac 1N	\right)$ (whose value will be chosen later on), as well as the sequence:
\begin{equation}
	k_{n} = \delta +\frac{\delta }{2^{n}}, \quad \forall n\geq 0,  \label{kn}
\end{equation}
where:
\begin{equation}
	\delta < k_{n+1} < k_{n} < 2\delta, \quad \forall n\geq 1, \quad k_{n}
	\rightarrow \delta \text{ as } n\rightarrow \infty.  \label{kn1}
\end{equation}
Then, for any $i=1,\ldots,N$, we define:
\begin{equation}
	\vphi_{n,i}(x,t) = (\vphi_i - k_{n})^{-}:=-\min\{0,\vphi_i - k_{n}\},  \label{phik0}
\end{equation}
and, for any integer $n\geq 0$, we introduce the set:
\begin{equation*}
	A_{n,i}(t) = \{x\in \Gamma : k_n-\vphi_i(x,t) \geq 0\}, \quad \forall t\in
	[\tau,+\infty).
\end{equation*}
Clearly, we have:
\begin{equation*}
	A_{n+1,i}(t) \subseteq A_{n,i}(t), \quad \forall n\geq 0, \quad \forall t\in
	[\tau,+\infty).
\end{equation*}
Let us now set
\begin{equation*}
	z_{n,i}(t) = \int_{A_{n,i}(t)} 1\d\Gamma, \quad \forall n\geq 0,
\end{equation*}
and fix $t\in [\tau,+\infty)$ (we will not repeat the
dependence on $t$ from now on). For any $n\geq 0$, we consider the test function $v =
-\vphi_{n,i}$, multiply the $i$-th component of equation \eqref{mu_i}, by $v$, and integrate over $\Gamma$.
After an integration by parts, we obtain, for any $i=1,\ldots,N$,
\begin{align}
	&b\varepsilon \norm{\nabla_\Gamma \vphi_{n,i}}^2- \dfrac b\varepsilon \int_\Gamma\psi'(\vphi_i)\vphi_{n,i}\d\Gamma+ \int_\Gamma (A\vphi)_i\vphi_{n,i}\d\Gamma \notag
\\&- \dfrac{2\kappa  \vlambda_i}{R^2}\int_\Gamma u\vphi_{n,i}\d\Gamma - \kappa \int_\Gamma (\vlambda\cdot \vphi) \vlambda_i\vphi_{n,i}\d\Gamma - \kappa \int_\Gamma \Delta_\Gamma u\vphi_{n,i} \, \vlambda_i\d\Gamma
=	\int_\Gamma \boldsymbol\mu_i \vphi_{n,i}\d\Gamma, \label{idssp}
\end{align}
in $[\tau,+\infty)$. Here, we used the identity:
\begin{equation}
	-\int_{A_{n,i}} \nabla_\Gamma \vphi_{i} \cdot \nabla_\Gamma \vphi_{n,i} \d\Gamma= \Vert \nabla_\Gamma
	\vphi_{n,i}\Vert^{2}.  \label{Fss}
\end{equation}
Also, for any $x\in A_{n,i}(t)$, it holds:
\begin{equation}
	\psi^{\prime}(\vphi_i(x,t)) = \psi^{\prime}(k_{n}) + \psi^{\prime \prime}(c(x,t))
	(\vphi_i(x,t) - k_{n}),  \label{c0a}
\end{equation}
with $c(x,t)\in [\vphi_i(x,t),k_n]$. Therefore, considering that $k_{n} <
2\delta$ and, by (\textbf{E0}), $\psi^{\prime \prime}(c(x,t)) \geq \zeta>0$ for
any $x\in A_{n,i}(t)$ and $t\geq \tau$, we can write:
\begin{align}
	& -\frac\beta\varepsilon\int_{\Gamma} \psi^{\prime}(\vphi_i) \vphi_{n,i} \d\Gamma=- \frac b\varepsilon\int_{A_{n,i}(t)}
	\psi^{\prime}(\vphi_{i}) \vphi_{n,i} \d\Gamma \notag \\
	& \geq -\frac b\varepsilon\psi^{\prime}(k_{n}) \int_{\Gamma} \vphi_{n,i}\d\Gamma +\zeta\frac b\varepsilon \int_{\Gamma}
	\vphi_{n,i}^{2}\d\Gamma  \notag \\
	& \geq -\frac b \varepsilon \psi^{\prime}(2\delta) \int_{\Gamma} \vphi_{n,i}\d\Gamma + \zeta\frac b\varepsilon \int_{\Gamma}
	\vphi_{n,i}^{2}\d\Gamma.  \label{c0}
\end{align}
Observe now that
\begin{equation}
	\left\vert\int_{\Gamma} (A\vphi)_i\vphi_{n,i} \d\Gamma\right\vert\leq  C\int_{\Gamma}
	\vphi_{n,i}\d\Gamma,  \label{c1}
\end{equation}
as well as, since $0\leq \vphi\leq1$,
\begin{align}
	\left\vert\kappa \int_\Gamma (\vlambda\cdot \vphi) \vlambda_i\vphi_{n,i}\right\vert\d\Gamma\leq \kappa\vert \vlambda\vert^2 \int_\Gamma \vphi_{n,i}\d\Gamma.
\end{align}
Recalling \eqref{kn1} and  \eqref{phik0}, as first introduced in \cite{P}, we deduce
\begin{align}
0\leq \vphi_{n,i}\leq 2\delta.
	\label{ess1}
\end{align}
Therefore, by \eqref{mu_i} and H\"{o}lder's inequality we get, for almost any $t\geq \tau$,
\begin{align}
	& \int_{\Gamma} \boldsymbol\mu_i \vphi_{n,i}\d\Gamma = \int_{A_{n,i}} \boldsymbol\mu_i \vphi_{n,i} \d\Gamma \notag \\
	& \leq \Vert \vphi_{n,i}\Vert_{L^{\infty}(\Gamma)} \Vert \boldsymbol\mu_i
	\Vert_{L^{p}(\Gamma)} \left(\int_{A_{n,i}} 1\d\Gamma\right)^{1-\frac{1}{p}}  \notag
	\\
	& \leq C(\tau)\delta  z_{n,i}^{1-\frac{1}{p}}, \quad
	\text{for } p\geq 2.  \label{c2}
\end{align}
Moreover, thanks to \eqref{e1}, and recalling the embedding $H^4(\Gamma)\hookrightarrow W^{2,\infty}(\Gamma)$, we obtain
\begin{align}
&\left \vert\dfrac{2\kappa  \vlambda_i}{R^2}\int_\Gamma u\vphi_{n,i}\d\Gamma + \kappa \int_\Gamma \Delta_\Gamma u\vphi_{n,i} \, \vlambda_i\d\Gamma\right\vert\leq \dfrac{2\kappa  \vert\vlambda\vert}{R^2}\norm{u}_{L^\infty(\Gamma)}\int_\Gamma\vphi_{n,i}\d\Gamma+\kappa\vert \vlambda\vert\norm{u}_{W^{2,\infty}(\Gamma)}\int_\Gamma\vphi_{n,i}\d\Gamma\nonumber\\&
\leq  C(\tau)\left(\dfrac{2\kappa  \vert\vlambda\vert}{R^2}+\kappa\vert \vlambda\vert\right)\int_\Gamma\vphi_{n,i}\d\Gamma.
\end{align}
Collecting the above estimates, from \eqref{idssp} we infer
\begin{align}
&\nonumber	\Vert \nabla_\Gamma \vphi_{n,i}\Vert^{2}+\zeta\frac b \varepsilon\Vert\vphi_{n,i}\Vert^2 - \left(\frac b \varepsilon\psi^{\prime}(2\delta)+C(\tau)\left(\dfrac{2\kappa  \vert\vlambda\vert}{R^2}+\kappa\vert \vlambda\vert\right)+\kappa\vert \vlambda\vert^2\right)
	\int_{\Gamma} \vphi_{n,i}\d\Gamma\\&\leq C(\tau)
	\delta \sqrt{p} z_{n,i}^{1-\frac{1}{p}}, \quad \text{for } p\geq 2.
	\label{ast}
\end{align}
By choosing $\delta>0$ sufficiently small, we ensure that (see (\textbf{E1}))
$$
-\frac b	\varepsilon\psi^{\prime}(2\delta)> C(\tau)\left(\dfrac{2\kappa  \vert\vlambda\vert}{R^2}+\kappa\vert \vlambda\vert\right)+\kappa\vert \vlambda\vert^2.
$$
Moreover, for
any $t\in [\tau,\infty)$ and for any $x\in A_{n+1,i}(t)$, observe that
\begin{align}
	& \vphi_{n,i}(x,t) \geq -\vphi_i(x,t) + \left[\delta +\frac{\delta }{2^{n}}\right]
	\notag \\
	& =\underbrace{-\vphi_i(x,t) + \left[\delta +\frac{\delta }{2^{n+1}}\right]}%
	_{\vphi_{n+1,i}(x,t)\geq 0} + \delta \left[\frac{1}{2^{n}}-\frac{1}{2^{n+1}}%
	\right] \geq \frac{\delta }{2^{n+1}},  \label{control}
\end{align}
which implies (see also \cite{GalPoia}):
\begin{equation*}
	\int_{\Gamma}|\vphi_{n,i}|^{3} \d\Gamma\geq \int_{A_{n+1,i}(t)}|\vphi_{n,i}|^{3}\d\Gamma \geq
	\left(\frac{\delta }{2^{n+1}}\right)^{3}\int_{A_{n+1,i}(t)}1\d\Gamma = \left(\frac{%
		\delta }{2^{n+1}}\right)^{3}z_{n+1,i}.
\end{equation*}
Applying H\"{o}lder's inequality, we get:
\begin{equation}
	\left(\frac{\delta }{2^{n+1}}\right)^{3}z_{n+1,i} \leq
	\int_{\Gamma}|\vphi_{n,i}|^{3}\d\Gamma
	= \int_{A_{n,i}(t)}|\vphi_{n,i}|^{3}\d\Gamma \leq
	\left(\int_{\Gamma}|\vphi_{n,i}|^{4}\d\Gamma\right)^{\frac{3}{4}}\left(%
	\int_{A_{n,i}(t)}1\d\Gamma\right)^{\frac{1}{4}}.  \label{est2}
\end{equation}
By means of Gagliardo-Nirenberg's inequality \eqref{M} with $p=4$,
and taking \eqref{ast} into account, we can write:
\begin{align*}
	& \int_{\Gamma}|\vphi_{n,i}|^{4}\d\Gamma \leq C\Vert
	\vphi_{n,i}\Vert_{V}^{2}\Vert \vphi_{n,i}\Vert^{2} \\
	& \leq C\left(\Vert \vphi_{n,i}\Vert^{2} + \Vert \nabla_\Gamma
	\vphi_{n,i}\Vert^{2}\right)\Vert \vphi_{n,i}\Vert^{2} \\
	& \leq C\left(C(\tau)\delta \sqrt{p}z_{n,i}^{1-\frac{1}{p}} +
	C(\tau)\delta \sqrt{p}z_{n,i}^{1-\frac{1}{p}}\right)C(\tau)\delta \sqrt{p}z_{n,i}^{1-%
		\frac{1}{p}} \\
	& \leq C(\tau)\delta^{2}pz_{n,i}^{2-\frac{2}{p}},
\end{align*}
where we used an equivalent norm in $V$. Returning to \eqref{est2}, we immediately obtain:
\begin{equation}
\left(\frac{\delta }{2^{n+1}}\right)^{3}z_{n+1,i} \leq
	\left(\int_{\Gamma}|\vphi_{n,i}|^{4}\d\Gamma\right)^{\frac{3}{4}}z_{n,i}^{\frac{1}{4}}
\leq 4^\frac34C(\tau)^{\frac{3}{4}}\delta^{\frac{3}{2}}z_{n,i}^{\frac{11%
		}{8}},  \label{est3}
\end{equation}
where we have chosen $p=4$. In conclusion, we find
\begin{equation}
	z_{n+1,i} \leq 2^{3n+\frac92}\delta^{-\frac{3}{2}}C(\tau)^{\frac{3}{4}}
	z_{n,i}^{\frac{11}{8}}.  \label{last0tt}
\end{equation}
We can thus apply Lemma \ref{conv}. In particular, we have $b=2^{3}>1$, $%
C=2^{\frac92}\delta ^{-\frac{3}{2}}C(\tau )^{\frac{3}{4}}>0$, and $%
\gamma =\frac38$, which
allows us to conclude that ${z}_{n,i}\rightarrow 0$, as long as, rearranging some constants,
\begin{equation}
	z_{0,i}\leq \dfrac{\delta ^{4}}{2^{\frac{76}{%
			3}}C(\tau )^{2}}. \label{last}
\end{equation}%
Thanks to assumption %
(\textbf{E0}), $\psi^{\prime }$ is monotone in a neighborhood of $%
0$, so that we infer that for any $q\geq 2$ and $\delta >0$ sufficiently small:
\begin{equation*}
	z_{0,i} =\int_{A_{0,i}(t)}1\d\Gamma\leq \int_{\{x\in \Gamma :\ \vphi_i (x,t)\leq
		2\delta \}}1\d\Gamma
	\leq \int_{A_{0}(t)}\frac{|\psi^{\prime }(\vphi_i )|^{q}}{(-\psi^{\prime }(2\delta ))^q}\d\Gamma\leq \frac{\int_{\Gamma }|\psi^{\prime }(\vphi_i )|^{q}\d\Gamma}{(-\psi^{\prime }(2\delta ))^{q}}
	\leq \frac{C_{1}(\tau )^{q}(\sqrt{q})^{q}}{{(-\psi^{\prime }(2\delta ))^{q}}},
\end{equation*}%
where we also used  \eqref{e2}.
If we ensure that
\begin{equation*}
	\frac{C_{1}(\tau )^{q}(\sqrt{q})^{q}}{{(-\psi^{\prime }(2\delta ))^{q}}}\leq
	\dfrac{\delta ^{4}}{2^\frac{76}3C(\tau)^2},
\end{equation*}%
then \eqref{last} holds. This can be obtained as in \cite{GalPoia}. Indeed, let us
recall that, by assumption (\textbf{E2}), there exists $C_{\psi}>0$ such that, for $%
\delta $ sufficiently small:
\begin{equation*}
	-\frac{1}{\psi^{\prime }(2\delta )}\leq \frac{C_{\psi}}{|\ln (\delta )|^{\iota }},\quad \iota>\frac12.
\end{equation*}%
We now make the final step and fix $\delta =e^{-q}$ for $q\geq 2$
sufficiently large, so we get
\begin{equation}
	\frac{C_{1}(\tau )^{q}(\sqrt{q})^{q}}{{(-\psi^{\prime }(2\delta ))^{q}}}\leq
	\frac{C_{1}(\tau )^{q}C_{\psi}^{q}(\sqrt{q})^{q}}{{|\ln (\delta )|^{\iota q}}}=%
	\frac{C_{1}(\tau )^{q}C_{\psi}^{q}(\sqrt{q})^{q}}{{q^{\iota q}}}=\frac{%
		C_{1}(\tau )^{q}C_{\psi}^{q}}{q^{q(\iota -\frac{1}{2})}}.  \label{bass}
\end{equation}%
Condition \eqref{bass}, with this choice of $\delta $, is thus ensured. Indeed, for $q>2$ sufficiently large, it holds
\begin{equation*}
	\frac{C_{1}(\tau )^{q}C_{F}^{q}}{q^{q(\iota -\frac{1}{2})}}\leq C(\tau )%
	\dfrac{e^{-{4q}}}{2^{\frac{76}3}C(\tau)^2},
\end{equation*}%
which is ensured by the condition $\iota>\frac12$ as in assumption (\textbf{E2}).

Therefore, by choosing $q$ as above (i.e. corresponding to a
sufficiently small $\delta =e^{-q}$), we can ensure the validity of
\eqref{bass}. This implies that $z_{n,i}\rightarrow 0$ as $n\rightarrow \infty $.
Passing to the limit in $z_{n,i}$ as $n\rightarrow \infty $, we have shown
that, for any $t\geq \tau$,
\begin{equation*}
	\Vert (\vphi_{i} (t)-\delta)^{-}\Vert _{L^{\infty }(\Gamma)}=0,
\end{equation*}%
since, as $n\rightarrow \infty $, we have
\begin{equation*}
	z_{n,i}(t)\rightarrow \left\vert \left\{ x\in \Omega :\vphi_i (x,t)\leq \delta
	\right\} \right\vert ,
\end{equation*}%
and $z_{n,i}(t)\rightarrow 0$ as $n\rightarrow \infty $, for any fixed $t\geq
\tau $. It is clear that $\delta $ depends on $\tau$, but it does not depend
neither on the specific $t\geq \tau$ nor on $i=1,\ldots,N$, so that the uniform (instantaneous) strict
separation holds.
Note that, since $\vphi_i\in L^\infty(\tau,\infty;W^{2,r}(\Gamma))\cap H^1(t,t+1;H^1(\Gamma))$, for any $r\geq 2$ and for any $t\geq\tau$,
we have that $\vphi_i\in C([\tau,+\infty);H^{3/2}(\Gamma))$, implying $\vphi_i\in C([\tau,+\infty);C({\Gamma}))$. Since $\sum_i^N\vphi_i\equiv 1$, we thus conclude that there exists $\delta (\tau)>0$ such that:
\begin{equation*}
 \delta\leq  \vphi_i (x,t)\leq 1-(N-1)\delta ,\quad \forall\, (x,t) \in \Gamma%
	\times \lbrack \tau ,+\infty).
\end{equation*}%
Since $\delta\in \left(0,\frac1N\right)$ and  does not depend on the component $i=1,\ldots,N$ we choose, we can repeat the same arguments for each $i=1,\ldots,N$, and deduce also that there exists $\delta\in  \left(0,\frac1N\right)$ such that \eqref{globalssp} holds.
Theorem \ref{mainb} is completely proven.

\section{Proof of Theorem \ref{maincon}}\label{convec}

In this section we adapt the results obtained in \cite{GGPS}, detailing the main steps and postponing the proofs of two technical lemmas to the Appendix.
 Given $(\vphi_0,{u}_{0})\in \mathcal{V}_{\al}$, we consider the unique (global) weak solution, given by Theorems \ref{maina} and \ref{mainb}.
In particular, we have that $\vphi\in L^{\infty }(\tau ,\infty ;\mathbf{H}^{2}(\Gamma ))$ for any $\tau >0$, as well as $u\in L^\infty(\tau,\infty; H^4(\Gamma))$. Hence by the relative compactness of the trajectories in $\mathbf{H}^{2r}(\Gamma )\times H^{4r}(\Gamma)$ we infer by standard results that $\omega
(\vphi_0,{u}_{0})$ is non-empty, compact and connected in $\mathbf{H}^{2r}(\Gamma )\times H^{4r}(\Gamma)$. Moreover, it is immediate to show that
\begin{equation}
\lim_{t \rightarrow +\infty }\text{dist}_{\mathbf{H}^{2r}(\Gamma )\times H^{4r}(\Gamma)}((\vphi(t),u(t)),\omega (\vphi_0,{u}_{0}))=0.
\label{conv1}
\end{equation}%
Let introduce the notion of stationary point. Given
\begin{equation*}
\pmb{f}_{1}\in \mathcal{G}:=\{\mathbf{v}\in \mathbf{L}^{\infty }(\Gamma ):\
\mathbf{v}(x)\in T\Sigma \text{ for almost any }x\in \Gamma \},
\end{equation*}%
we say that $(\vphi,u)\in \mathcal{V}_{\al} \cap (\HHH^2(\Gamma)\times H^4(\Gamma))$
is a stationary point if it solves the following system
\begin{align}
\label{steady2b}
\begin{split}
  - b\varepsilon \Delta_\Gamma \vphi+\frac{b}{\varepsilon}\mathbf{P}\boldsymbol\Psi'(\vphi) &= \pmb{f}_1+\mathbf{P}\left(\dfrac b\varepsilon  A\vphi-\dfrac{2\kappa u \vlambda}{R^2}
    - \kappa (\vlambda \cdot \vphi) \vlambda - \kappa \Delta_\Gamma u \, \vlambda\right) \\
    0 &=  -\kappa \Delta_\Gamma^2 u + \left(\sigma - \dfrac{2\kappa}{R^2} \right) \Delta_\Gamma u + \dfrac{2\sigma u}{R^2} - \kappa \vlambda \cdot \Delta_\Gamma \vphi - \dfrac{2\kappa \vlambda\cdot(\vphi-\pmb \alpha)}{R^2}
    \end{split}
\end{align}
almost everywhere on $\Gamma$.

Denote by $\mathcal{W}$ the set of all the stationary points and observe that
$$
 - b\varepsilon \Delta_\Gamma \vphi+\frac{b}{\varepsilon}\mathbf{P}\boldsymbol\Psi'(\vphi) = \pmb f,
$$
almost everywhere on $\Gamma$, where
$$
\pmb f:= \pmb{f}_1+\mathbf{P}\left(\dfrac b\varepsilon  A\vphi-\dfrac{2\kappa u \vlambda}{R^2}  -
\kappa (\vlambda \cdot \vphi) \vlambda - \kappa \Delta_\Gamma u \, \vlambda\right)\in \mathcal{G},
$$
since $u\in H^4(\Gamma)\hookrightarrow W^{2,\infty}(\Gamma)$. Thus, by \cite[Theorem 8.1]{GGPS} (which also holds, by the very same proof, in the case of a surface $\Gamma$),
we infer that $\vphi$ is \textit{strictly separated} from
the pure phases, i.e., there exists $0<{\delta }=\delta (\pmb f)<\frac{1}{N}$
such that
\begin{equation}
\delta <\vphi(x)  \label{prop2}
\end{equation}%
for any $x\in \Gamma$. Thus all the stationary points in $%
\mathcal{W}$ are strictly separated, possibly not uniformly. However,
it can be proven that  $\omega (\vphi_0,{u}_{0})\subset \mathcal{W}$
and that the first component of any element of  $\omega (\vphi_0,{u}_{0})$ is actually \textit{uniformly} strictly
separated from the pure phases. Indeed, we have (see Section \ref{se} below for the proof).

\begin{lem}
\label{separation} For any $(\vphi_0,u_0)\in \mathcal{V}_{\al}$ it holds
$\omega(\vphi_0,u_0)\subset \mathcal{W}$, namely, each element $(\vphi_\infty,u_\infty)\in \omega(\vphi_0,u_0)$ is a solution to the steady-state equation %
\eqref{steady2b}, with
$$
\mathbf{f}_1=\frac{b}{\varepsilon}\dashint_\Gamma{\mathbf{P}\boldsymbol\Psi'(\vphi)}-\dashint_\Gamma{\mathbf{P}\left(\dfrac b\varepsilon  A\vphi-\dfrac{2\kappa u \vlambda}{R^2}
- \kappa (\vlambda \cdot \vphi) \vlambda - \kappa \Delta_\Gamma u \, \vlambda\right)}.
$$
Moreover, it holds $\dashint_\Gamma{\vphi}_\infty=%
\dashint_\Gamma{\vphi}_0$, so that $\omega(\vphi_0,{u}_0)\subset \mathcal{V}_{\al}$, and there exists $\delta>0$ so that
\begin{equation*}
\delta<\vphi_\infty, \quad\forall\, (\vphi_\infty,u_\infty)\,\in \omega(\vphi_0,{u}_0),
\end{equation*}
for any $x\in \Gamma$, i.e., the first component of the any element of the $\omega$-limit set is uniformly strictly separated from the pure phases.
\end{lem}

\begin{oss}
As already noticed, thanks to the constraints on $(\vphi_\infty,u_\infty)\in \omega(\vphi_0,{u}_0)$, the strict separation property also implies that
\begin{equation*}
\vphi_\infty(x)<1-(N-1)\delta
\end{equation*}
for any $x\in \Gamma$.
\end{oss}
Recalling the choice $r\in (\tfrac{1}{2},1)$, we deduce that $\ol$ is
compact in $\mathbf{L}^{\infty }(\Gamma )\times L^\infty(\Gamma)$, so we can take $\xi
>0$ such that, e.g., $\xi <\frac{\delta }{2}$, where $\delta $ is
given by Lemma \ref{separation}, and we can find a finite number $M_{0}$ of $\mathbf{L}^{\infty }\times L^\infty$-open balls $B_{\xi ,n}$ of
radius $\xi $ such that
\begin{equation*}
\omega (\vphi_{0},u_0)\subset \bigcup_{n=1}^{M_{0}}B_{\xi
,n}=:U_{\xi }\subset \mathbf{L}^{\infty }(\Gamma )\times L^\infty(\Gamma),
\end{equation*}%
and $\ol\cap B_{\xi ,n}\not=\emptyset $ for any $%
n=1,\ldots ,M_{0}$. Note that $U_\xi$ is open in $\mathbf{L}^{\infty }(\Gamma )\times L^\infty(\Gamma)$%
. Therefore, thanks to Lemma \ref{separation}, we infer that, for almost any $x\in \Gamma $,
\begin{equation}
0<\delta -\xi \leq \vphi(x)\leq 1-((N-1)\delta -\xi
)<1,\quad \forall\, (\vphi,u)\in U_{\xi }.
\label{deltasep}
\end{equation}%
Furthermore, by \eqref{conv1} and the embeddings $\mathbf{H}^{2r}(\Omega
)\hookrightarrow \mathbf{L}^{\infty }(\Omega )$, $H^{4r}(\Gamma)\hookrightarrow L^\infty(\Gamma)$, we deduce that there exists $%
t^{\ast }\geq 0$ such that $(\vphi(t),{u}(t))\in
U_{\xi }$ for any $t\geq t^{\ast }$.
This means that a strict separation property holds asymptotically, even without assumption (\textbf{E2}), namely,

\begin{thm}
Let the assumptions {\color{black} listed in Subsection \ref{constassum} hold, but (\textbf{E2})}. Then, for any $%
\al\in (0,1)$, $\al\in\Sigma$, and for any initial datum $(\vphi_0,{u}_{0})\in \mathcal{V}_{%
\al}$, there exists $\delta >0$ and $t^{\ast }=t^{\ast }(\vphi_{0},{u}%
_{0})$ such that the corresponding (unique) solution $(\vphi,u)$ satisfies:
\begin{equation}
\delta <\vphi(x,t) < 1-(N-1)\delta,\quad \an{\text{ for any }}(x,t)\in \Gamma \times
[t^{\ast },+\infty ).  \label{sep1}
\end{equation}
\end{thm}

Due to \eqref{sep1}, the singularities of $\psi $ and its derivatives no
longer play a role, as we are only interested in the
behavior of the solution $\vphi(t)$, as $t\rightarrow +\infty $. Therefore we
can modify the function $\psi $ outside the interval $I_{\xi
}=[\delta -\xi ,1-((N-1)\delta -\xi )]$ in such a way that
the extension $\widetilde{\psi }$ is of class $C^{3}(\mathbb{R}^{N})$ and
additionally $|\widetilde{\psi }^{(j)}(s)|$, $j=1,2,3$, are uniformly
bounded on $\mathbb{R}$. We then define $\widetilde{\Psi }(%
\mathbf{s}):=\sum_{i=1}^{N}\widetilde{\psi }(s_{i})-\frac{1}{2}%
A\mathbf{s}\cdot \mathbf{s}$. Observe that $\widetilde{\psi }%
_{\vert I_{\xi }}=\psi $ and $\psi $ is analytic in $I_{\xi }$ by
assumption (\textbf{E3}).
 We now introduce the following operators and spaces: $\mathbb{A}:\VVV_0\to\VVV_0'$ is the (component-wise) Laplace-Beltrami operator restricted to the space of functions $\vvv$ such that $\sum_{i=1}^N\vvv_i\equiv 0$ and $\dashint_\Gamma\vvv\equiv \mathbf0$, whereas $\widetilde{\mathbb{A}}: \mathfrak{D}(\mathbb{A})\subset \HHH_0\to \HHH_0$ is the corresponding unbounded operator, i.e., $-\Delta_\Gamma$ applied component-wise (see also \cite[Appendix 8.1]{GGPS} for a precise description of these operators). We also introduce the spaces
$$
H^2_0:=\{w\in H^2(\Gamma):\ w\perp^{L^2} \{1,\nu_1,\nu_2,\nu_3\}\},
$$
endowed with the $H^2$ topology, and
$$
L^2_0:=\{w\in L^2(\Gamma):\ w\perp^{L^2} \{1,\nu_1,\nu_2,\nu_3\}\},
$$
together with the Hilbert triplet $H^2_0\hookrightarrow L^2_0\hookrightarrow (H_0^2)'$.
We then define $\mathbb{B}: H^2_0\to (H^2_0)'$ as the operator:
\begin{align*}
	\left\langle
	\mathbb{B}v,w\right\rangle_{(H^2_0)',H^2_0}:=\int_\Gamma \left(\kappa\Delta_\Gamma v\Delta_\Gamma w+\left(\sigma-\frac{2\kappa}{R^2}\right)\nabla_\Gamma v\cdot \nabla_\Gamma w-\frac{2\sigma}{R^2}vw\right)\d\Gamma.
\end{align*}
Observe that, by Lax-Milgram's Lemma, this operator is invertible and, by elliptic regularity, $\mathbb{B}^{-1}f\in H^4(\Gamma)\cap H^2_0$ for any $f\in L^2_0$. Therefore, we can define the operator $\widetilde{\mathbb B}: \mathfrak{D}(\widetilde{\mathbb B})\subset L^2_0\to L^2_0$ as the restriction of $\mathbb B$ on $\mathfrak{D}(\widetilde{\mathbb B}):=H^2_0\cap H^4(\Gamma)$. Clearly this operator has compact inverse, like $\widetilde{\mathbb{A}}$.

\noindent We then introduce the \textquotedblleft reduced\textquotedblright energy
$\widetilde{\mathcal{E}} :\mathbf{V}_{0}\times H^2_0\rightarrow \mathbb{R}$
by setting
\begin{align*}
    \widetilde{\mathcal E}(\vphi, u) &= \int_{\Gamma_0} \left(\dfrac{\kappa}{2} (\Delta_{\Gamma} u)^2 + \dfrac 12 \left(\sigma - \dfrac{2\kappa}{R^2}\right)|\nabla_{\Gamma} u|^2 - \dfrac{\sigma u^2}{R^2} +
    \kappa (\vlambda \cdot (\vphi+\al)) \Delta_{\Gamma}u
    \right.\\&\left.
    +\dfrac{2\kappa u \vlambda \cdot (\vphi+\al)}{R^2} + \dfrac{b\varepsilon}{2}|\nabla_{\Gamma} \vphi|^2 + \dfrac{b}{\varepsilon} \widetilde{\Psi}(\vphi+\al) +\dfrac{\kappa (\vlambda\cdot(\vphi+\al))^2}{2}\right)\d\Gamma.
\end{align*}
Observe that $\widetilde{\mathcal{E}}(\vphi_0-\al,u_0)=\mathcal{E}(u_0,\vphi_0)$ for all $%
(\vphi_{0},u_0)\in \mathcal{V}_{\al}\cap U_{\xi}$, thanks to
\eqref{sep1} and to the definition of the extension $\widetilde{\psi }$. We
then recall the following lemma whose proof is based on \cite{Chill1}
(see Appendix \ref{loj} below).

\begin{lem}[{\L}ojasiewicz-Simon inequality]
Let $(\vphi,u)$ be the global solution in the sense of
Theorems \ref{maina} and \ref{mainb}, with $(\vphi_0,{u}_{0})\in \mathcal{V}_{%
\al}$, and suppose $(\vphi_\infty,{u}_{\infty })\in \ol$.
Then, under the additional assumption (\textbf{E3}), there exist $\theta \in (0,\tfrac{1}{2}]$, $C,\sigma >0$ such that
\begin{equation*}
\left\vert \widetilde{\mathcal{E}}(\mathbf{v},\tilde{v})-\widetilde{\mathcal{E}}(%
\vphi_{\infty }-\al,u_\infty)\right\vert ^{1-\theta
}\leq C\Vert \widetilde{\mathcal{E}}^{\prime }(\mathbf{v},\tilde{v})\Vert _{(\mathbf{V}_{0} \times H^2_0)'},
\end{equation*}%
whenever $\Vert (\mathbf{v}-\vphi_{\infty }+\al,\tilde{v}-u_\infty)%
\Vert _{\mathbf{V}_{0}\times H^2_0}\leq \sigma $. \label{Loja}
\end{lem}

As already noticed, for $(\vphi_0,{u}_{0})\in
\mathcal{V}_{\al}$, exploiting the fact that, by \eqref{deltasep}, there
exists $\widehat{U}\subset I_{\xi }$ such that, for any $(\vphi_\infty,{u}
_{\infty })\in \ol\subset U_\xi$,  it holds $\vphi_{\infty }(x)\in \widehat U$
for any $x\in \Gamma$, and recalling the definition of $\widetilde{%
\psi }$, we have that $\widetilde{\mathcal{E}}|_{\ol-(\al,0)}=\mathcal{E}|_{\ol}$. Arguing as
in \eqref{en}, $\mathcal{E}({u}_{\infty },\vphi_\infty)=E_{\infty
}=\lim_{s\rightarrow \infty }\mathcal{E}(u(s),\vphi(s))$ for any $%
(\vphi_{\infty },u_\infty)\in \ol$, so that  $\widetilde{\mathcal{E}}|_{\ol-(\al,0)}$ is equal to $%
E_{\infty }$. By Lemma \ref{Loja}, the Lojasiewicz-Simon inequality is valid
for any $(\vphi_\infty,{u}_{\infty })\in \ol$. This means that there exist
constants $\theta \in (0,\tfrac{1}{2}],\ C>0,\ \sigma >0$ such that
\begin{equation*}
\left\vert \widetilde{\mathcal{E}}(\mathbf{v},\tilde{v})-E_{\infty }\right\vert
^{1-\theta }\leq C\Vert \widetilde{\mathcal{E}}^{\prime }(\mathbf{v},\tilde{v})\Vert _{%
(\mathbf{V}_{0}\times H^2_0)^{\prime }},
\end{equation*}%
for any $(\mathbf{v},\tilde{v})\in \mathbf{V}_{0}\times H^2_0$ such that $\Vert (\mathbf{v}+\al%
-\vphi_{\infty },\tilde{v}-u_\infty)\Vert _{\mathbf{V}_{0}\times H^2_0}\leq
\sigma $. Clearly, this can be restated as
\begin{equation}
\left\vert \widetilde{\mathcal{E}}(\pmb\xi -\al,\tilde{v}%
)-E_{\infty }\right\vert ^{1-\theta }\leq C\Vert \widetilde{\mathcal{E}}%
^{\prime }(\pmb\xi-\al,\tilde{v})\Vert _{\mathbf{V}%
_{0}^{\prime }},
\label{e}
\end{equation}%
for any $(\pmb\xi,\tilde{v}) \in (\mathbf{V}_{0}+\al)\times H^2_0$ such
that $\Vert (\pmb\xi -\vphi_{\infty },\tilde{v}-u_\infty)\Vert _{\mathbf{V}_{0}\times H_0^2}\leq \sigma $%
. Since $\mathbf{H}^{2r}(\Gamma )\hookrightarrow \hookrightarrow \VVV$ and $H^{4r}(\Gamma)\hookrightarrow\hookrightarrow H^2(\Gamma)$, $\ol$ is compact in $\VVV\times H^2(\Gamma)$, thus we can find a finite number $M_{1}$ of $\VVV\times H^2(\Gamma)$-open balls $B_{m}$, $m=1,\ldots ,M_{1}$ of radius $\sigma $, centered at $(\vphi_m,u_m)\in \ol$,
such that
\begin{equation*}
\ol\subset \widetilde{U}:=\bigcup_{m=1}^{M_{1}}B_{m}.
\end{equation*}%
\an{Note that since $(\vphi_m,{u}_m)$ are in finite number, we can easily deduce that \eqref{e} holds \textit{uniformly} for any $(\pmb\xi,\widetilde v) \in (\mathbf{V}_{0}+\al)\times H^2_0\cap \widetilde{U}$.}
From \eqref{conv1}, we
deduce that there exists $\tilde{t}>0$ such that $(\vphi(t),{u}(t))\in \widetilde{U}$ for any $t\geq \widetilde{t}$. Recalling the
definition of $U_{\xi }$ given above in Section \ref{convec}, we have
\begin{equation*}
(\vphi(t),{u}(t))\in U_{\xi },\quad \forall\, t\geq
t^{\ast },
\end{equation*}%
therefore we can choose $\overline{t}:=\max \{\widetilde{t},t^{\ast }\}$ and
$\mathbf{U}=\widetilde{U}\cap U_{\xi}$ such that $(\vphi(t),{u}(t))\in
\mathbf{U}$ for any $t\geq \overline{t}$, entailing:
\begin{equation*}
\Vert (\vphi(t)-\vphi_{\infty }(t),u(t)-u_\infty)\Vert _{\mathbf{V}_{0}\times H^2_0}\leq
\sigma, \quad \forall\, t\geq \overline{t}.
\end{equation*}%
Since $(\vphi(t),u(t))\in (\mathbf{V}_{0}+\al)\times H^2_0\cap \widetilde{U}%
$, it holds
\begin{equation*}
\left\vert {\mathcal{E}}({u}(t),\vphi(t))-E_{\infty }\right\vert ^{1-\theta
}\leq C\Vert \widetilde{\mathcal{E}}^{\prime }(\vphi(t)-\al,u(t))\Vert _{(\mathbf{V}_{0}\times H^2_0)^{\prime }},\quad \forall\, t\geq
\overline{t},
\end{equation*}%
indeed $(\vphi(t),{u}(t))\in U_{\xi }$ for any $t\geq
\overline{t}$, thanks to \eqref{sep1} and the
definition of $\widetilde{\psi }$. Observe now that, for any $t\geq \overline{t}$, it holds
\begin{align*}
&\langle \widetilde{\mathcal{E}}^{\prime }(\vphi(t)-\al,u(t)),(\pmb\xi,v)\rangle _{%
({\VVV}_{0}\times H^2_0)',{\VVV}_{0}\times H^2_0}\\&
=\int_\Gamma \left({b\varepsilon}\nabla_\Gamma \vphi(t)\cdot \nabla_\Gamma \pmb\xi+\frac{b}{\varepsilon}({\boldsymbol\psi}(\vphi(t))-A\vphi(t))\cdot \pmb\xi+\kappa(\vlambda\cdot \pmb\xi)\Delta_\Gamma u(t)+\kappa(\vlambda\cdot \vphi(t))(\vlambda\cdot \pmb\xi)\right)\d\Gamma
\\&+\int_\Gamma \left(\kappa\Delta_\Gamma u(t)\Delta_\Gamma v+\left(\sigma-\frac{2\kappa}{R^2}\right)\nabla u(t)\cdot \nabla v-\frac{2\sigma}{R^2}u(t)v+\kappa(\vlambda \cdot \vphi(t)) \Delta_{\Gamma}v+\frac{2\kappa v\vlambda\cdot \vphi(t)}{R^2}\right)\d\Gamma\\&
=\int_\Gamma \left(-{b\varepsilon}\Delta_\Gamma \vphi(t)+\textbf{P}_0\mathbf{P}\left(\frac{b}{\varepsilon}({\boldsymbol\psi}(\vphi(t))-A\vphi(t)))+\kappa\Delta_\Gamma u(t)\vlambda+\kappa(\vlambda\cdot \vphi(t))\vlambda\right)\right)\cdot \pmb\xi\d\Gamma
\\&+\int_\Gamma \left(\kappa\Delta_\Gamma^2 u(t)-\left(\sigma-\frac{2\kappa}{R^2}\right)\Delta u(t)-\frac{2\sigma}{R^2}u(t) +\kappa(\vlambda \cdot \Delta_\Gamma \vphi(t)) +\frac{2\kappa \vlambda\cdot \textbf{P}_0\vphi(t)}{R^2}\right)v \d\Gamma\\
&
=(\mathbf{w}(t)-\dashint_\Gamma{\mathbf{w}(t)},\pmb\xi)_{\HHH_0}+(\beta\partial_t u(t),v)_{L^2_0} \\
& \leq \Vert \nabla_\Gamma \mathbf{w}(t)\Vert \Vert \pmb\xi\Vert+\beta\Vert \partial_t u(t)\Vert \Vert v\Vert\\& \leq C\left(\sqrt{(\mathbf{L}\nabla_\Gamma \mathbf{w}(t),\nabla_\Gamma \mathbf{w}(t))}\Vert \pmb\xi\Vert _{%
\mathbf{V}_{0}}+\sqrt{\beta}\Vert \partial_t u(t)\Vert \Vert v\Vert_{H^2_0}\right),\quad \forall\, (\pmb\xi,v)\in \mathbf{V}_{0}\times H^2_0,
\end{align*}%
where we exploited
Poincar\'{e}'s inequality and \eqref{pos}. This means that
\begin{equation}
\Vert \widetilde{\mathcal{E}}^{\prime }(\vphi(t)-\al,u(t))\Vert _{(\mathbf{V}_{0}\times H^2_0)^{\prime }}\leq C\left(\sqrt{(\mathbf{L} \nabla_\Gamma
\mathbf{w}(t),\nabla_\Gamma \mathbf{w}(t))}+\sqrt{\beta}\Vert \partial_t u(t)\Vert \right),\quad \forall\, t\geq \overline{t}.
\label{ep}
\end{equation}%
Recalling the energy identity and setting $\HH(t):=\left\vert {%
\mathcal{E}}(u(t),\vphi(t))-E_{\infty }\right\vert ^{\theta }$, by \eqref{ep},
we have that
\begin{align*}
& -\dfrac{d}{dt}\HH(t)=-\theta \dfrac{d{\mathcal{E}}(u(t),\vphi(t))}{dt}%
\left\vert {\mathcal{E}}(u(t),\vphi(t))-E_{\infty }\right\vert ^{\theta -1}
\\
& \geq \theta \dfrac{(\mathbf{L} \nabla_\Gamma \mathbf{w}(t),\nabla_\Gamma \mathbf{w}(t))+\beta\Vert \partial_t u(t)\Vert^2}{%
C\Vert \widetilde{\mathcal{E}}^{\prime }(\vphi(t)-\al,u(t))\Vert _{(\mathbf{V}_{0}\times H^2_0)^{\prime }}} \\
& \geq C\left(\sqrt{(\mathbf{L} \nabla_\Gamma \mathbf{w}(t),\nabla_\Gamma \mathbf{w}(t))}+\sqrt{\beta}\Vert \partial_t u(t)\Vert\right), \quad \forall\, t\geq \overline{t},
\end{align*}%
where we used the inequality $(a+b)\geq \frac{1}{ 2}\left(a^\frac 1 2+b^\frac 1 2\right)^2$, for $a,b\geq0$.
Since $\HH$ is a non nonincreasing nonnegative function such that $%
\HH(t)\rightarrow 0$ as $t \rightarrow +\infty $, we can integrate from $%
\overline{t}$ to $+\infty $ and infer (see \eqref{pos}),
\begin{equation*}
\int_{\overline{t}}^{\infty }\Vert \nabla_\Gamma \mathbf{w}(t)\Vert dt+\int_{\overline{t}}^\infty \Vert \partial_t u(t)\Vert dt\leq C\HH(\overline{t})<+\infty ,
\end{equation*}%
i.e., $\nabla_\Gamma \mathbf{w}\in L^{1}(\overline{t},+\infty ;\HHH)$ and $\partial_t u\in L^1(\overline{t},+\infty ;H)$, entailing by comparison $\partial _{t}\vphi\in L^{1}(%
\overline{t},+\infty ;\VVV^{\prime })$. Hence, there exists $(\vphi_\infty,{{u}}_\infty)\in \ol$  such that
\begin{equation*}
\vphi(t)=\vphi(\overline{t})+\int_{\overline{t}}^{t}\partial _{t}%
\vphi(\tau )d\tau {\longrightarrow }%
{\vphi}_\infty\quad \text{ in }\VVV^{\prime },\quad \text{ as } t \to + \infty,
\end{equation*}%
and \begin{equation*}
{u}(t)={u}(\overline{t})+\int_{\overline{t}}^{t}\partial _{t}%
{u}(\tau )d\tau {\longrightarrow }%
{{u}}_\infty\quad \text{ in }H,\quad \text{ as } t \to + \infty,
\end{equation*}%
and, by uniqueness of the limit, we conclude that  $\omega(\vphi_0)=\{(\vphi_\infty,{{u}}_\infty)\}$.
The proof is finished, since \eqref{sing} can now be obtained from \eqref{conv1}.

\appendix
\section{Appendix}

\subsection{Proof of Theorem \ref{5_thm:energy}}

Since the calculations are lengthy and they have already obtained in \cite{EllHat21,LukeThesis}, we divide the proof into several steps, and we only sketch the main arguments. We wish to calculate the first and second variations of the functional $\mathcal L_\rho$: this requires in particular calculating the variations of $\mathcal W, \mathcal A, \mathcal V$. These are standard and a proof can be checked e.g. in \cite{LukeThesis}. We collect them in the next result:

\begin{lem}\label{5_lem:variations_standard}
The following hold for the variations of $\mathcal W, \mathcal A, \mathcal V$:
\begin{align*}
    \mathcal W'(\Gamma_0) [u\nu_0] &= 0, & \mathcal W''(\Gamma_0)(u\nu_0, u\nu_0) &= \int_{\Gamma_0} (\Delta_\Gamma u)^2 - \dfrac{2}{R^2}|\tgrad  u|^2\d\Gamma, \\
    \mathcal V'(\Gamma_0) [u\nu_0] &= \int_{\Gamma_0} u\d\Gamma, & \mathcal V''(\Gamma_0)(u\nu_0, u\nu_0) &= \int_{\Gamma_0} H_{\Gamma_0} u^2\d\Gamma, \\
    \mathcal A'(\Gamma_0) [u\nu_0] &= \int_{\Gamma_0} H_{\Gamma_0} u\d\Gamma, & \mathcal A''(\Gamma_0)(u\nu_0, u\nu_0) &= \int_{\Gamma_0} |\tgrad  u|^2 + \dfrac{2u^2}{R^2}\d\Gamma.
\end{align*}
\end{lem}

The term involving the first variation of $\mathcal F_2$ is multiplied by $\rho$ and thus vanishes when evaluated at $\rho=0$. It thus suffices to compute the first variations of $\mathcal F_1$ and $\mathcal C$, which we do in the next results.

\begin{lem}
The first variation of $\mathcal F_1$ is given by
\begin{align*}
    \mathcal F_1'(\Gamma_0; \vphi)[u\nu_0] = - \int_{\Gamma_0} -(\mathbf \Lambda \cdot \vphi) \Delta_\Gamma u + \dfrac{2}{R^2} (\mathbf \Lambda\cdot\vphi) u \d\Gamma.
\end{align*}
\end{lem}

\begin{proof}
For $\mathcal F_1$:
\begin{align*}
\dfrac{d}{d\rho} \Big\lvert_{\rho=0} \mathcal F_1(\Gamma_\rho; \phi) &= \dfrac{d}{d\rho} \Big\lvert_{\rho=0} \left[ - \int_{\Gamma_\rho} H_\rho \mathbf\Lambda \cdot \vphi\d\Gamma\right] = -\int_{\Gamma_0} \partial_\rho^\bullet H_\rho \lvert_{\rho=0} \mathbf\Lambda\cdot\vphi\d\Gamma - \int_{\Gamma_0} H_0 \mathbf\Lambda\cdot\vphi \tgrad \cdot(u\nu_0)\d\Gamma,
\end{align*}
where $\Gamma_\rho(u) := \{x + \rho u(x)\nu_0(x) : x \in \Gamma_0\}$, $H_\rho$ is the curvature corresponding to $\Gamma_\rho$, and $\partial^\bullet_\rho$ is the total derivative with respect to $\rho$.
Using $\partial^\bullet_\rho H_\rho = -\Delta_\Gamma u - | H_\rho|^2 u$ (see, e.g., \cite[Lemma A.1.4]{LukeThesis}) and $| H_0|^2 = 2/R^2$ (recall that $\Gamma_0$ is a sphere of radius $R$) we obtain for the first term
\begin{align*}
    \partial^\bullet_\rho H_\rho \mathbf\Lambda\cdot\vphi \Big\lvert_{\rho=0} = -(\mathbf \Lambda \cdot \vphi) \Delta_\Gamma u - (\mathbf\Lambda\cdot\vphi) | H_0|^2 u = -(\mathbf \Lambda \cdot \vphi) \Delta_\Gamma u - \dfrac{2}{R^2} (\mathbf \Lambda\cdot\vphi) u,
\end{align*}
and for the last term, using $\tgrad  \cdot \nu_0 = H_0 = 2/R$ and the fact that $\nabla_\Gamma u\in T\Gamma_0$,
\begin{align*}
    H_\rho \mathbf\Lambda\cdot \vphi_\rho \nabla_{\Gamma_\rho} \cdot (u \nu_\rho) \Big\lvert_{\rho=0} = \dfrac{2}{R} \mathbf\Lambda\cdot \vphi \tgrad  u \cdot \nu_0 + \dfrac{2}{R} (\mathbf\Lambda\cdot \vphi)  u \tgrad  \cdot \nu_0 = \dfrac{4}{R^2}(\mathbf\Lambda\cdot\vphi) u,
\end{align*}
leading to
\begin{align*}
    \mathcal F_1'(\Gamma_0; \vphi)[u\nu_0] &= - \int_{\Gamma_0} -(\mathbf \Lambda \cdot \vphi) \Delta_\Gamma u - \dfrac{2}{R^2} (\mathbf \Lambda\cdot\vphi) u +  \dfrac{4}{R^2}(\mathbf\Lambda\cdot\vphi) u\d\Gamma \\
    &= - \int_{\Gamma_0} -(\mathbf \Lambda \cdot \vphi) \Delta_\Gamma u + \dfrac{2}{R^2} (\mathbf \Lambda\cdot\vphi) u\d\Gamma
\end{align*}
as desired.
\end{proof}

We now prove Theorem \ref{5_thm:energy}.

\begin{proof}
We aim to perform a second order Taylor expansion in $\rho$:
\begin{align*}
    \mathcal L_\rho(\Gamma_\rho, \lambda_\rho; \vphi) &= \mathcal L_0(\Gamma_0, \lambda_0; \vphi) + \rho \dfrac{d}{d\rho} \Big\lvert_{\rho=0} \mathcal L_\rho(\Gamma_\rho, \lambda_\rho; \vphi) + \dfrac{\rho^2}{2} \dfrac{d^2}{d\rho^2}\Big\lvert_{\rho=0} \mathcal L_\rho(\Gamma_\rho, \lambda_\rho; \vphi) + \mathcal O(\rho^3).
\end{align*}
As for the first term, we have
\begin{align*}
    \mathcal L_0(\Gamma_0, \lambda_0; \vphi) = \kappa \mathcal W(\Gamma_0) + \sigma \mathcal A(\Gamma_0) = \left( \dfrac{2\kappa}{R^2} + \sigma\right) |\Gamma_0|.
\end{align*}
For the second term, we have, recalling that $(\Gamma_0, \lambda_0)$ is a critical point of $\mathcal L$,
\begin{align*}
    \dfrac{d}{d\rho} \Big\lvert_{\rho=0} \mathcal L_\rho(\Gamma_\rho,
    \lambda_\rho; \vphi) &= \mathcal F(\Gamma_0,
    \lambda_0; \vphi) = \kappa \mathcal F_1(\Gamma_0; \vphi) = -\dfrac{2\kappa}{R} \mathbf \Lambda\cdot \al.
\end{align*}
Finally, for the second order term, we have at $\rho=0$
\begin{align*}
    \dfrac{d^2}{d\rho^2} \Big\lvert_{\rho=0} \mathcal L_\rho(\Gamma_\rho,
    \lambda_\rho; \vphi) &= \kappa \mathcal W''(\Gamma_0)[u\nu_0, u\nu_0] + \sigma \mathcal A''(\Gamma_0)[u\nu_0, u\nu_0] + \lambda_1 \mathcal V'(\Gamma_0)[u\nu_0] \\
    &\hskip 5mm +\lambda_0 \mathcal V''(\Gamma_0)[u\nu_0, u\nu_0] + 2 \kappa \mathcal F_1'(\Gamma_0; \vphi)[u\nu_0] + 2\mathcal F_2(\Gamma_0; \vphi)[u\nu_0].
\end{align*}
Using the formulas in Lemma \ref{5_lem:variations_standard} and recalling $\lambda_0 = -2\sigma/R$, this becomes (see also \cite[Theorem 2.3.1]{LukeThesis} for similar computations)
\begin{align*}
    \dfrac{d^2}{d\rho^2} \Big\lvert_{\rho=0} \mathcal L_\rho(\Gamma_\rho, \vphi_\rho,
    \lambda_\rho, \mu_\rho)
    &= 2 \mathcal E(\vphi, u).
\end{align*}
All in all, the calculations lead to
\begin{align*}
    \mathcal L_\rho(\Gamma_\rho, \lambda_\rho; \vphi) &= \left( \dfrac{2\kappa}{R^2} + \sigma\right) |\Gamma_0| -\dfrac{2\kappa}{R} \mathbf \Lambda\cdot \al  \, \rho + \rho^2 \mathcal E(\vphi, u) + \mathcal O(\rho^3),
\end{align*}
so that \eqref{5_eq:taylor_expansion} holds with
\begin{align*}
    C_1 = \left( \dfrac{2 \kappa}{R^2} + \sigma \right) |\Gamma_0|, \quad C_2 = -\dfrac{2\kappa}{R} \mathbf \Lambda \cdot \al.
\end{align*}
\end{proof}

\subsection{Proof of Lemma \protect\ref{separation}}

\label{se}
Let us consider $(\vphi_\infty,{u}_{\infty })\in \omega (\vphi_0,u_0)$. By definition of $\omega $-limit set there exists a sequence $t_{n}\rightarrow +\infty $ such that $(\vphi(t_n),{u}(t_{n}))\rightarrow (\vphi_\infty,{u}_{\infty })$ in $\mathbf{H}^{2r}(\Gamma)\times H^{4r}(\Gamma)$ as $%
n\rightarrow \infty $. We then define the sequences $\vphi_{n}(t):=\vphi(t+t_{n})$, $\mathbf{w}_{n}(t):=\mathbf{w}(t+t_{n})$, and $u_n(t):=u(t+t_n)$.
By Theorem \ref{mainb}, we immediately infer that $\vphi_{n}$ is
uniformly (in $n$) bounded in $L^{\infty }(0,\infty ;\mathbf{H}^{2}(\Gamma
))\cap H^{1}(0,T ;\VVV)$, for any $T>0$, $\mathbf{w}_{n}$ is uniformly bounded in $L^{\infty }(0,\infty ;%
\VVV)$, $%
\boldsymbol\psi' (\vphi_n)$ is uniformly bounded in $L^{\infty }(0,\infty ;%
\HHH)$, and $u_n$ is uniformly bounded in the spaces $L^\infty(0,\infty;H^4(\Gamma))$ and $H^1(0,T;H^2(\Gamma))$, for any $T>0$. From these bounds, by passing to the limit, up to
subsequences, in the equations solved by $({u}_{n},\vphi_n)$, we infer the
existence of $(\vphi^\ast,{u}^{\ast })$ which is a solution to Problem \ref{5_pr}. In particular, concerning the initial datum, $\vphi^{\ast
}(0)=\lim_{n\rightarrow \infty }\vphi_{n}(0)=\lim_{n\rightarrow \infty }%
\vphi(t_{n})={\vphi_{\infty }}$, where the limit is intended in
the sense, e.g., of $\VVV$. Analogously, we have $u^{\ast
}(0)=\lim_{n\rightarrow \infty }u_{n}(0)=\lim_{n\rightarrow \infty }%
u(t_{n})={u_{\infty }}$, in
the sense, e.g., of ${H}^{2}(\Gamma )$. We thus have $\lim_{n\rightarrow
\infty }\mathcal{E}(u_n(t),\vphi_n(t))=\mathcal{E}({u}^{\ast }(t),\vphi^\ast(t))$
for all $t\geq 0$. Thanks to the energy inequality, we see that the energy $%
\mathcal{E}({u}(t),\vphi(t))$ is nonincreasing in time, therefore there exists $%
E_{\infty }$ such that $\lim_{t\to +\infty }{\mathcal{E}}({u}(t),\vphi(t))=E_{\infty }$. Thus we have
\begin{equation}
{\mathcal{E}}({u}^{\ast }(t),\vphi^\ast(t))=\lim_{n\rightarrow \infty }{\mathcal{E}%
}({u}_{n}(t),\vphi_n(t))=E_{\infty },\quad \forall\, t\geq 0,  \label{en}
\end{equation}%
entailing that ${\mathcal{E}}({u}^{\ast }(t),\vphi^\ast(t))$ is constant in time.
Passing then to the limit as $n\rightarrow \infty $ in the energy inequality, we infer
\begin{equation*}
E_{\infty }+\int_s^t \beta\Vert \partial_t u^\ast(\tau)\Vert^2d\tau+\int_{s}^{t}(\mathbf{L} \nabla_\Gamma \mathbf{w}^{\ast }(\tau ),\nabla_\Gamma
\mathbf{w}^{\ast }(\tau ))d\tau \leq E_{\infty },\quad \forall\, t\geq s,\quad\text{for a.a. }s\geq 0,
\end{equation*}%
where $\mathbf{w}^{\ast }$ is the chemical potential corresponding to $%
\vphi^{\ast }$, implying that $\nabla_\Gamma \mathbf{w}^{\ast }=\textbf{0}$ and $\partial_t u^\ast=0$ almost
everywhere in $\Gamma \times (0,\infty) $, and thus, by comparison, $%
\partial _{t}\vphi^{\ast }=\textbf{0}$ almost everywhere in $\Gamma \times (0,\infty)$.
As a consequence, we infer that $\vphi^{\ast }$ and $u^\ast$
are constant in time, namely $\vphi^{\ast }(t)\equiv \vphi^{\ast
}(0)=\vphi_{\infty }$ and ${u}^{\ast }(t)\equiv {u}^{\ast
}(0)={u}_{\infty }$ for all $t\geq 0$, and $\mathbf{w}^{\ast }$ is
constant in space and time. This means that $(\vphi_\infty,{u}_{\infty })$
satisfies \eqref{steady2b1}.
It is then easy to see that $\dashint_\Gamma{\vphi%
^{\ast }}(t)=\lim_{n\rightarrow \infty }\dashint_\Gamma{\vphi}_{n}(t)\equiv
\dashint_\Gamma{\vphi_{0}}$ for any $t\geq 0$, thus $\dashint_\Gamma{\vphi}%
_{\infty }=\dashint_\Gamma \vphi_0$. Analogously, it is immediate to verify that $\sum_{i=1}^N\vphi_{\infty,i}\equiv 1$ and $u_\infty\in \text{span}\{1,\nu_1,\nu_2,\nu_3\}^{\perp}$. Moreover, for further use we notice that
\begin{equation*}
\mathcal{E}({u}_{\infty },\vphi_\infty,)=E_{\infty }=\lim_{s\rightarrow +\infty }%
\mathcal{E}({u}(s),\vphi(s))=\inf_{s\geq 0}\mathcal{E}({u}(s),\vphi(s))\leq
\mathcal{E}({u}(t),\vphi(t)),\quad \forall\, t\geq 0.
\end{equation*}%
We can thus infer that $(\vphi_\infty,{u}_{\infty })\in \mathcal{W}\cap \mathcal{V}_{\al}$. Therefore,  $\ol\subset \mathcal{W}\cap \mathcal{V}_{\al}$.
In conclusion, the fact that all the elements $\vphi$ of the pairs $(\vphi,u)$ in $\ol$ are \textit{uniformly} strictly separated can be proved exactly as in \cite[Lemma 3.11]{GGPS}, to which we refer. The proof is finished.

\subsection{Proof of Lemma \protect\ref{Loja}}
\label{loj}
The proof is similar to \cite [Lemma 3.15]{GGPS}. The first Fr\'{e}chet derivative of $\widetilde{\mathcal{%
E}}$ reads as follows (recall that $\widetilde{\Psi }$ is smooth):
\begin{align*}
&\langle \widetilde{\mathcal{E}}^{\prime }(\vphi,u),(\pmb\xi,v)\rangle _{%
({\VVV}_{0}\times H^2_0)',{\VVV}_{0}\times H^2_0}\\&=
\int_\Gamma \left({b\varepsilon}\nabla_\Gamma \vphi\cdot \nabla_\Gamma \pmb\xi+\frac{b}{\varepsilon}\widetilde{\Psi}_{,\vphi}(\vphi+\al)\cdot \pmb\xi+\kappa(\vlambda\cdot \pmb\xi)\Delta_\Gamma u+\kappa(\vlambda\cdot (\vphi+\al))(\vlambda\cdot \pmb\xi)\right)\d\Gamma
\\&+\int_\Gamma \left(\kappa\Delta_\Gamma u\Delta_\Gamma v+\left(\sigma-\frac{2\kappa}{R^2}\right)\nabla u\cdot \nabla v-\frac{2\sigma}{R^2}uv+\kappa(\vlambda \cdot (\vphi+\al)) \Delta_{\Gamma}v+\frac{2\kappa v\vlambda\cdot (\vphi+\al)}{R^2}\right)\d\Gamma,
\end{align*}%
where $\widetilde{\Psi }_{,\vphi}(\uuu):=(\tilde{\psi}^{\prime }(%
\uuu_i)-(A\mathbf{\uuu})_{i})_{i=1,\ldots ,N}.$ Notice that $(\vphi%
_{\infty }-\al,u_\infty)\in \omega (\vphi_0,u_0)-%
(\al,0)$ is a critical point for $\widetilde{%
\mathcal{E}}$. Indeed, for $(\vphi_0,u_{0})\in \KK$, with $0\leq \vphi_0\leq 1$ and $\sum_{i=1}^N\vphi_{0,i}=1$,
thanks to the fact that, for any $(\vphi_\infty,{u}_{\infty })\in \omega (\vphi_0,{u}%
_{0})$, by Lemma \ref{separation}, there exists a set $\tilde{U}\subset
I_{\xi }$ ($I_{\xi }$ is defined in Section \ref{convec}) such
that $\vphi_{\infty }(x)\in \tilde{U}$ for any $x\in \Gamma$
and due to the definition of $\widetilde{\psi }$, we have $\widetilde{%
\mathcal{E}}|_{\omega (\vphi_0,{u}%
_{0})-(\al,0)}=%
\mathcal{E}|_{\omega (\vphi_0,{u}%
_{0})}$. Therefore, we have
\begin{align*}
&\langle \widetilde{\mathcal{E}}^{\prime }(\vphi_{\infty }-\al,u_\infty),(\pmb\xi,v)\rangle _{%
({\VVV}_{0}\times H^2_0)',{\VVV}_{0}\times H^2_0}\\&
=\int_\Gamma \left({b\varepsilon}\nabla_\Gamma \vphi_\infty\cdot \nabla_\Gamma \pmb\xi+\frac{b}{\varepsilon}(\boldsymbol\psi'(\vphi_\infty)-A\vphi_\infty)\cdot \pmb\xi+\kappa(\vlambda\cdot \pmb\xi)\Delta_\Gamma u_\infty+\kappa(\vlambda\cdot \vphi_\infty)(\vlambda\cdot \pmb\xi)\right)\d\Gamma
\\&+\int_\Gamma \left(\kappa\Delta_\Gamma u_\infty\Delta_\Gamma v+\left(\sigma-\frac{2\kappa}{R^2}\right)\nabla u_\infty\cdot \nabla v-\frac{2\sigma}{R^2}u_\infty v+\kappa(\vlambda \cdot \vphi_\infty) \Delta_{\Gamma}v+\frac{2\kappa v\vlambda\cdot \vphi_\infty}{R^2}\right)\d\Gamma\\&
=\int_\Gamma \left(-{b\varepsilon}\Delta_\Gamma \vphi_\infty+\textbf{P}_0\mathbf{P}\left(\frac{b}{\varepsilon}(\boldsymbol\psi'(\vphi_\infty)-A\vphi_\infty)+\kappa\Delta_\Gamma u_\infty\vlambda+\kappa(\vlambda\cdot \vphi_\infty)\vlambda\right)\right)\cdot \pmb\xi \d\Gamma
\\&+\int_\Gamma \left(\kappa\Delta_\Gamma^2 u_\infty-\left(\sigma-\frac{2\kappa}{R^2}\right)\Delta u_\infty-\frac{2\sigma}{R^2}u_\infty +\kappa(\vlambda \cdot \Delta_\Gamma \vphi_\infty) +\frac{2\kappa \vlambda\cdot \textbf{P}_0\vphi_\infty}{R^2}\right)v \d\Gamma\\&
=0,\quad
\forall\, (\pmb \xi,v)\in \VVV_0\times H^2_0,
\end{align*}%
where $\textbf{P}_0$ is the $\textbf{L}^{2}$- projector onto the subspace with zero spatial average.
Recall that $(\vphi_{\infty },u_\infty)$ satisfies \eqref{steady2b} (see Lemma \ref{separation}).
Concerning the second Fr\'{e}chet derivative, it is easy to show that
\begin{align*}
&\langle \widetilde{\mathcal{E}}^{\prime\prime }(\vphi,u)(\pmb\xi_1,v_1),(\pmb\xi_2,v_2)\rangle _{%
({\VVV}_{0}\times H^2_0)',{\VVV}_{0}\times H^2_0}\\&=
\int_\Gamma \left({b\varepsilon}\nabla_\Gamma \pmb\xi_1\cdot \nabla_\Gamma \pmb\xi_2+\frac{b}{\varepsilon}\sum_{i=1}^N\widetilde{\psi}''(\vphi_i+\al_i)\pmb\xi_{1,i}\pmb\xi_{2,i}+\kappa(\vlambda\cdot\pmb\xi_1)(\vlambda\cdot \pmb\xi_2)\right)\d\Gamma
\\&+\int_\Gamma \left(\kappa\Delta_\Gamma v_1\Delta_\Gamma v_2+\left(\sigma-\frac{2\kappa}{R^2}\right)\nabla v_1\cdot \nabla v_2-\frac{2\sigma}{R^2}v_1v_2\right)\d\Gamma,
\end{align*}%
for all  $(\vphi,u), (\pmb\xi_{1},v_1),(\pmb\xi_{2},v_2)\in \VVV_{0}\times H^2_0 $.
Let us set $\mathcal{L}:=\widetilde{\mathcal{E}}
^{\prime \prime }\in \mathcal{B}(\VVV_{0}\times H_0^2,\VVV_{0}'\times (H_0^2)')$, omitting the dependence on $(\vphi,u)$, which will be pointed out if
necessary, and recalling the standard result that the space $(\VVV_{0}\times H_0^2)'$ is isomorphic to $\VVV_{0}'\times (H_0^2)'$. We also consider the Hilbert triplet $$\VVV_{0}\times H_0^2\hookrightarrow\hookrightarrow \HHH_0\times L^2_0\cong (\HHH_0\times L^2_0)'\cong \HHH_0'\times (L^2_0)'\hookrightarrow (\VVV_{0}\times H_0^2)'\cong \VVV_{0}'\times (H_0^2)'.$$ We then observe that, for all $(\mathbf{z},w)\in Ker(%
\mathcal{L})\subset \VVV_{0}\times H_0^2$,
setting
$$
\mathbf{v}:=(\widetilde{\psi }^{\prime \prime }(%
\vphi_i+\al_{i})\mathbf{z}_{i}-(A\mathbf z%
)_{i})_{i=1,\ldots ,N}\in \mathbf{L}^{2}(\Gamma ),
$$
we have
\begin{equation*}
\langle \mathbb{A}\mathbf{z},\mathbf{h}\rangle _{\VVV_0^{\prime },%
\VVV_0}+\langle\mathbb{B}w,q\rangle_{(H^2_0)',H_0^2}=-(\mathbf{v},\mathbf{h})=-(\textbf{P}_{0}\mathbf{P}\mathbf{v},\mathbf{%
h}), \quad\forall\,(\mathbf{h},q)\in \VVV_{0}\times H_0^2.
\end{equation*}%
This means, recalling the identification $\VVV_0\hookrightarrow \hookrightarrow
\HHH_0\equiv \HHH_0^{\prime }\hookrightarrow \VVV%
_{0}^{\prime }$ and reasoning component-wise,
\begin{equation*}
\mathbb{A}\mathbf{z}=-\textbf{P}_{0}\mathbf{P}\mathbf{v}\in \HHH_0,\quad w\equiv 0,
\end{equation*}%
implying that $\mathbf{z}\in \mathfrak{D}(\widetilde{\mathbb{A}})$. This entails $Ker(\mathcal{L})\subset \mathfrak{D}(%
\widetilde{\mathbb{A}})\times\{0\}$. Since the second component of the elements of $Ker(\mathcal{L})$ is zero, we can basically follow the same argument of \cite[Lemma 3.15]{GGPS} treating the first component. We thus only sketch the main steps.  Let us introduce the operator $\mathbb{Q}\in \mathcal{B}(%
\VVV_0^{\prime })$ such that, for any $\mathbf{z}\in \VVV%
_{0}^{\prime }$,
\begin{equation*}
\langle \mathbb{Q}\mathbf{z},\mathbf{w}\rangle _{\VVV_0^{\prime },%
\VVV_0}=\left\langle \mathbf{z},\textbf{P}_{0}\mathbf{P}\left( \widetilde{%
\psi }^{\prime \prime }(\vphi_i+\al_{i})\mathbf{w}%
_{i}-(\mathbf{Aw})_{i}\right) _{i=1,\ldots ,N}\right\rangle _{\VVV%
_{0}^{\prime },\VVV_0},\quad \forall\, \mathbf{w}\in \VVV_0,
\end{equation*}%
well defined since $\widetilde{\psi }$ is of class $C^{3}(\mathbb{R})$
and $\mathbf{u}\in \VVV_0$. Therefore, for any $(\mathbf{z},v)\in
\mathfrak{D}(\mathcal{L})=\mathfrak{D}(\mathbb{A})\times \mathfrak{D}(\mathbb{B})=\VVV_0\times H^2_0$, we have
\begin{equation*}
\mathcal{L}(\mathbf{z},v)=
(\mathbb{A}\mathbf{z}+\mathbb{Q}\mathbf{z},\mathbb{B}v)=\begin{bmatrix}
    \mathbb{A}+\mathbb{Q}& 0
    \\ \mathbf{0}&\mathbb{B}
\end{bmatrix}\begin{bmatrix}
    \mathbf{z}\\ v
\end{bmatrix},
\end{equation*}%
since, for $\mathbf{z}\in \HHH_0$, being $A$ symmetric,
\begin{equation*}
\mathbb{Q}\mathbf{z}=\textbf{P}_{0}\mathbf{P}\left( \widetilde{\psi }^{\prime \prime
}(\vphi_i+\al_{i})\mathbf{z}_{i}-(A\mathbf z%
)_{i}\right) _{i=1,\ldots ,N}\in \HHH_0,
\end{equation*}%
as in the definition of $\mathcal{L}$. Note that, thanks to the regularity of $\vphi\in
\VVV_0$ and $\widetilde{\psi }^{\prime \prime }$, if $\mathbf{z}\in
\VVV_0$ we also have $\mathbb{Q}\mathbf{z}\in \VVV_0$.
Now, ${%
\mathbb{A}}^{-1}:\VVV_0^{\prime }\rightarrow \mathfrak{D}(\mathbb{A}%
)=\VVV_0\hookrightarrow \VVV_0^{\prime }$, $\widetilde{%
\mathbb{A}}^{-1}:\HHH_0\rightarrow \mathfrak{D}(\widetilde{\mathbb{A}%
})$, ${%
\mathbb{B}}^{-1}:(H^2_0)^{\prime }\rightarrow \mathfrak{D}(\mathbb{B}%
)=H^2_0\hookrightarrow (H^2_0)^{\prime }$, $\widetilde{%
\mathbb{B}}^{-1}:L_0^2\rightarrow \mathfrak{D}(\widetilde{\mathbb{B}%
})$  are compact in $\VVV_0^{\prime }$, $\HHH_0$, $(H^2_0)'$ and $L^2_0$,
respectively.
In particular, we have that both $\mathbb{A}^{-1}\mathbb{Q}$ and $\mathbb{Q}\mathbb{A%
}^{-1}$ are compact operators on $\VVV_0^{\prime }$.
Thus we can apply \cite[Theorem 8.2]{GGPS} to $A=\mathcal{L}:\mathfrak{D}(\mathbb{A%
})\times \mathfrak{D}(\mathbb{B%
})\hookrightarrow \VVV_0^{\prime }\times (H^2_0)'\rightarrow \VVV_0^{\prime
}\times (H^2_0)'$, with the bounded operator (in the notation of the theorem) $${T}=\begin{bmatrix}
      \mathbb{A}^{-1}& 0\\
      \mathbf{0}&\mathbb{B}^{-1}
\end{bmatrix}
 \in \mathcal{B}(\VVV_0^{\prime }\times (H_0^2)'),$$ to
deduce that $\mathcal{L}$ is a Fredholm operator, implying that $Ker(%
\mathcal{L})\subset \VVV_0^{\prime }\times\{0\}$ is finite dimensional and $%
Range(\mathcal{L})$ is closed in $\VVV_0^{\prime }\times (H^2_0)'$ (and thus in $(\VVV_0\times H^2_0)'$). Moreover, since it is easy to see (cf. \cite{GGPS}) that the operators $\mathbb{A}$ and $\mathbb{B}$ are selfadjoint with respect to
the Hilbert adjoint, we can write the Hilbert
adjoint of $\mathcal{L}$ as $\mathcal{L}^{\ast }=\begin{bmatrix}
\mathbb{A}+\mathbb{Q}^{\ast}&0\\
\mathbf{0}&\mathbb{B}
    \end{bmatrix}
$. Now, by the Closed Range Theorem, recalling that $\mathcal{L}^{\prime }=%
\begin{bmatrix}
\mathbb{A}^{-1}&0\\
\mathbf{0}&\mathbb{B}^{-1}
    \end{bmatrix}\mathcal{L}^{\ast }\begin{bmatrix}
\mathbb{A}&0\\
\mathbf{0}&\mathbb{B}
    \end{bmatrix}$ ($\mathcal{L}^{\prime }$ being
the adjoint of $\mathcal{L}$), observe that
\begin{align}
&Range(\mathcal{L}) \notag\\
&=\{\mathbf{y}^{\ast }
\in \VVV_0^{\prime }\times (H^2_0)':\
\langle \mathbf{y}^{\ast }_1,\mathbf{x}_1\rangle _{\VVV_0^{\prime },%
\VVV_0}+\langle \mathbf{y}^{\ast }_2,\mathbf{x}_2\rangle _{(H^2_0)',H^2_0}=0,\quad \forall\, \mathbf{x}=(\mathbf{x}_1,\mathbf{x}_2)\in Ker(\mathcal{L}^{\prime })\}
\notag \\
& =\{\mathbf{y}^{\ast }\in \VVV_0^{\prime }\times (H_0^2)':\ \langle \mathbf{y}^{\ast }_1,\mathbb{A}^{-1}\mathbf{z}_1\rangle _{\VVV_0^{\prime },%
\VVV_0}+\langle \mathbf{y}^{\ast }_2,\mathbb{B}^{-1}\mathbf{z}_2\rangle _{(H^2_0)',H^2_0}=0,\quad \forall\, \mathbf{z}=(\mathbf{z}_1,\mathbf{z}_2)\in Ker(\mathcal{L}^{\ast })\}  \notag \\
& =\{\mathbf{y}^{\ast }\in \VVV_0^{\prime }\times (H^2_0)':\ \langle \mathbf{y}_1%
^{\ast },\mathbf{q}\rangle _{\VVV_0^{\prime },\VVV_0}=0,\quad
\forall\, \mathbf{q}\in (Ker(\mathcal{L}))_1\}=(Ker(\mathcal{L}))^{\perp },
\label{vprime}
\end{align}
where $\perp $ is intended to be the annihilator of the set and $(S)_i$ represents the $i$-th component of the elements of the set $S\subset \VVV_0\times H^2_0$. The last isomorphism is due to the fact that it holds
\begin{equation*}
\mathbb{Q}^{\ast }\mathbf{w}=\mathbb{A}(\mathbb{Q}(\mathbb{A}^{-1}\mathbf{w}%
)),
\end{equation*}
and thus (see \cite{GGPS})
\begin{equation*}
Ker(\mathcal{L}^{\ast })=\begin{bmatrix}
\mathbb{A}&0\\
\mathbf{0}&\mathbb{B}
    \end{bmatrix}Ker(\mathcal{L})=\mathbb A(Ker(\LL))_1	\times\{0\}.
\end{equation*}%
Let us denote by $\widetilde{\mathbb{Q}}\in \mathcal{B}(%
\HHH_0)$  the restriction to $\HHH_0$ of $\mathbb{Q}$ and note that, for any $z\in
\mathfrak{D}(\widetilde{\mathbb{A}})\times \mathfrak{D}(\widetilde{\mathbb{B}})\hookrightarrow \HHH_0\times L^2_0$, it holds
\begin{equation*}
\mathcal{L}_{|\mathfrak{D}(\widetilde{\mathbb{A}})\times \mathfrak{D}(\widetilde{\mathbb{B}})}(\mathbf{z},v)=(\widetilde{%
\mathbb{A}}\mathbf{z}+\widetilde{\mathbb{Q}}\mathbf{z},\widetilde{%
\mathbb{B}}v),
\end{equation*}%
and both $\widetilde{\mathbb{A}}^{-1}\widetilde{\mathbb{Q}}$ and $%
\widetilde{\mathbb{Q}}\widetilde{\mathbb{A}}^{-1}$ are compact operators on $%
\HHH_0$ since they are compositions of a compact and a bounded operator on $\HHH%
_{0}$. We can thus apply again \cite[Theorem 8.2]{GGPS} with $T=\begin{bmatrix}
      \widetilde{\mathbb{A}}^{-1}& 0\\
      \mathbf{0}&\widetilde{\mathbb{B}}^{-1}
\end{bmatrix}\in \mathcal{B}(\HHH_0\times L^2_0)$ and deduce that also $%
\mathcal{L}_{|{\mathfrak{D}(\widetilde{\mathbb{A}})\times \mathfrak{D}(\widetilde{\mathbb{B}})}}$ is a Fredholm
operator. Namely, since clearly $\mathcal{L}_{|{\mathfrak{D}(%
\widetilde{\mathbb{A}})}\times \mathfrak{D}(\widetilde{\mathbb{B}})}$ is selfadjoint, we immediately deduce that
\begin{equation}
Range(\mathcal{L}_{|{\mathfrak{D}(\widetilde{\mathbb{A}})\times \mathfrak{D}(\widetilde{\mathbb{B}})}})=Ker(\mathcal{L}%
_{|{\mathfrak{D}(\widetilde{\mathbb{A}})\times \mathfrak{D}(\widetilde{\mathbb{B}})}})^{\perp _{{\HHH_0}\times L^2_0}}=Ker(%
\mathcal{L})^{\perp _{{\HHH_0}\times L^2_0}}=Ker(\mathbb{A}+\mathbb{Q})^{\perp _{\HHH_0}}\times L^2_0,  \label{range}
\end{equation}%
where the last identity is due to the fact that $Ker(\mathcal{L})\subset
\mathfrak{D}(\widetilde{\mathbb{A}})\times\{0\}$. We can now
apply \cite[Corollary 3.11]{Chill1}. To this aim, let us fix $(\vphi%
_{\infty },u_\infty)\in \omega (\vphi_0,{u}_{0})$ and set $\vphi_{\infty }^{0}=%
\vphi_{\infty }-\al\in \mathfrak{D}(%
\widetilde{\mathbb{A}})$. With the notation of \cite[Corollary 3.11]%
{Chill1}, we set $V_{0}:=Ker(\mathcal{L}(\vphi_{\infty }^{0}, u_\infty))=Ker(\mathbb{A}+\mathbb{Q})\times\{0\}:=\widetilde{V}_0\times\{0\}$ and
define the projection $P_1\in \mathcal{B}(\VVV_0)$ as to be the
orthogonal $\VVV_0$-projection on $\tV_{0}$. Set also $\tV_{1}=Ker(P_1)=\tV_0^{\perp_{\textbf{V}_0}}$,
so that we have the direct (orthogonal in $\textbf{V}_0$) sum  $\VVV_0=\tV_{0}\oplus\tV_{1}$. By setting $V_1=\tV_1\times H^2_0$ and $P=(P_1,P_2)$, with $P_2\equiv 0$  (which is an orthogonal $\VVV_0\times H^2_0$-projector), we have, clearly, $\VVV_0\times H^2_0=V_0 \oplus V_{1}$. In this way, \cite[%
Hypothesis 3.2]{Chill1} is verified. Let us now verify \cite[Hypothesis 3.4]%
{Chill1}. We set $W:=(\VVV_0\times H^2_0)'$, $U=\VVV_0\times H^2_0$
and notice that the adjoint of $P$, $P^{\prime }:(\VVV_0\times H^2_0)'\rightarrow (\VVV_0\times H^2_0)'$, is such that
\begin{equation*}
Range(P^{\prime })=V_{1}^{\perp }=V_{0}^{\prime },\quad Ker(P^{\prime
})=V_{0}^{\perp }=V_{1}^{\prime },\quad (\VVV_0\times H^2_0)'=V_{0}^{\prime }\oplus V_{1}^{\prime }.
\end{equation*}%
Therefore, we have:

\begin{itemize}
\item[(I)] $W=(\VVV_0\times H^2_0)^{\prime}\hookrightarrow V'$ by
construction;

\item[(II)] $P^{\prime}W=P^{\prime}(\VVV_0\times H^2_0)'=V_0^{\prime}\subset
W$;

\item[(III)]  $\widetilde{\mathcal{E}}^{\prime}\in C^1(U,\VVV_0^{\prime }\times (H^2_0)')$ since $\widetilde{\Psi} \in C^3(\mathbb{R})$, by recalling once more the isomorphism between $\VVV_0^{\prime }\times (H^2_0)'$ and $(\VVV_0\times H^2_0)'$;

\item[(IV)] by \eqref{vprime}, $Range(\mathcal{L}(\vphi%
_\infty^0,u_\infty))=V_0^\perp=V_1^{\prime}\cap (\VVV_0\times H^2_0)'=V_1^{\prime}%
\cap W$.
\end{itemize}
We are only left to verify the assumptions of \cite[Corollary 3.11]{Chill1}.
Set $X:=\mathfrak{D}(\widetilde{\mathbb{A}})\times \mathfrak{D}(\widetilde{\mathbb{B}})$ and $Y:=\HHH_0\times L^2_0$ (
 that the space $Y$ is endowed with the inner product $((\mathbf{u},\tilde{u}),(\mathbf{v},\tilde{v}))_Y=(\mathbf{u},\mathbf{v})_{\mathbf{H}_0}+(\tilde{u},\tilde{v})_{L^2_0}$).

\begin{itemize}
\item[(1)] Since $V_0=Ker(\mathcal{L}(\vphi_\infty,{u}_\infty^0))\subset
\mathfrak{D}(\widetilde{\mathbb{A}})\times \mathfrak{D}(\widetilde{\mathbb{B}})$, $PX\subset X$. Moreover, $P^{\prime}Y=P^{\prime}(\HHH_0\times L^2_0)\subset
Y$. Indeed, $(\mathbf{h},w)\in \HHH_0\times L^2_0\hookrightarrow\mathbf{%
V}_0^{\prime}\times(H_0^2)'\cong(\mathbf{%
V}_0\times H_0^2)' $ can be written as the sum of
$$
(\mathbf{h}_1,0)=(P_{\HHH_0} \mathbf{h},0)\in Ker(\mathcal{L}(\vphi_\infty^0,u_\infty)),
$$
where $P_{\HHH%
_0}\in \mathcal{B}(\HHH_0)$ is the $\HHH_0$-orthogonal
projection onto $\tV_0$, and
$$
(\mathbf{h}_2,w)=((I-P_{%
\HHH_0})\mathbf{h},w)\in Range(\mathcal{L}(\mathbf{u}_\infty^0)_{|%
\mathfrak{D}(\widetilde{\mathbb{A}})})$$
by \eqref{range}, so that $(\mathbf{h%
}_1,\mathbf{h}_2)_{\HHH_0}=0$ (clearly \eqref{range} is necessary
to reach this conclusion, since $P_1$ is  a projection with respect to the inner product in $\VVV_0$ and not in $%
\HHH_0$). Then we have
\begin{align*}
&\langle P^{\prime}(\mathbf{h},w),(\mathbf{v},z)\rangle_{(\VVV_0\times H^2_0)',%
\VVV_0\times H^2_0}=\langle \mathbf{h},P_1\mathbf{v}\rangle_{\VVV_0^{\prime},%
\VVV_0}=(\mathbf{h},P_1\mathbf{v})_{\HHH_0} \\
&=(\mathbf{h}_1,P_1%
\mathbf{v})_{\HHH_0}=\langle P_1^{\prime}\mathbf{h}_1,\mathbf{v}\rangle_{%
\VVV_0^{\prime},\VVV_0}=(\mathbf{h}_1,\mathbf{v})_{\textbf{H}_0}=\langle \mathbf{h}_1,\mathbf{v}\rangle_{\VVV_0^{\prime},%
\VVV_0}=\langle (\mathbf{h}_1,0),(\mathbf{v},z)\rangle_{(\VVV_0\times H^2_0)',%
\VVV_0\times H^2_0},
\end{align*}
since $(\mathbf{h}_2,w)$ is orthogonal to $Ker(\mathcal{L})$ with respect to the
inner product in $\HHH_0\times L^2_0$ and due to the fact that $\mathbf{h}_1\in
\tV_1^\perp=Range(P_1^{\prime})$ and thus $P_1^{\prime}%
\mathbf{h}_1=\mathbf{h}_1$. Then, $P^{\prime}(\mathbf{h},w)=%
(\mathbf{h}_1,0)\in \HHH_0\times L^2_0$, i.e., $P^{\prime}Y\subset Y$. To conclude, $X$
and $Y$ are invariant under the action of $P$ and $P^{\prime}$, respectively.

\item[(2)]  For any $(\vphi_\infty,{u}_\infty)\in \omega(\vphi_0,u_0)$, by
Lemma \ref{separation}, there exists a set $\tilde{U}\subset I_\xi$ (%
$I_\xi$ defined in Section \ref{convec}) such that $\vphi%
_\infty(x)\in \tilde{U}$ for any $x\in \overline{\Omega}$ and $\psi$ is
analytic in $\tilde{U}$ (recall that $\widetilde{\psi}_{|I_\xi}=\psi$%
). Since all the other terms are linear in $X$ with values in $Y$ (possibly after an integration by parts), it thus holds that for each $(\vphi_\infty^0,u_0)\in  \omega (\vphi_0,u_\infty)-%
(\al,0)$ the function $\widetilde{%
\mathcal{E}}^{\prime}$ is real analytic in a neighborhood of $(\vphi
_\infty^0,u_\infty)$ in $X$ (see, e.g., \cite[Proof of Corollary 4.6]{Chill1}) with
values in $Y=\HHH_0\times L^2_0$.

\item[(3)]  $Ker(\mathcal{L}(\vphi_\infty^0,u_\infty))\subset
X$ is finite-dimensional.

\item[(4)] Recalling $Ker P^{\prime}=V_1^{\prime}=V_0^\perp$, let us assume $%
(\mathbf{w},\tilde{w})\in V_0^\perp\cap \HHH_0\times L^2_0$. Then, for any $(\mathbf{z},\tilde{z})\in V_0=Ker(%
\mathcal{L}(\vphi_\infty^0,u_\infty))$, since clearly $\tilde z=0$, it holds $(\tilde{w},\tilde z)_H=0$, and
\begin{equation*}
0=\langle\mathbf{w},\mathbf{z}\rangle_{\VVV_0^{\prime},\VVV_0}+\langle\tilde{w},\tilde{z}\rangle_{(H^2_0)^{\prime},H^2_0}=(%
\mathbf{w},\mathbf{z})_{\mathbf{H}_0},
\end{equation*}
i.e., $(\mathbf{w},\tilde{w})\in V_0^{\perp_{Y}}$. On the other hand, by the same
argument, $V_0^{\perp_{Y}}\subset {V}_0^\perp\cap Y$.
Thus we deduce $V_0^{\perp_{Y}}= {V}_0^\perp\cap Y$. Then, from \eqref{range}, we immediately infer
\begin{equation*}
Range(\mathcal{L}(\vphi_\infty^0,u_\infty)_{|\mathfrak{D}(\widetilde{\mathbb{A}}%
)\times \mathfrak{D}(\widetilde{\mathbb{B}}%
)})=V_0^{\perp_{Y}}={V}_0^\perp\cap Y=Ker
P^{\prime}\cap Y,
\end{equation*}
as desired.
\end{itemize}
Therefore, all the assumptions of \cite[Corollary 3.11]{Chill1} are
satisfied and the proof is finished.

\section*{Acknowledgements}
The research of D.  Caetano was
funded by an EPSRC grant under the MASDOC centre for doctoral training at the University of
Warwick. M. Grasselli and A. Poiatti have been partially funded by MIUR-PRIN Grant 2020F3NCPX "Mathematics for Industry 4.0 (Math4I4)". M. Grasselli and A. Poiatti are also members of Gruppo Nazionale per l’Analisi Matematica, la Probabilità e le loro Applicazioni (GNAMPA), Istituto Nazionale di Alta Matematica (INdAM). The research of M. Grasselli is part of the activities of ``Dipartimento di Eccellenza 2023-2027''.
This research was funded in part by the Austrian Science Fund (FWF) \href{https://doi.org/10.55776/ESP552}{10.55776/ESP552}. For open access purposes, the authors have applied a CC BY public copyright license to any author accepted manuscript version arising from this submission.


\bigskip
\begin{thebibliography}{}

\bibitem{AGP2024} H. Abels, H. Garcke, A. Poiatti, \textit{Mathematical analysis of a diffuse interface model for multi-phase flows of incompressible viscous fluids with different densities}, J. Math. Fluid Mech. \textbf{26} (2024), no. 29, 51 pp.

\bibitem{AGP1} H. Abels, H. Garcke, A. Poiatti, \textit{Diffuse interface model for two-phase flows on evolving eurfaces with different densities: local well-posedness}, arXiv:2407.14941.

\bibitem{AGP2} H. Abels, H. Garcke, A. Poiatti, \textit{Diffuse interface model for two-phase flows on evolving surfaces with different densities: global well-posedness}, Calc. Var. Partial Differential Equations \textbf{64} (2025), 141 .

\bibitem{ArrWalTorr18} M. Arroyo, N. Walani, A. Torres-Sanchez, D. Kaurin, D., \textit{Onsager’s variational principle in soft matter: introduction and application to the dynamics of adsorption of proteins onto fluid membranes}, in: The Role of Mechanics in the Study of Lipid Bilayers, D. Steigmann, Ed., Springer, Cham, 2018, pp. 287–332.

\bibitem{BHW2003} T. Baumgart, S.T. Hess, W.W. Webb, \textit{Imaging coexisting fluid domains in biomembrane models coupling curvature and line tension}, Nat. \textbf{425} (2003), 821-824.


\bibitem{CE}D. Caetano, C.M. Elliott, \textit{Cahn-Hilliard equations on an evolving surface}, European J. Appl. Math. \textbf{32} (2021), 937-1000.




\bibitem{CEGP}D. Caetano, C.M. Elliott, M. Grasselli, A. Poiatti, \textit{regularization and
	separation for evolving surface Cahn-Hilliard equations}, SIAM J. Math. Anal. \textbf{55} (2023), 6625-6675.

\bibitem{Chill1} R. Chill, \textit{On the {\L}ojasiewicz-Simon gradient
inequality}, J. Funct. Anal. \textbf{201} (2003), 572-601.

\bibitem{doi11} M. Doi, \textit{Onsager’s variational principle in soft matter},
J. Condens. Matter Phys. \textbf{23} (2011), 284118, 8 pp.

\bibitem{EllFriHob17-a} C.M. Elliott, H. Fritz, G. Hobbs, \textit{Small deformations of Helfrich energy minimising surfaces with applications to
biomembranes}, Math. Models Methods Appl. Sci. \textbf{27} (2017), 1547-1586.

\bibitem{EllGar97-b} C.M. Elliott, H.Garcke \textit{Diffusional Phase Transitions in Multicomponent Systems with a Concentration Dependent Mobility Matrix}, Physica D \textbf{109} (1997), 242-256

\bibitem{EllHat21} C.M. Elliott, L. Hatcher, \textit{Domain formation via phase separation for spherical biomembranes with small deformations},
European J. Appl. Math. \textbf{32} (2021), 1127-1152.
	
\bibitem{EllLuc91} C.M. Elliott, S. Luckhaus, \textit{A generalized diffusion
		equation for phase separation of a multi component mixture with interfacial
		free energy}, IMA Preprint Series \# 887, 1991.
        
\bibitem{ES} C.M. Elliott, T. Sales, \textit{{N}avier-{S}tokes-{C}ahn-{H}illiard equations on evolving surfaces}, Interfaces Free Bound. \textbf{27} (2025), 285-348.


\bibitem{GGGP2023} C.G. Gal, A. Giorgini, M. Grasselli, A. Poiatti, \textit{Global well-posedness and convergence to equilibrium for the Abels-Garcke-Gr\"{u}n model with nonlocal free energy}, J. Math. Pures Appl. (9) \textbf{178} (2023), 46-109.

\bibitem{GGPS} C.G. Gal, M. Grasselli, A. Poiatti,  J.L. Shomberg, \textit{Multi-component Cahn-Hilliard systems with singular potentials: Theoretical results},  Appl. Math. Optim.  \textbf{88} (2023), Paper No. 73, 46 pp.


\bibitem{GalPoia} C.G. Gal, A. Poiatti, \textit{Unified framework for the separation property in binary phase segregation processes with singular entropy densities}, European J. Appl. Math. \textbf{36} (2025),40-67.
	
\bibitem{Gar05} H. Garcke, \textit{On a Cahn-Hilliard model for phase
		separation with elastic misfit}, Ann. Inst. H. Poincar\'{e} Anal. Non Lin\'{e}aire \textbf{22} (2005), 165-185.

\bibitem{AC2023} M. Grasselli, A. Poiatti, \textit{Multi-component conserved Allen-Cahn equations,} Interfaces Free Bound. \textbf{26} (2024), 489-541.

\bibitem{LukeThesis} L. Hatcher, \textit{Phase field models for small deformations of biomembranes arising as Helfrich
energy equilibria}, September 2020. URL http://wrap.warwick.ac.uk/152741/

\bibitem{HF2011} F.A. Heberle, G.W. Feigenson, \textit{Phase Separation in Lipid Membranes}, Cold Spring Harb. Perspect Biol. 2011,
3:a004630, 14 pp.

\bibitem{Lady} O.A. Lady\v{z}enskaja, V.A. Solonnikov, N. Ural'ceva, \textit{ Linear and quasilinear equations of parabolic type},
Translations of Mathematical Monographs \textbf{23}, American Mathematical Society, Providence, R.I. 1968.

\bibitem{LHXF2024} L. Liu, S. Huang, X. Xiao, X. Fenga, \textit{Mathematical modeling and numerical simulation of the N-component
Cahn-Hilliard model on evolving surfaces}, J. Comp. Phys. \textbf{513} (2024), 113189.

\bibitem{ons31a} L.  Onsager, \textit{Reciprocal relations in irreversible processes. I}, Phys. Rev.
\textbf{37} (1931), 405–426.

\bibitem{P} A. Poiatti, \textit{The 3D strict separation property for the nonlocal Cahn--Hilliard
equation with singular potential}, Anal. PDE \textbf{18} (2025), 109-139.

\bibitem{RB2024} U. Rana, K. Xu, A. Narayanan, M.T. Walls, A.Z. Panagiotopoulos, J.L. Avalos, C.P. Brangwynne, \textit{Asymmetric oligomerization state and sequence patterning can tune multiphase condensate miscibility}, Nat. Chem. (2024), https://doi.org/10.1038/s41557-024-01456-6, 24 pp.

\bibitem{RS} E. Rocca, G. Schimperna, \textit{Universal attractor for some
	singular phase transition systems}, Phys. D \textbf{192} (2004), 279-307.

\bibitem{SB2021} K. Shrinivas, M.P. Brenner, \textit{Phase separation in fluids with many interacting components}, Proc. Natl. Acad. Sci. USA
\textbf{118} (2021), e2108551118, 8 pp.

\bibitem{Temam} R. Temam, \textit{Infinite-Dimensional Dynamical Systems in
	Mechanics and Physics}, Springer-Verlag, New York, 1997.

\bibitem{WBV2012} T. Witkowski, R. Backofen, A. Voigt, \textit{The influence of membrane bound proteins on phase separation and
coarsening in cell membranes}, Phys. Chem. Chem. Phys. \textbf{14} (2012), 14509-14515.

\bibitem{ZH1989} O.-Y. Zhong-can, W. Helfrich,  \textit{Bending energy of vesicle membranes: General expressions for the first, second, and third
variation of the shape energy and applications to spheres and cylinders}, Phys. Rev. A \textbf{39} (1989), 5280-5288.


\bibitem{ZL2022} D. Zwicker, L. Laan, \textit{Evolved interactions stabilize many coexisting phases in multicomponent liquids}, Proc. Natl. Acad. Sci. USA \textbf{119} (2022), e2201250119, 8 pp. Correction in: Proc. Natl. Acad. Sci. USA. \textbf{119} (2022), e2219263119, 2 pp.




\end{thebibliography}

\end{document}